%% file: sln.tex
\documentclass [11pt] {amsart} 
\usepackage {amsmath, amssymb, amscd, mathrsfs, graphicx, color, url}

\newtheorem {theorem}{Theorem}[section]
\newtheorem {definition}[theorem]{Definition}
\newtheorem {lemma}[theorem]{Lemma}
\newtheorem {proposition}[theorem]{Proposition}

\newtheorem {conjecture}[theorem]{Conjecture}

\newtheorem {remark}[theorem]{Remark}
\newtheorem {example}[theorem]{Example}

\def\zz {{\mathbb{Z}}}
\def\rr {{\mathbb{R}}}
\def\cc {{\mathbb{C}}}
\def\qq {{\mathbb{Q}}}
\def\pp {{\mathbb{P}}}

\def\bb {{\mathcal{B}}}

\def\oo {{\mathcal{O}}}
\def\ss {{\mathcal{S}}}

\def\lll {{\mathcal{L}}}

\def\ss {{\mathcal{S}}}
\def\ccc {{\mathcal{C}}}
\def\yy {{\mathcal{Y}}}
\def\nnn {{\mathcal{N}}}
\def\ua {{\mathcal{U}}}

\def\hh {{\mathfrak{h}}}
\def\gg {{\mathfrak{g}}}
\def\sl {{\mathfrak{sl}}}
\def\so {{\mathfrak{so}}}

\def\gl {{\mathfrak{gl}}}
\def\lfrak {{\mathfrak{L}}}

\def\ggg {{\mathfrak{g}}}
\def\hhh {{\mathfrak{h}}}
\def\ppp {{\mathfrak{p}}}

\def\cpn {{\cc \pp^{n-1}}}

\def\lk {{\kappa}}
\def\nn {{(n)}}
\def\ymnt {{\yy_{m, n, \tau}}}
\def\lambdav {{\vec{\lambda}}}
\def\muv {{\vec{\mu}}}

\def\ff {{\mathscr{F}}}
\def\ee {{\mathscr{E}}}
\def\kr {{\mathscr{H}}}
\def\bb {{\mathscr{B}}}
\def\xx {{\mathscr{X}}}
\def\tt {{\mathscr{T}}}

\def\Reg {{\operatorname {Reg} \ }}
\def\symp {{\scriptscriptstyle{\mathrm{symp}}}}
\def\mapp {{q}}
\def\krs {{\kr_{\nn \symp}}}
\def\krss {{\kr_{\nn \symp}^*}}
\def\krsk {{\kr_{\nn \symp}^k}}
\def\ggp {{\hat \ggg^{\pi}}}
\def\gp {{\tilde \ggg^{\pi}}}
\def\gpi {{\ggg^{\pi}}}
\def\sln {{\sl(n)}}
\def\slan {{\sl(n, \cc)^{(1(n-1))}}}
\def\chin {{\chi^{(1(n-1))}}}

\def\fin {{\hfill \square}}

\def\im {{\hspace{2pt} \text{Im }}}
\def\ker {{\hspace{2pt} \text{Ker }}}

\def\mat {{\mathfrak{m}}}
\def\lfrak {{\mathcal{L}}}
\def\pim {{\operatorname{proj}^{\scriptscriptstyle{\mathrm{im}}}}}

\def\jm {{\scriptscriptstyle{\mathrm{JM}}}}
\def\sjm {{\ss^{\jm}}}

\def\reg {{\scriptscriptstyle{\mathrm{reg}}}}
\def\red {{\scriptscriptstyle{\mathrm{red}}}}
\def\sub {{\scriptscriptstyle{\mathrm{sub}}}}
\def\resc {{\scriptscriptstyle{\mathrm{resc}}}}
\def\sym {{\operatorname{Sym}}}
\def\inv {{\operatorname{Inv}}}

\def\hom {{\operatorname{Hom}}}
  
\def\Id {{I}}
\def\iso {{\ \cong \ }}

\def\conf {{Conf}}
\def\confb {{BConf}}

\begin{document}

\title [Link homology theories from symplectic geometry] 
{Link homology theories from symplectic geometry} \author 
[Ciprian Manolescu]{Ciprian Manolescu} 
\thanks {The author was supported by a Clay Research Fellowship.}

\begin {abstract} For each positive integer $n,$ Khovanov and Rozansky constructed an 
invariant of links in the form of a doubly-graded cohomology theory whose Euler 
characteristic is the $\sln$ link polynomial. We use Lagrangian Floer cohomology on 
some suitable affine varieties to build a similar series of link invariants, and we 
conjecture them to be equal to those of Khovanov and Rozansky after a 
collapse of 
the bigrading. Our work is a generalization of that of Seidel and Smith, who treated 
the case $n=2.$ \end {abstract}

\address {Department of Mathematics, Columbia University\\ 2990 Broadway, New York, 
NY 10027}
\email {cm@math.columbia.edu}

\maketitle

\section {Introduction}

For any $n > 0,$ the quantum $\sln$ polynomial invariant $P_{\nn}$ of an oriented link 
$\lk \subset S^3$ is uniquely determined by the skein relation: 
\begin{equation}
\label {skey}
q^n P_{\nn} \left( \, \includegraphics[scale=0.3]{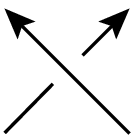} \, \right)
- q^{-n} P_{\nn}  \left( \, \includegraphics[scale=0.3]{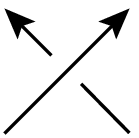} \,  \right)
= ( q - q^{-1}) P_{\nn} \left(\includegraphics[scale=0.3]{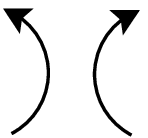} \right),
\end{equation}
together with the value on the unknot 
$P_{\nn}(\text{unknot}) = (q^n - q^{-n})/(q-q^{-1}).$ The invariant can also be defined 
in terms of the representation theory of the quantum group $U_q(\sln),$ hence the 
name.  When $n=2$ we obtain the Jones polynomial, up to a $q+q^{-1}$ factor. The same skein 
relation for $n=0$ with the normalization $P_{(0)}(\text{unknot})=1$ gives the 
Alexander polynomial, but the representation theoretic story is 
somewhat different in this case. The polynomials $P_{\nn}$ are all different 
specializations of a single link invariant, the two variable HOMFLY polynomial 
\cite{HOMFLY}, \cite{PT}.

Khovanov and Rozansky \cite{KR1} associated to every link $\lk$ a series of 
bigraded cohomology theories $\kr^{i, j}_{\nn}(\lk)$ for $n > 0$ and showed 
that they are link invariants. Their theories can be interpreted as 
categorifications of $P_{\nn},$ in the sense that $$ P_{\nn}(\lk) = \sum_{i, 
j\in \zz} (-1)^i q^j \dim_{\qq} \kr^{i, j}_{\nn} (\lk).$$ When $n=2,$ they 
recover the older categorification of the Jones polynomial due to Khovanov 
\cite{Kh1}.

Khovanov-Rozansky homology is particularly interesting because it is conjectured to 
be related to the knot Floer homology of Ozsv\'ath-Szab\'o and Rasmussen \cite{OS}, 
\cite{R}. Knot Floer homology is an invariant defined used Lagrangian Floer homology, 
and an important question is to find a way to compute it algorithmically. Its graded 
Euler characteristic is the Alexander polynomial corresponding to $n=0$ above. 
On the other hand, Khovanov-Rozansky homology is algorithmically computable by definition, 
and the hope is to be able to extract the case $n=0$ from the $n > 0$ theories. A 
precise conjecture in this direction was made by Dunfield, Gukov and Rasmussen in 
\cite{DGR}, and a potentially useful triply graded categorification of the HOMFLY 
polynomial was constructed by Khovanov and Rozansky in the sequel \cite{KR2}.

In this paper we construct a sequence of link invariants $\krss(\lk)$ using 
Lagrangian Floer theory. This has been done by Seidel and Smith in the case $n=2$ 
\cite{SS}, and our work is a generalization of theirs. We conjecture that our 
invariants are related to Khovanov-Rozansky homology of the mirror link $\lk^!$ in 
the following way:

\begin {conjecture}
\label {conj}
$$\krsk(\lk) \otimes \qq = \bigoplus_{i + j =k} \kr^{i, j}_{\nn}(\lk^!).$$
\end {conjecture}

Our construction is inspired from that of Seidel and Smith, with some differences coming from 
the fact that the standard (quantum) representation $V$ of $\sln$ is not self-dual for $n > 
2.$ As in \cite{SS}, we start by presenting the link $\lk$ as the closure of an $m$-stranded 
braid $\beta \in Br_m.$ The rough idea is to find a symplectic manifold $(M, \omega)$ with an 
action of the braid group by symplectomorphisms $\phi: Br_m \to \pi_0(Symp (M, \omega)),$ to 
take a specific Lagrangian $L \subset M$ and to consider the Floer cohomology of $L$ and 
$\phi(\beta)L$ in $M.$ Following the ideas of Khovanov from \cite{Kh2}, we would like the 
Grothendieck group of the derived Fukaya category of $M$ to be related to the space $\inv 
(m,n)$ of invariants in the representation $V^{\otimes m} \otimes (V^{\otimes m})^*$ of 
$\sln.$ The reason is that $\inv (m,n)$ naturally occurs in the representation theoretic 
definition of the polynomial $P_{\nn}.$ On the other hand, the Grothendieck group of the 
derived Fukaya category is not a well-understood object in general, but it is related (and in 
some special cases equal) to the middle dimensional homology of $M.$ In the case $n=2$ for 
example, Seidel and Smith worked with a symplectic manifold whose middle dimensional Betti 
number is the $m$-th Catalan number, the same as the dimension of $\inv (m,2).$ Their 
construction uses the geometry of the adjoint quotient and nilpotent slices in the Lie algebra 
$\sl(2m),$ and our manifolds $M = \ymnt$ below are a natural extension of theirs, obtained by 
looking at the Lie algebra $\sl(mn)$ instead.

We define the bipartite configuration space $$ \confb^0_m = \Bigl\{ \bigl( \{ \lambda_1, 
\dots, \lambda_m \}, \{ \mu_1, \dots, \mu_m \} \bigr ) \ | \ \lambda_i, \mu_j \in \cc 
\text{ distinct}, \ \sum \lambda_i + (n-1)\sum \mu_j =0  \Bigr \}.$$

The elements in each of the two sets $\lambdav = \{ \lambda_1, \dots, \lambda_m \}$ and $\muv 
= \{ \mu_1, \dots, \mu_m \}$ are not ordered, but the pair $\tau = (\lambdav, \muv)$ is 
ordered. The fundamental group of $\confb_m^0$ is the colored braid group on two colors 
$Br_{m,m}.$ This has a (noncanonical) subgroup isomorphic to $Br_m,$ which corresponds to 
keeping $\muv$ fixed.

For each $m, n > 0$ and $\tau = (\lambdav, \muv) \in \confb^0_m,$ we construct a complex 
affine 
variety $\ymnt$ as follows. Let $N_{m, n}$ be a nilpotent element in 
$\sl(mn)$ with $n$ Jordan blocks of size $m.$ After a change of basis, we can write 
$$ N_{m,n} =  \left( \begin {array}{ccccc} 0 & I &  &  & \\
 & 0 & I & \ & \\
 & \ & \ & \ldots & \ \\
 & \ & \  & \  &  I \\
 & \ & \ & \ & 0
\end {array} \right ), $$
where $I$ and $0$ are the $n$-by-$n$ identity and zero matrix, respectively.

The following is a transverse slice to the adjoint orbit of $N_{m, n}$ in $\sl(mn):$

$$ \ss_{m, n} = \left \{ \left( \begin {array}{ccccc} Y_1 & I &  &  & \\
Y_2 & 0 & I & \ & \\
\ldots & \ & \ & \ldots & \ \\
 & \ & \  & \  &  I \\
Y_m & \ & \ & \ & 0
\end {array} \right ): Y_i \in \gl(n), \ \operatorname{trace}(Y_1) =0  \right \}.$$

Consider also the diagonal matrix $D_{\tau} \in \sl(mn)$ with eigenvalues $\lambda_1, \dots, 
\lambda_m$ with multiplicity $1$ (we call these {\it thin} eigenvalues), and $\mu_1, \dots, 
\mu_m,$ each with multiplicity $n-1$ (we call these {\it thick} eigenvalues). Let 
$\oo_{\tau}$ be the adjoint orbit of $D_{\tau}$ in $\sl(mn),$ and set
\begin {equation}
\label {ot}
\ymnt = \ss_{m,n} \cap \oo_{\tau}.
\end {equation}

We will show that, as $\tau$ varies over $\confb^0_m,$ the spaces $\yy_{m, n, \tau}$ form a 
symplectic fibration that admits good parallel transport maps. Moreover, for a specific 
$\tau$ we build a Lagrangian $L \subset \ymnt$ by iterating a relative vanishing $\cpn$ 
construction. The vanishing projective spaces which we introduce in this paper 
share many properties with the usual vanishing cycles in symplectic geometry, and may be of 
interest on their own. Indeed, Huybrechts and Thomas \cite{HT} pointed out that Lagrangian 
projective spaces can play a role in homological mirror symmetry, akin to the role of 
Lagrangian spheres.

The local model for the vanishing $\cpn$ construction is the space $Z$ of $n$-by-$n$ 
traceless matrices which have an eigenvalue of multiplicity at least $n-1.$ Let us assume 
$n > 2.$ (When $n=2,$ our vanishing $\cpn$ is a usual vanishing cycle.) Consider the map 
$\chi: Z \to \cc$ which takes a matrix to the corresponding high multiplicity eigenvalue, 
and denote by $Z_t$ the fiber over $t \in \cc.$ Note that $Z$ has a singularity at the 
origin. Equip (the smooth strata of) $Z$ and the fibers $Z_t$ with the restriction of the 
standard K\"ahler form on the space of all $n$-by-$n$ matrices, viewed as $\cc^{n^2}.$ We 
can then consider parallel transport in $Z$ with respect to $\chi,$ as long as we stay away 
from $Z_0.$ Taking a linear path from $t \in \cc^*$ to the origin, we let $L_t$ be the set 
of points in $Z_t$ which are taken to $0 \in Z_0$ in the limit of parallel transport along 
that path. It is not hard to describe $L_t$ explicitly.

\begin {lemma}
\label {rex}
Consider the diagonal matrix $E_t=diag(t,t,\dots,t, (1-n)t) \in Z_t.$ Then:
$$ L_t = \{ UE_tU^{-1} \ | \ U \in U(n) \}. $$
\end {lemma}

Thus $L_t$ is diffeomorphic to $U(n)/(U(n-1) \times U(1)) \cong \cpn.$ It is also easy to 
check that it is a Lagrangian subspace of $Z_t.$ This is the basic model for a 
\textit {vanishing projective space}. 

In Section~\ref{sec:vvo} we establish a more general version of Lemma~\ref{rex}, in which 
the smooth stratum of $Z$ is equipped with any K\"ahler form satisfying certain real 
analyticity and proportionality conditions. We define $L_t$ just as above, and show that it 
is a Lagrangian $\cpn$ in $Z_t.$ This is done by observing that $Z$ is the GIT quotient of 
$\cc^{2n}$ by the linear $\cc^*$-action with weights $1$ and $-1,$ each with multiplicity 
$n.$ We can then lift the spaces $L_t$ to the affine space $\cc^{2n},$ where we argue that 
they are vanishing cycles in a certain singular metric. Thus, the vanishing projective 
spaces are quotients of vanishing cycles by a circle action.

The spaces $Z_t$ are the same as the fibers $\yy_{1,n, \tau}$ from (\ref{ot}), for $m=1.$ 
Therefore, for $m=1$ Lemma~\ref{rex} gives a Lagrangian $\cpn$ in each $\yy_{1,n, \tau}.$ 
We then proceed to construct a Lagrangian $L\subset \ymnt$ inductively in $m,$ by using a 
relative version of the vanishing $\cpn$'s. The resulting $L$ is diffeomorphic to the 
product of $m$ copies of $\cpn.$ Given an element $b \in Br_m \subset Br_{m,m}$ whose 
closure is a link $\lk,$ we use parallel transport along the corresponding loop $\beta$ in 
$ \confb^0_m$ to construct a second Lagrangian $L'=h^{\resc}_{\beta}(L)$ in the same 
$\ymnt.$ Let $w$ be the writhe of the braid $b.$ Our main result is:

\begin {theorem}
\label {the}
Up to isomorphism of graded abelian groups, the shifted Floer cohomology groups
$$ \krss(\lk) = HF^{*+(n-1)(m+w)}(L, L')$$
depend only on the oriented link $\lk.$ 
\end {theorem}

The proof of the theorem involves checking invariance under the Markov moves $I$ and $II$ 
which relate braids with the same closure. 

We managed to compute the groups $\krs$ in a few examples. For the unknot we have 
$\krss(\text{unknot}) = H^{*+n-1}(\cpn),$ while for the unlink of $p$ components we get 
the tensor product of $p$ copies of the same group. The first nontrivial computation is for the 
trefoil, for which we have the following result, consistent with the formula in 
\cite[Proposition 6.6]{DGR}.

\begin {proposition}
\label {3f}
When $\lk$ is the left-handed trefoil, the groups $\krs$ are given by:
$$ \krss(\lk) = H^{*-n+1}(\cpn) \oplus H^{*-n-1}(UT\cpn),$$
where $UT\cpn$ is the unit tangent bundle to $\cpn.$
\end {proposition}

One advantage that our theory has over that of Khovanov and Rozansky is that it 
produces abelian groups rather than vector spaces over $\qq$. For example, there is a 
$\zz/n\zz$ torsion group appearing in the computation of $\krs$ of the trefoil, and 
that group is invisible in $\kr_{\nn}.$

Nevertheless, there is also an obvious shortcoming of our theory, the fact that it does 
not come with a bigrading. In the case $n=2,$ a bigrading for the Seidel-Smith cochain 
complex was constructed by the author in \cite{M}, but it is not yet clear whether it 
descends to cohomology. The bigrading was built using an open holomorphic embedding of the 
manifold $\yy_{m,2,\tau}$ into a Hilbert scheme. It would be interesting to study whether 
similar embeddings exist for $n > 2.$

It is worth noting here that there are several alternate descriptions of the spaces 
$\ymnt,$ and these could lead to further insights into our construction. (This 
observation was made by Seidel and Smith in their introduction to \cite{SS}, for $n=2.$) 
First, the spaces $\ymnt$ can be viewed as quiver varieties of type $A_{2m-1},$ cf. 
Nakajima's Conjecture 8.6 in \cite{N1}, proved by Maffei in \cite{Maf}. Via an ADHM 
transform, these can also be viewed as moduli spaces of rank $n$ instantons on an ALE 
space (cf. \cite{KN}), or as moduli spaces of solutions to Nahm's equations (cf. 
\cite{AB}, \cite{Kr2}, \cite{N1}). To some extent these alternate descriptions are 
conjectural, because most of the works cited in this paragraph only deal with nilpotent 
orbits. However, we expect the respective results to generalize to our situation.

This paper is organized as follows. In Section~\ref{sec:gaos} we study the general 
properties of intersections between transverse slices and adjoint orbits. In 
Section~\ref{sec:av} we apply these properties to study our objects of interest, the spaces 
$\ymnt,$ as well as their degenerations. In Section~\ref{sec:vvo} we present the 
construction of vanishing projective spaces. In Section~\ref{sec:fiba2} we study in detail 
a geometric situation that is key to the proof of Markov $II$ invariance. In 
Section~\ref{sec:floer} we review the definition of Floer cohomology and discuss some 
relevant properties. Section~\ref{sec:maind} contains the construction of the Lagrangians 
$L, L' \subset \ymnt.$ We prove Theorem~\ref{the} in Section~\ref{sec:markov} by showing 
invariance under the two Markov moves, and then do the trefoil computation in 
Section~\ref{sec:trefoil}. In the last section we speculate on the existence of other 
classes of link invariants: some that correspond to various other Lie algebras and 
representations, and some that we expect to arise by considering a particular involution on 
the spaces $\ymnt.$ The latter invariants should form a series parametrized by integers $n 
\geq 2,$ such that the $n=2$ case gives the Heegaard Floer homology of the double branched 
cover.

We should note that many of the arguments in this paper follow the ones given by Seidel and 
Smith in \cite{SS} closely, and we sometimes refer to their article for full details. 
However, there are also several aspects that are fundamentally new in our work as compared 
with the $n=2$ case, and these are the ones which we choose to emphasize in our exposition. 
The first is the distinction between thin and thick eignevalues, stemming from the fact 
that $V \not \cong V^*$ for the standard representation of $\sl(n)$ when $n > 2.$ A 
consequence of the appearance of thick eigenvalues is that the orbit $\oo_{\tau}$ from 
(\ref{ot}) is no longer maximal, and we need a careful study of the geometry of 
intersections between slices and non-maximal orbits. One new phenomenon is that the total 
space consisting of all $\ymnt$'s and their degenerations is singular for $n > 2.$ It is in 
this space where the vanishing $\cpn$ construction is made, and as far as we know this 
construction is new. We deal with vanishing projective spaces by viewing them as the 
quotients of ordinary vanishing cycles by an $S^1$-action. The corresponding vanishing 
cycles appear in a singular metric though, and in order to justify their existence we have 
to resort to a real analyticity condition on the metric. Also, the parallel transport 
estimates for vanishing $\cpn$'s in Section~\ref{sec:fiba2}, although similar in spirit to 
those in \cite{SS}, are somewhat more involved here. In particular, some care is needed in 
making sure that the K\"ahler metrics on our local models in singular spaces are equivalent 
to standard metrics. Finally, the discussion of orientations in Floer cohomology depends on 
the parity of $n,$ because for $n$ odd the complex projective space $\cpn$ is not spin.

\medskip \noindent \textbf{Acknowledgements.} I am grateful to Mikhail Khovanov, Peter 
Kronheimer, Duong Phong and Jacob Rasmussen for several valuable discussions, to Edward 
Bierstone, Krzysztof Kurdyka and Alec Mihailovs for helpful email correspondence, and to Paul 
Seidel and Ivan Smith for their interest in this work.

\section {The geometry of adjoint orbits and slices}
\label {sec:gaos}

This section parallels Section 2 in \cite{SS}. We collect some facts about partial 
Grothendieck resolutions and their intersections with transverse slices. Our main 
reference is the work of Borho and MacPherson \cite{BM}, where partial Grothendieck 
resolutions are studied in detail. We also drew inspiration from \cite{Kr2}, \cite{N1}, 
\cite{Ro}, \cite{Sl} and \cite{SS}. Some of the results there were formulated only for 
the full Grothendieck resolution, but admit straightforward generalizations to the 
partial case.

The discussion below can be made more general, but for the sake of concreteness we 
restrict our attention to the case of the group $G= SL(N, \cc)$ and its Lie algebra 
$\ggg = \sl(N, \cc).$ We also fix the standard basis for $\cc^N.$ We denote by $\hh 
\cong \cc^{N-1}$ the corresponding Cartan subalgebra of traceless diagonal matrices, and 
by $W = S_{N}$ its Weyl group.

\subsection {Partial Grothendieck resolutions}
Let $\pi=(\pi_1, \dots, \pi_s)$ be a partition of $N = \pi_1 + \dots + \pi_s,$ with $\pi_1 
\geq \pi_2 \geq \dots \geq \pi_s > 0.$ We denote by $m_{\pi}(k)\geq 0$ the number of 
times $k$ appears among the $\pi_j,$ for $k=1, \dots, N.$ Thus $\sum m_{\pi}(k) = s$ 
and $\sum k m_{\pi}(k) = N.$ Sometimes the partition $\pi$ is also written as $\pi = 
(1^{m_{\pi}(1)} 2^{m_{\pi}(2)} \dots N^{m_{\pi}(N)}).$ There is a dual partition $\pi^* 
= (\pi^*_1, \pi^*_2, \dots)$ with $\pi^*_j = m_{\pi}(j) + m_{\pi}(j+1) + \dots + 
m_{\pi}(N).$

Associated to $\pi$ is a partial flag variety 
$$ \ff^{\pi} = \{0=F_0 \subset F_1 \subset \dots \subset F_s=\cc^N \ | \ \dim (F_j) - 
\dim (F_{j-1}) = \pi_j \}.$$

We denote by $F^{st} \in \ff^{\pi}$ the standard flag $0 \subset \cc^{\pi_1} \subset 
\cc^{\pi_1 + \pi_2} \subset \dots \subset \cc^N.$ Note that for every flag $F=(F_0, 
\dots, F_s) \in \ff^{\pi}$ there is a corresponding parabolic subalgebra $\ppp(F) 
\subset \ggg$ consisting of those matrices $x\in \ggg$ preserving the flag, i.e. such 
that $x(F_j) \subset F_j$ for all $j.$ Thus we can identify $\ff^{\pi}$ with the space 
of parabolic Lie subalgebras of $\ggg$ conjugate to $\ppp(F^{st}).$

The adjoint quotient map $\chi: \ggg \to \ggg/G = \hh/W$ takes a matrix to the set of its 
generalized eigenvalues. Note that $\hh/W$ can be identified with $\cc^{N-1}$ via 
symmetric polynomials. Set
$$\ggp = \{(x, F) \ | \ F \in \ff^{\pi}, x \in \ppp(F) \}.$$

Consider also the subgroup $W_{\pi} \subset W$ given by $S_{\pi_1} \times \dots 
\times S_{\pi_s} \subset S_N.$ This is the Weyl group corresponding to the Levi subalgebra of 
$\ppp(F^{st}).$ The {\it partial simultaneous resolution} of $\chi$ associated to 
the partition $\pi$ consists of the commutative diagram: 
\begin {equation}
\label {psr}
\begin {CD}  \ggp @>>> \ggg 
\\ @V{\hat \chi^{\pi}}VV @V{\chi}VV \\ \hh/ W_{\pi} @>>> \hh/W 
\end {CD}
\end {equation}
Observe that for every $(x, F) \in \ggp$ there is an induced action of $x$ on each
quotient $F_j/F_{j-1},$ for $j=1, \dots, s.$ This gives elements $x_j^F \in
\gl(F_j/F_{j-1}).$  The map 
$\hat \chi^{\pi}$ is defined to take a pair $(x, F)$ to the sets of generalized 
eigenvalues of $x_j^F$ for all $j.$

The map $\tilde \chi = \hat \chi^{\pi}$ for $\pi=(1^N)$ was the one 
originally studied by Grothendieck. In that case $W_{\pi} = 1,\ \tilde \ggg = \ggp$ is a 
smooth manifold and $\tilde \chi$ is a {\it simultaneous resolution} in the following 
sense: $\tilde \chi$ is a submersion with the property that for each $\tilde t \in \hh, 
\tilde \chi^{-1}(\tilde t)$ is a resolution of singularities for $\chi^{-1}(t),$ where 
$t$ is the image of $\tilde t$ in $\hh/W.$ Diagramatically,
\begin {equation}
\label {grothen}
\begin {CD}  \tilde \ggg @>>> \ggg
\\ @V{\tilde \chi}VV @V{\chi}VV \\ \hh @>>> \hh/W.
\end {CD}
\end {equation}

For general $\pi,$ the variety $\ggp$ is not smooth. $\hat \chi^{\pi}$ is called a 
partial resolution 
because the map $\tilde \ggg \to \ggg$ factors through $\ggp.$ However, as explained 
below, if we restrict $\ggp \to \ggg$ to a certain subset of $\ggp$ we do get an honest 
simultaneous resolution.

\subsection {Restricted partial Grothendieck resolutions} \label{sec:rpr} Set 
$$\gp = \{(x, F) \in \ggp \ | \ x_j^F \in Z(\gl(F_j/F_{j-1})) \text{ for all } j \}. $$

Consider the subspace $\hh^{\pi} \subset \hh/W_{\pi}$ made of traceless diagonal 
matrices of the type $\alpha = diag(\alpha_1, \dots, \alpha_1, \alpha_2, \dots, 
\alpha_2, \dots, \alpha_s),$ where each $\alpha_j$ appears exactly $\pi_j$ times. 
For every $(x, F) \in 
\gp$ we have $x_j^F = \alpha_j \cdot I$ for some $\alpha_j \in \cc.$ Therefore, there 
is a induced map $\tilde \chi^{\pi}: \gp \to \hh^{\pi}.$ This fits into a commutative diagram 
\begin {equation}
\label {rpr}
\begin {CD}
\gp @>>> \gpi \\ @V{\tilde \chi^{\pi}}VV @V{\chi^{\pi}}VV \\ \hh^{\pi} @>>> 
\hh^{\pi}/W^{\pi} \end {CD}
\end {equation}

We call this diagram the {\it restricted partial simultaneous resolution} associated to 
$\pi.$ The subset $\gpi \subset \ggg$ is defined to be the image of $\gp$ under the map 
$\ggp \to \ggg$ from (\ref{psr}), and the {\it partial Weyl group} $W^{\pi}$ is defined as 
$W^{\pi} = S_{m_{\pi}(1)} \times \dots \times S_{m_{\pi}(N)} \subset S_s.$ The vertical map 
$\chi^{\pi}$ is naturally induced from the adjoint quotient map $\chi.$ Note that 
$\hh^{\pi} \cong \cc^{s-1},$ while $\hh^{\pi}/W^{\pi}$ is naturally the quotient by $\cc$ 
of a product of symmetric spaces $\sym^{m_{\pi}(1)}(\cc) \times \dots \times 
\sym^{m_{\pi}(N)}(\cc).$ (This quotient can also be identified with $\cc^{s-1},$ using 
symmetric polynomials.)

\begin {example}
\label {ex0}
For $\pi = (N-1,1),$ the space $\gpi$ consists of traceless $N$-by-$N$ matrices having
an eigenspace of dimension at least $N-1.$ (This is the space $Z$ mentioned in the 
introduction.) The map $\chi^{\pi}$ takes a 
matrix to the corresponding high multiplicity eigenvalue in $\hh^{\pi}/W^{\pi} \cong \cc.$ 
\end {example}

In general, it is easy to see that $\gpi$ is a closed subvariety of $\ggg.$ Also, the 
adjoint action $Ad$ of $G$ on $\ggg$ induces actions of $G$ on all the spaces in 
(\ref{rpr}), and all the maps there are $G$-equivariant. It follows that $\gpi$ is a union 
of adjoint orbits in $\ggg.$

The proof of the following proposition is analogous to that of Lemma 4 in \cite{SS}:

\begin {proposition}
\label {aboutrpr}
The diagram (\ref{rpr}) is a simultaneous resolution of $\chi^{\pi},$ and the map $\tilde 
\chi^{\pi}$ is naturally a differentiable fiber bundle.
\end {proposition}

In fact, we can exhibit an explicit trivialization of the fiber bundle $\tilde 
\chi^{\pi}.$ Every flag $F \in F^{\pi}$ is of the form $u(F)F^{st},$ where $u(F) \in 
U(N)$ is unique up to right translation by an element in $U^{\pi} = U(\pi_1) \times 
\dots \times U(\pi_s).$ The map
\begin {equation}
\label {trivia}
(\tilde \chi^{\pi})^{-1}(0) \times \hh^{\pi} \to \gp \ , \ \ (x, F, \alpha) 
\to (x+Ad(u(F))\alpha, F) 
\end {equation}
is a trivialization for $\tilde \chi^{\pi}.$ Note that this trivialization is not natural, 
because it depends on our chosen basis for $\cc^N.$

Now consider the open subset $\hh^{\pi, \reg} \subset \hh^{\pi}$ where $\alpha_1, \dots, 
\alpha_s$ are pairwise distinct. Its image $\conf^0_{\pi} = \hh^{\pi, \reg}/W^{\pi}$ is 
called the {\it $\pi$-colored configuration space}. The superscript $0$ stands for 
traceless. If we don't impose the condition $\sum \pi_j \alpha_j =0,$ we get a larger 
but homotopy equivalent configuration space $\conf_{\pi}.$

\begin {definition}
The {\it $\pi$-colored braid group} is the fundamental group $Br^{\pi} = 
\pi_1(\conf^0_{\pi}).$ 
\end {definition}

Observe that the $(1^N)$-colored braid group is the usual braid group 
on $N$ strands $Br_N,$ while at the other extreme the $(N)$-colored braid group is the 
pure braid group $PBr_N.$ (In the literature $PBr_N$ is sometimes called the colored braid 
group.) In general, $Br^{\pi}$ is intermediate between these two cases, in the sense 
that the natural map $PBr_s \to Br_s$ factors through $Br^{\pi},$ in correspondence to a 
factoring of covering spaces. 

Let $\ggg^{\pi, \reg}$ be the preimage of $\hh^{\pi, \reg}/W^{\pi}$ under the map 
$\chi^{\pi}.$ For $x \in \ggg^{\pi, \reg},$ the choice of a partial flag $F$ with $(x, 
F) \in \gp$ is unique once we fix an ordering of the $\alpha_j$'s with the same $\pi_j.$ 
A quick consequence of this is the following:

\begin {proposition}
\label {regular}
The points $x \in \ggg^{\pi, \reg}$ are regular for the maps $\chi^{\pi},$ 
and the restriction of $\chi^{\pi}$ to $\ggg^{\pi, \reg}$ is naturally a differentiable 
fiber bundle over $\conf^0_{\pi}.$
\end {proposition}

Since $\gg^{\pi}$ is typically singular, we should clarify some terminology. In 
Proposition~\ref{regular} and everywhere below, a point $x \in \gg^{\pi}$ is called 
regular for 
the map $\chi^{\pi}$ if $\gg^{\pi}$ is smooth at $x,$ and the differential $(d\chi^{\pi})_x$ 
is onto; $x$ is called critical otherwise. (The same goes for points on any singular 
variety mapped to a smooth one.)

\subsection {Reorderings of a partition}
By a \textit{reordering} $\pi_{\tau}$ of the partition $\pi = (\pi_1, \dots, \pi_s)$ we 
simply mean the ordered $s$-tuple $(\pi_{\tau(1)}, \dots, \pi_{\tau(s)})$ for some 
permutation $\tau$ of $\{1,2,\dots, s\}.$ Given such a reordering, we can define partial 
flags of type $\pi_{\tau}$ so that the differences in consecutive dimensions are 
$\pi_{\tau(1)}, \dots, \pi_{\tau(s)}$ in that order. We can then define $\tilde 
\gg^{\pi_{\tau}}, \ggg^{\pi_{\tau}}, \hh^{\pi_{\tau}}$ and $W^{\pi_{\tau}}$ as in the 
previous subsection, an adjoint quotient map $\chi^{\pi_{\tau}}$ and an analogue of 
the diagram (\ref{rpr}).  Note that there are isomorphisms $\hh^{\pi}  \cong 
\hh^{\pi_{\tau}}$ (from reordering the diagonal elements) and $W^{\pi} \cong 
W^{\pi_{\tau}}.$ We denote by $\iota_{\tau}$ the induced isomorphism 
$\hh^{\pi}/W^{\pi} \to \hh^{\pi_{\tau}}/W^{\pi_{\tau}}.$ The following result will 
prove useful: 

\begin {lemma}
\label {reorder}
The subvariety $\ggg^{\pi_{\tau}} \subset \gg$ is independent of $\tau,$ and in 
particular identical to $\gg^{\pi}.$ The map $\chi^{\pi_{\tau}} : \ggg^{\pi_{\tau}} 
\to \hh^{\pi_{\tau}}/W^{\pi_{\tau}}$ equals the composite $\iota_{\tau} \circ 
\chi^{\pi}.$   
\end {lemma}

\noindent \textbf{Proof.} We first deal with the case when $\pi = (a, b)$ with 
$a \geq b, \ a+b=N$ and $\tau$ is the transposition, so that $\pi_{\tau} = (b,a).$ An 
element $x$ is in $\gg^{\pi}$ if it fixes some subspace $V= F_1  \subset 
\cc^{N}$ of dimension $a,$ and $x$ acts on $V$ diagonally by $\alpha \cdot \Id$ and on 
$\cc^N/V$ diagonally by $\beta \cdot I,$ for some $\alpha, \beta \in \cc.$ We seek to 
show that $x$ is in $\gg^{\pi_{\tau}}$ by constructing a subspace $W \subset \cc^N$ of 
dimension $b$ such that $x$ acts by $\beta \cdot \Id$ on $W$ and by $\alpha \cdot I$ on 
$\cc^N/W.$ If $\alpha \neq \beta,$ we simply take $W$ to be the $\beta$-eigenspace of $x.$ If 
$\alpha = \beta,$ then $(x - \alpha\Id)$ acts trivially on both $V$ and $\cc^N/V$ and 
therefore is determined by an induced map $\cc^N/V \to V.$ In this case we can 
choose $W$ to be any $b$-dimensional subspace of $V$ containing the image of $(x - 
\alpha \Id).$

Conversely, if we know that $X \in \gg^{\pi_{\tau}},$ then there exists a 
$b$-dimensional subspace $W \subset \cc^N$ as in the previous paragraph, and we want to 
construct $V.$ 
If $\alpha \neq \beta,$ we take $V$ to be the $\alpha$-eigenspace of $x.$ If $\alpha = 
\beta,$ then $x$ is determined by a map $\cc^N/W \to W$ induced by $(x-\alpha \Id),$ and 
we can take $V$ to be any $a$-dimensional subspace of the kernel of $(x - \alpha \Id),$ 
containing $W.$ 

The case of a general reordering follows from this by induction, using the fact that 
every permutation $\pi$ is a product of transpositions. $\hfill \fin$

\subsection {Slices at semisimple elements}
\label {sec:semis}

Let $x$ be an element of $\ggg.$ We denote its stabilizer subgroup by $G_x$ and the 
corresponding Lie subalgebra by $\ggg_x = \{y\in \ggg \ | [x,y]=0\}$. A local {\it 
transverse slice} to the adjoint orbit $\oo(x) = Gx$ of $x$ is a local complex submanifold 
$\ss \subset \ggg$ such that $x\in \ss$ and the tangent spaces at $x$ of $\oo(x)$ and 
$\ss$ are complementary. Note that $$ T_x \oo(x) = ad(\ggg)x = \{[y,x] \ | \ y \in 
\ggg\}.$$

Let $x$ be semisimple. The dimensions of the eigenspaces of $x$ form a partition 
$\sigma=(\sigma_1, \dots, \sigma_l)$ of $N.$ The Lie algebra $\ggg_x$ splits into a direct 
sum of its center $Z(\ggg_x)$ and a reduced part $\ggg_x^{\red}=[\ggg_x, \ggg_x],$ which is 
a product of factors $\ggg_x^{\red} [ i ]$ isomorphic to $\sl(\sigma_i, \cc).$ There is 
an 
adjoint quotient map for $\ggg_x^{\red},$ denoted by 
\begin {equation}
\label {chix}
\chi_x^{\red}: \ggg^{\red}_x \to \ggg^{\red}_x/G_x = \hhh^{\red}_x/W^{\red}_x.
\end {equation}
Here $\hhh^{\red}_x \subset \hhh$ corresponds to block diagonal matrices such that each 
block has trace zero, and $W^{\red}_x= W_{\sigma}.$

Using transverse slices, it can be shown that a fibered version of $\chi_x^{\red}$ gives a 
a local model for the adjoint quotient map $\chi$ near the orbit $\oo(x).$ We sketch this 
here (following \cite[Section 2(B)]{SS}), and then explain how to get a local model 
for the restricted map $\chi^{\pi}$ near $\oo(x)$ when that orbit lies in $\gg^{\pi}.$

A canonical transverse slice for $x$ is $\ss^{ss} = x + \ggg_x.$ 
This has the property that, given any other slice $\ss,$ there is a canonical local 
isomorphism between $\ss$ and $\ss^{ss}$ obtained as follows. We take a 
neighborhood $V$ of zero inside $ad(\ggg)x.$ Then $exp(V) \subset G$ is tranverse to $G_x$ 
and we can use the composition
\begin {equation}
\label {compose}
\begin {CD}
\ss \hookrightarrow \ggg @>\operatorname{local \ \cong}>> exp(V) \times \ss^{ss} 
@>\operatorname{projection}>> \ss^{ss}.
\end {CD}
\end {equation}

The slice $\ss^{ss}$ is $G_x$-invariant, and therefore the map
\begin {equation}
\label {fibered}
 G \times_{G_x} \ss^{ss} \to \ggg, \ (g, y) \to Ad(g)y
\end {equation}
is a local isomorphism between neighborhoods of $G/G_x \times \{x\}$ and $\oo(x).$
Using the fact that $G_x$ acts trivially on the first summand in the splitting 
$\ggg_x = Z(\ggg_x) \oplus \ggg_x^{\red},$ we obtain a commutative diagram
\begin {equation}
\label {cairo}
\begin {CD}
Z(\ggg_x) \times (G \times_{G_x} \ggg_x^{\red}) @>>> \ggg \\
@VVV @V{\chi}VV \\
Z(\ggg_x) \times \hhh^{\red}_x/W^{\red}_x @>>> \hhh/W.
\end {CD}
\end {equation}

The top $\rightarrow$ is (\ref{compose}), and the 
bottom one is the induced local isomorphism on quotients, in both cases identifying 
$\ss^{ss}$ with $\ggg_x  = Z(\ggg_x) \times \ggg_x^{\red}$ by translation. The vertical map 
on the left is identity times the quotient map, which means that $\chi$ looks locally 
like a fibered version of (\ref{chix}).

Let us go back to the restricted partial resolution (\ref{rpr}) and assume that $x \in 
\ggg^{\pi}.$ Since $x$ is semisimple, it must be conjugate to a matrix $diag(\alpha_1, 
\dots, \alpha_1, \alpha_2, \dots, \alpha_2, \dots, \alpha_s) \in \hh^{\pi}.$ (Note that 
this diagonal matrix is not typically unique, because of the $W^{\pi}$ symmetry.) The 
set of indices $\{1,2,\dots, s\}$ has a decomposition into a disjoint union $A_1 
\amalg \dots \amalg A_l$ such that $\alpha_j = \alpha_k$ if and only if $j$ and $k$ are 
in the same $A_i.$ The sizes $\sigma_i$ of the sets $A_i$ form the partition $(\sigma_1, 
\dots, \sigma_l)$ of $N$ that was described for any semisimple element. Since $x \in 
\ggg^{\pi},$ we can be more specific and say that the partition $\pi$ of $N$ ``breaks'' in 
the following sense:

\begin {definition}
Let $\pi=(\pi_1, \dots, \pi_s)$ and $\sigma = (\sigma_1, \dots, \sigma_l)$ be two 
partitions of the same positive integer $N.$ A breaking $b = (b[1], b[2], \dots, b[l])$ 
of $\pi$ according to $\sigma$ consists of partitions $b[i]$ of each $\sigma_i$ such that 
their concatenation is a reordering of the partition $\pi.$ The set of breakings of $\pi$ 
according to $\sigma$ is denoted $\bb_{\pi \sigma}.$
\end {definition} 

In our situation, $\pi$ breaks into $b = (b[1], b[2], \dots, b[l]),$ where $b[i]$ is 
composed of all the $\pi_j$'s with $j \in A_i.$ By varying the diagonal matrix in 
$\hh^{\pi}$ conjugate to $x,$ all the breakings in $\bb_{\pi \sigma}$ can occur.

\begin {example}
\label {ex1}
If $\pi=(2,1,1)$ and $\sigma=(2,2),$ then there are two possible breakings of $\pi$ 
according to $\sigma: \ \bigl( b[1]=(2),\ b[2]=(1,1) \bigr)$ and $\bigl( b[1]=(1,1),\  
b[2]=(2) \bigr).$ More informally, we say that the latter breaking, for example, consists 
of the first $2$ in $\sigma=(2,2)$ being broken as $1+1$ and the second $2$ being broken 
trivially.
\end {example}

\begin {example}
\label {ex2}
For $m,n \geq 1,$ there is a unique breaking $b$ of $\pi =
(1^m(n-1)^m)$ according to $\sigma=(1^{m-1}(n-1)^{m-1}n),$ namely $n$ breaks as $(n-1)+1$ 
and all the $1$'s and $(n-1)$'s break trivially.
\end {example}

\begin {example}
\label {ex3}
For $m, n \geq 1,$ there are exactly $m$ breakings of $\pi=(1^{m+n-1}(n-1)^{m-1})$ 
according to $\sigma = (1^{m-1}(n-1)^{m-1}n):$ either the $n$ summand breaks as $1+ \dots 
+1$ and the rest break trivially, or the $n$ summand breaks as $(n-1)+1,$ one of the 
$(n-1)$'s breaks as a sum of $1$'s and the rest break trivially.
\end {example}

\begin {example}
\label {ex4}
For $m \geq 2$ and $n > 3,$ there is a unique breaking of $\pi =
(1^m(n-1)^m)$ according to $\sigma=(1^{m-2}2(n-1)^m),$ namely $2$ breaks as $1+1$ and the 
rest trivially.
\end {example}

\begin {example}
\label {ex5}
For $m \geq 2$ and $n \geq 1,$ there is again a unique breaking of $\pi =
(1^m(n-1)^m)$ according to $\sigma = (1^{m-2}(n-1)^{m-1}(n+1)),$ namely $(n+1) = 
(n-1) + 1 +1$ and the rest trivial. 
\end {example}

With respect to transverse slices at $x \in \gg^{\pi}$, the first observation which needs 
to be made is that the local isomorphism (\ref{compose}) moves points only inside their 
adjoint orbits, and therefore induces a canonical local isomorphism
\begin {equation}
\label {caniso}
 \ss \cap \ggg^{\pi} \cong \ss^{ss} \cap \ggg^{\pi}.
\end {equation}

Second, the same reasoning applies to the map (\ref{fibered}), giving a local isomorphism
\begin {equation}
\label {fibered2}
G \times_{G_x} (\ss^{ss} \cap \ggg^{\pi}) \to \ggg^{\pi}.
\end {equation}

The left hand side of (\ref{fibered}) can be made more explicit as follows. Recall that 
$\ss^{ss} = x + \ggg_x$ and $\ggg_x$ decomposes as $Z(\ggg_x) \times \prod_i 
\ggg_x^{\red}[i],$ with $\gg_x^{\red}[i] \cong \sl(\sigma_i, \cc).$ We have a restricted 
partial simultaneous resolution associated to each factor $\ggg_x^{\red}[i]$ and the 
corresponding partition $b[i]$ in a breaking of $\pi$ according to $\sigma.$ The role of 
$\chi^{\pi}$ in the left half of the diagram (\ref{rpr}) is played by maps 
$$\chi^{b[i]} : \gg_x^{\red}[i]^{b[i]} \to  \hh_x^{\red}[i]^{b[i]}/W^{b[i]}$$
for each $i.$

\begin {lemma}
\label {hello}
Translation $y \to (y-x)$ induces an identification of 
$\ss^{ss} \cap \ggg^{\pi}$ with
$$Z(\gg_x) \times \bigcup_{b \in \bb_{\pi \sigma}} \prod_i \gg_x^{\red}[i]^{b[i]}.$$  
\end {lemma}

\noindent \textbf{Proof.} First observe that an element in $\gg_x$ preserves each eigenspace 
$E_i$ of $x.$ The direct sum of all $E_i$ is $\cc^N,$ the sizes of $E_i$ form the partition 
$\sigma,$ and $ \gg_x^{\red}[i]$ is the Lie alegbra of traceless endomorphisms of $E_i.$

Now choose some $y$ such that $y -x $ is in $Z(\gg_x) \times \prod_i 
\gg_x^{\red}[i]^{b[i]}$ for some $b \in  \bb_{\pi \sigma}.$ We have that $y$ is in $\gg_x,$
and $y$ acts on each $E_i$ by a central element plus something in $\gg_x^{\red}[i]^{b[i]}$ 
This means that we can find partial flags of type $b[i]$ made of subspaces of $E_i$ that are 
preserved by $y$ and such that $y$ acts diagonally on consecutive quotients. Putting these flags 
together (by taking direct sums) we form a partial flag made of subsets of $\cc^N.$ Its type is 
a reordering of the partition $\pi.$ Lemma~\ref{reorder} implies that $y$ must be in 
$\gg^{\pi}.$

Conversely, take some $y \in (x + \gg_x) \cap \ggg^{\pi},$ sufficiently close to $x.$ Since 
$y$ is in $\gg^{\pi}$ it must preserve a partial flag $F$ of type $\pi$ and act diagonally 
on the consecutive quotients. We also know that $y$ preserves each $E_i,$ and therefore the 
intersections of $F$ with each $E_i$ are also preserved, with diagonal actions on 
consecutive quotients. The resulting partial flags are of type $b[i],$ for some breaking of 
$\pi$ according to $\sigma.$ We get that the traceless part of $y|_{E_i}$ is in  
$\gg_x^{\red}[i]^{b[i]}$ for all $i. \ \hfill \fin$

\medskip

Using the result of Lemma~\ref{hello}, we get a restricted version of (\ref{cairo}):
\begin {equation}
\label {cairo2}
\begin {CD}
Z(\ggg_x) \times \bigl (G \times_{G_x}  \bigcup_{b \in \bb_{\pi \sigma}}(\prod_i  
\gg_x^{\red}[i]^{b[i]}) \bigr ) @>>> \ggg^{\pi} \\
@VVV @V{\chi^{\pi}}VV \\
Z(\ggg_x) \times  \bigcup_{b \in \bb_{\pi \sigma}} (\prod_i  \hh_x^{\red}[i]^{b[i]}/W^{b[i]}) 
@>>> \hhh^{\pi}/W^{\pi}.
\end {CD}
\end {equation}

The map going vertically on the left is identity times a quotient map (coming from the 
product of all $\chi^{b[i]}$'s), and the horizontal arrows are local isomorphisms. It 
follows that locally near $\oo(x) \subset \gg^{\pi},$ the map $\chi^{\pi}$ looks like a 
fibered version of the union of products of the $\chi^{b[i]}$'s, taken over all $b \in \bb_{\pi 
\sigma}.$

\subsection {Invariant slices at nilpotent elements}

Now let $x \in \ggg$ be nilpotent. One way to construct transverse slices at $x$ with nice 
global properties is using the Jacobson-Morozov lemma, which claims the existence of a 
triple $(n^+=x, n^-, h)$ of elements of $\ggg$ such that
\begin {equation}
\label {nueva}
[h, n^+]=2n^+, \ [h, n^-]=-2n^-, \ [n^+, n^-]=h.
\end {equation}

This gives a representation of $\sl(2,\cc)$ into $\ggg,$ which in turn 
produces a splitting $\ggg = \im ad(n^+) \oplus \ker ad(n^-).$ It follows that:

\begin {lemma}
The affine subspace $\sjm = n^+ + \ker ad(n^-)$ is a local transverse slice at $x.$  
\end {lemma}

Furthermore, we have a linear $\cc^*$-action on $\ggg$
\begin {equation}
\label {act1}
 \lambda_r: \cc^* \to \operatorname{Aut}(\ggg), \ \ r\in \cc^* \ : \  y \to 
r^2Ad(r^{h})y.
\end {equation}

This restricts to an action on the Jacobson-Morozov slice $\sjm,$ and has the property 
that it contracts $\sjm$ to $0$ as $r\to 0.$ Using this action we can show that the 
behavior of $\sjm$ globally is determined by the behavior near $x.$ In particular, 
$\sjm$ intersects all adjoint orbits transversely.

In \cite{SS}, Seidel and Smith use a more general notion:

\begin {definition}
Given a triple $(n^+, n^-, h)$ satisying (\ref{nueva}), a $\lambda$-invariant slice is an 
affine subspace $\ss \subset \ggg$ invariant under $\ggg$ and such that it represents a local 
transverse slice at $n^+=x.$ 
\end {definition}

In fact, $\lambda$-invariant slices share many properties with Jacobson-Morozov slices:

\begin {proposition}[Lemma 11(i) and Lemma 14 in \cite{SS}] 
\label {invslice}
(i) A $\lambda$-invariant slice $\ss$ has transverse intersections with all adjoint 
orbits. (ii) For $\lambda$-invariant $\ss$ and Jacobson-Morozov slice $\sjm$ (possibly 
coming from a  
different choice of a triple with $n^+=x$), there exists a (noncanonical) 
$\cc^*$-equivariant isomorphism $\ss \cong \sjm,$ which moves points only in their 
adjoint orbits.
\end {proposition}

Now assume that $x \in \gpi.$ Because the diagram (\ref{rpr}) is $G$-equivariant, it 
behaves well under intersecting with a $\lambda$-invariant transverse slice $\ss.$ 
More precisely, we have:

\begin {proposition}
\label {big}
(i) The isomorphism in Proposition~\ref{invslice} (ii) maps $\ss \cap \gpi$ into $\sjm 
\cap \gpi.$ (ii) A point of $\ss \cap \gpi$ is a critical point of $\chi^{\pi}$ iff it 
is a critical point of $\chi^{\pi}|_{\ss}.$ (iii) Let $\tilde \ss$ be the 
preimage of $\ss$ in $\tilde \ggg.$ Then 
\begin {equation}
\label {slice}
\begin {CD}
\tilde \ss \cap \gp @>>> \ss \cap \gpi \\ @V{\tilde \chi^{\pi}}VV @V{\chi^{\pi}}VV \\ 
\hh^{\pi} 
@>>> \hh^{\pi}/W^{\pi} \end {CD}
\end {equation}
is a simultaneous resolution. (iv) $\tilde \chi^{\pi}: \tilde \ss \cap \gp \to 
\hh^{\pi}$ is naturally a differentiable fiber bundle. 
\end {proposition}

Point (i) is a direct implication of Proposition~\ref{invslice} and the fact that 
$\gpi$ is a union of adjoint orbits. The proofs of the other assertions use the 
properties of $\lambda,$ and are completely analogous to those in \cite[Lemma 11 and 
Section 3(D)]{SS}.

Putting propositions \ref{regular} and \ref{big}(ii) together we get that the points of 
$\ss \cap \ggg^{\pi, \reg}$ are regular for the restriction of $\chi^{\pi}$ to $\ss.$ 
In fact, we also have:

\begin {proposition}
\label {luna}
The restriction of $\chi^{\pi}$ to $\ss \cap \ggg^{\pi, \reg}$ is naturally a 
differentiable fiber bundle over $\conf^0_{\pi}.$
\end {proposition}

\subsection {An example} \label {sec:exa} Consider the following nilpotent in $\sl(n+1, 
\cc):$
$$ x=n^+ = \left( \begin {array}{ccccc} 0 & 1 &  &  & \\
 & 0 & \ & \ & \\
 & \ & 0 & \ & \ \\
 & \ & \  & \ldots  &  \ \\
 & \ & \ & \ & 0
\end {array} \right ) .$$
We complete it to a Jacobson-Morozov triple with:
$$n^- =  \left( \begin {array}{ccccc} 0 &  &  &  & \\
1 & 0 & \ & \ & \\
 & \ & 0 & \ & \ \\
 & \ & \  & \ldots  &  \ \\
 & \ & \ & \ & 0
\end {array} \right ) \ , \  
h = \left( \begin {array}{ccccc} 1 &  &  &  & \\
 & -1 & \ & \ & \\
 & \ & \ 0 & \ \\
 & \ & \  & \ldots  &  \ \\
 & \ & \ & \ & 0
\end {array} \right ).$$

The associated slice $\sjm$ consists of all matrices of the form
\begin {equation}
\label {sjmm}
  \left( \begin {array}{ccccc} \alpha & 1 & 0 & \cdots & 0\\
a_{11} & \alpha & a_{12} & \cdots & a_{1n}\\
a_{21} & 0 & a_{22} & \cdots & a_{2n} \\
 & \cdots & \  & \cdots  &  \ \\
a_{n1} & 0 & a_{n2} & \cdots & a_{nn}
\end {array} \right ),
\end {equation}
with $\alpha, a_{ij} \in \cc$ such that $2\alpha + a_{22} + \cdots + a_{nn} = 0.$

Consider the partition $\pi = (n-1,1,1)$ of $n+1.$ Then $\gg^{\pi} \subset \gg=\sl(n+1, 
\cc)$ consists of the traceless matrices which admit an $(n-1)$-dimensional eigenspace. For 
future reference, we will denote by $\xx_n$ the intersection $\sjm \cap 
\gg^{\pi}$ in this case, 
and by $\mapp$ the restriction of $\chi^{\pi}$ to $\sjm.$  More precisely,
\begin {equation}
\label {mapp}
\mapp : \xx_n \to \hh^{(n,1)}/W^{(n,1)} \cong \sym^2(\cc),
\end {equation}
where the identification with $\sym^2(\cc)$ comes form taking the diagonal matrix 
$diag(z_1, z_2, 
z_3,$ $\cdots, z_3) \in \hh^{\pi}$ to the pair $(z_1, z_2).$ (Note that $z_3 = -(z_1 + 
z_2)/(n-1).$) Of course, we can also further identify $\sym^2(\cc)$ with $\cc^2$ by taking 
$(z_1, z_2)$ to $(z_1 + z_2, z_1z_2).$ 

\begin {remark}
\label {act}
More canonically, we define a linear quadruple $T=(F_1, T_1, T_2, f)$ to be the data 
consisting of complex vector spaces $F_1, T_1, T_2$ of dimensions $1, n-1, 1,$ 
respectively, and a linear isomorphism $f: T_2 \to F_1.$ To every linear quadruple we can 
associate an $(n+1)$-dimensional complex vector space $E = F_1 \oplus T_2 \oplus T_1$ and a 
flag $F$ of type $1+(n-1)+1$ given by
\begin {equation}
\label {drapel}
 0 \subset F_1 \subset F_2 \subset E, \ \ F_2 = F_1 \oplus T_1.
\end {equation}

The isomorphism $f$ produces a nilpotent $n^+=x \in \sl(E)$ with $F_1 = \im (n^+), \ F_2 = 
\ker (n^+),$ while its inverse $f^{-1}$ gives rise to a nilpotent $n^- \in \sl(E)$ with 
$T_2 = \im (n^-), \ T_2 \oplus T_1 = \ker(n^-).$ We can complete this to a Jacobson-Morozov 
triple by setting $h=[n^+, n^-].$ It follows that to $T$ we can naturally associate a space 
$\xx(T)$ isomorphic to $\xx_n.$ There is still a natural map $q: \xx(T) \to \sym^2(\cc).$ 
\end {remark}

\subsection {More on nilpotents} The adjoint orbits of nilpotents in $\ggg$ are 
classified by the partitions of $N.$ Given the partition $\pi,$ the parabolic subalgebra 
$\ppp(F^{st})$ decomposes into a Levi part and a nilpotent part. The nilpotent part 
$\mathfrak{n}_{\pi}$ has a unique dense orbit $\oo(x_{\pi}),$ where the nilpotent 
element $x_{\pi}$ can be taken to be the Jordan matrix with Jordan cells whose sizes 
give the dual partition $\pi^*.$ For example, when $\pi=(N),$ the corresponding 
parabolic subalgebra is all of $\ggg,$ the nilpotent $x_{(N)}$ has $N$ Jordan blocks of 
size one, and $(x_{(N)}) = \{0\}.$ At the other extreme, when $\pi = (1^N),$ the 
corresponding parabolic subalgebra consists of the upper triangular matrices, the 
nilpotent $x_{(1^N)}$ is just one Jordan block of size $N,$ and $\oo(x_{(1^N)})$ is the 
so-called regular nilpotent orbit. In general, a nilpotent in $\oo(x_{\pi})$ is called 
of type $\pi.$

There is a well-known partial ordering on partitions of $N.$ Consider two partitions 
$\pi= (\pi_1, \pi_2, \dots)$ and $\rho = (\rho_1,\rho_2, \dots),$ with $\pi_1 \geq 
\pi_2 \geq \dots$ and $\rho_1 \geq \rho_2 \geq \dots$ as before, both adding up to $N.$
We write $\pi \preccurlyeq \rho$ if $\pi_1 + \dots + \pi_j \leq \rho_1 + \dots + \rho_j$ for 
all $j \geq 1.$ The following lemma is elementary:

\begin {lemma}
\label {lemm}
The closure of $\oo(x_{\pi})$ in $\ggg$ consists of the union of all orbits 
$\oo(x_{\rho})$ with $\pi \preccurlyeq \rho.$ These are exactly the nilpotent orbits that 
appear in $\gpi.$
\end {lemma}

\subsection {Slices at general points} \label {sec:gen} A general element $x \in \gg$ has 
a unique decomposition $x=x_s + x_n$ into a semisimple and a nilpotent part. Just as in 
Section~\ref{sec:semis}, $\gg_{x_s}$ decomposes into the direct sum of $Z(\gg_{x_s})$ and 
$\gg_{x_s}^{\red} = \oplus_i \gg_{x_s}^{\red}[i],$ and each $\gg_{x_s}^{\red}[i]$ is 
composed of the traceless endomorphisms of an eigenspace $E_i$ of $x_s.$ The sizes of 
these eigenspaces form a partition $\sigma$ of $N.$ The nilpotent part $x_n$ always lies 
in $\gg_{x_s}^{\red};$ denote by $x_n[i]$ its piece in $\gg_{x_s}^{\red}[i].$ We choose 
Jacobson-Morozov slices $\ss^{\jm, \red}[i]$ at $x_n[i]$ in $\gg_{x_s}^{\red}[i],$ and let 
$\ss^{\jm, \red}$ be their direct sum.

In \cite[Section 2(D)]{SS} it is proved that $x_s + (Z(\gg_{x_s}) \times \ss^{\jm, \red})$
is a transverse slice at $x$ in $\ggg.$ Using this, it follows that locally near $x,$ the 
adjoint map $\chi$ looks like a linear projection times the product of the 
restriction of the adjoint map of  $\gg_{x_s}^{\red}$ to $\ss^{\jm, \red}.$ This local model 
enables one to describe explicitly which points of $\gg$ are regular for the map $\chi:$ 
namely, the ones for which $x_n$ is a regular nilpotent in $\gg_{x_s}^{\red}.$

When $x \in \gg^{\pi},$ by Lemma~\ref{hello} we have that the nilpotents $x_n[i]$ are in 
$\gg_{x_s}^{\red}[i]^{b[i]}$ for some breaking $b \in \bb_{\pi \sigma}.$ We can restrict the 
local 
model to $\gg^{\pi}$ along the lines of Section~\ref{sec:semis}, with the following result:

\begin {lemma}
\label {locmodel}
A point $x = x_s + x_n \in \gg^{\pi}$ is regular for the map $\chi^{\pi}$ if and only if 
there exists a breaking $b = (b[1], \dots, b[l])$ of $\pi$ according to $\sigma,$ such that 
$x_n[i] \in \gg_{x_s}^{\red}[i]$ is a nilpotent of type $b[i].$  
\end {lemma}

Furthermore, note that a fiber of the map $\chi^{\pi}:\gg^{\pi} \to \hh^{\pi}/W^{\pi}$ is 
uniquely determined by the semisimple orbit $\oo(x_s)$ that it contains, corresponding to some 
partition $\sigma.$ Using Lemma~\ref{lemm}, we get that:

\begin {lemma}
\label {structure}
The fiber of $\chi^{\pi}$ containing $\oo(x_s)$ is composed of the orbits of $x_s + x_n,$ 
where $x_n$ is a nilpotent decomposing into $x_n[i] \in \gg_{x_s}^{\red}[i]$ for some 
breaking $b = (b[1], \dots, b[l]) \in \bb_{\pi \sigma},$ and with the types $\rho[i]$ of 
$x_n[i]$ satisfying $b[i] \preccurlyeq \rho[i].$ 
\end {lemma}

\subsection {The topology of the fibers} 
\label{sec:top}
Our main objects of interest are the fibers of 
the bundle map $\chi^{\pi}|_{\ss \cap \ggg^{\pi, \reg}}$ from Proposition~\ref{luna}. By 
Proposition~\ref{big}, they are homeomorphic to an arbitrary fiber of the 
restriction of $\tilde \chi^{\pi}$ to $\tilde \ss \cap \gp.$ We pick the nilpotent fiber 
\begin {equation}
\label {nnn}
\nnn_{\rho \pi} = \bigl(\tilde \chi^{\pi} |_{\tilde \ss \cap \gp}\bigr)^{-1}(0).
\end {equation}

Here we assumed that the nilpotent $x$ where we took the slice $\ss = \ss_{\rho}$ was 
in the orbit of $x_{\rho},$ for some partition $\rho$ with $\pi \preccurlyeq \rho.$ Observe 
that by Proposition~\ref{invslice}(ii), the topology is independent of which invariant 
slice we choose.

Using Lemma~\ref{lemm} we see that there is a natural map
$$ \Pi: \nnn_{\rho \pi} \to \ss_{\rho} \cap \overline{\oo(x_{\pi})}, \ \ (x, F) \to x. $$

Using the $\cc^*$-action $\lambda,$ it can be shown that $\nnn_{\rho \pi}$ deformation 
retracts into its ``compact core''
$$ \Pi^{-1}(x) = \{F \in \ff^{\pi} \ | \ x \text{ fixes } F \text { and acts 
trivially on all quotients } F_j/F_{j-1} \}. $$

$\Pi^{-1}(x)$ is called the {\it Spaltenstein variety} and appeared first in \cite{Sp}. 
It is a projective variety of half the dimension of $\nnn_{\rho \pi}.$
 
\begin {example}
When $\pi=(1^N), \ \Pi^{-1}(x)$ is the Springer variety of complete flags fixed by a 
given nilpotent element. 
\end {example}

\begin {example}
When $\rho = (N),$ we have $x = 0$ and $\Pi^{-1}(x)$ is the partial flag variety 
$\ff^{\pi}.$ In this case $\nnn_{\rho \pi}$ can be identified with the cotangent bundle 
$T^*\ff^{\pi}.$
\end {example}

In general, the real dimension $d_{\rho \pi}$ of the Spaltenstein variety  
$\Pi^{-1}(x)$ is given by the formula
\begin {equation}
\label {dimension}
d_{\rho \pi} = \sum_{j\geq 1} \rho_j^2 - \sum_{j\geq 1} \pi_j^2.
\end {equation}

In \cite[Sections 3.4 and 3.5]{BM}, Borho and MacPherson studied the rational cohomology of 
a general Spaltenstein variety. For example, they found that the number of $d_{\rho 
\pi}$-dimensional irreducible components of $\Pi^{-1}(x)$ (which is the middle dimensional 
Betti number of $\nnn_{\rho \pi}$) is equal to the Kostka number $K_{\rho \pi}$ which 
counts semistandard Young tableaux of shape $\rho$ and weight $\pi.$ For more information 
on the combinatorics of Kostka numbers and their relevance to representation theory, we 
refer to \cite{McD}. A different viewpoint on these cohomology computations is in 
\cite{N2}.

\section {The relevant affine variety and its degenerations}
\label{sec:av}

From now on we specialize the discussion in the previous section to the case of interest to 
us, by setting $\pi = (1^m(n-1)^m)$ and $\rho = (n^m),$ with $N = nm.$ Note that $W^{\pi} = 
S_m \times S_m$ for $n > 2.$ The case $n=2$ is slightly different, because  $W^{\pi} = 
S_{2m}.$ Since that case was treated by Seidel and Smith in \cite{SS}, we shall restrict our 
attention to $n > 2.$ (The case $n=1$ is trivial, and the resulting link invariant is $\zz$ 
for any link.)

\subsection {Motivation} \label {sec:motiv}
Let us briefly explain why our choice of $\pi = (1^m(n-1)^m)$ and 
$\rho = (n^m)$ is natural if we aim at constructing an analogue of Khovanov-Rozansky 
homology. As mentioned in the introduction, the Euler characteristic of Khovanov-Rozansky 
homology is the link polynomial $P_{\nn}.$ Let $V$ be the standard representation of 
$\sl(n, \cc).$ The polynomial $P_{\nn}$ can be defined as follows: we present 
the link as the closure of a braid $b$ on $m$ strands, associate to $b$ a map 
$F_{b}(q): V^{\otimes m} \to V^{\otimes m}$ depending on a quantum parameter $q,$ and 
then take its trace. This is equivalent to looking at the image of $1$ under a map 
\begin {equation}
\label {poly}
\begin {CD}
\cc \to  V^{\otimes m} \otimes (V^{\otimes m})^* @>{F_{b}(q) \times id}>> V^{\otimes m} 
\otimes (V^{\otimes m})^* \to \cc.
\end {CD} 
\end {equation}

The composite factors through the space
of invariants 
$$\inv(m,n) = \hom_{\sl(n, \cc)} (\cc, V^{\otimes m} \otimes (V^{\otimes 
m})^*).$$  

The dimension of $\inv(m,n)$ is equal to $d(m,n)=$ the number of permutations of $m$ 
elements with longest increasing subsequence of length $\leq n.$ (A table of the values of 
$d(m,n)$ can be found in \cite{O}.) According to Khovanov's principles for 
categorification \cite{Kh2}, we should look for a triangulated category whose Grothendieck 
group has dimension $d(m,n)$ or, more geometrically, for an exact symplectic manifold 
whose middle dimensional Betti number is $d(m,n).$ As explained at the end of 
Section~\ref{sec:top}, the middle dimensional Betti numbers of the smooth fibers of 
$\chi^{\pi}$ are the Kostka numbers. By the Schensted correspondence \cite{Sc}, $d(m,n)$ 
is the same as the Kostka number $K_{(n^m), (1^m(n-1)^m)},$ which explains our choice for 
$\pi$ and $\rho.$

\subsection {A few properties}
\label {sec:props}

From the orbit of $x_{\rho}$ we choose the following nilpotent element, written as a $m 
\times m$ matrix made of $n \times n$ blocks:
$$ N_{m,n} =  \left( \begin {array}{ccccc} 0 & I &  &  & \\
 & 0 & I & \ & \\
 & \ & \ & \ldots & \ \\
 & \ & \  & \  &  I \\
 & \ & \ & \ & 0
\end {array} \right ). $$
Here $I$ and $0$ are the $n$-by-$n$ identity and zero matrix, respectively.

We complete this to a Jacobson-Morozov triple $(N^+=N_{m,n}, N^-, H)$ with
$$ H =  \left( \begin {array}{ccccc} (m-1)I &  &  &  & \\
 & (m-3)I & \ & \ & \\
 & \ & (m-5)I & \ & \ \\
 & \ & \  & \dots  &  \ \\
 & \ & \ & \ & (-m+1)I
\end {array} \right ), $$
$$ N^- =  \left( \begin {array}{cccccc} 0 & & &  &  & \\
(m-1)I & 0 & \ & & \\ & \\
 & 2(m-2)I & 0 & \ & \ & \ \\
 & \ & 3(m-3)I & \  & \  &  I \\
& \ & \ & \dots & \ & \\
 & \ & \ & \ & (m-1)I & 0
\end {array} \right ), $$
written in the same form. The induced $\cc^*$ action on $\ggg = \sl(mn)$ is:
\begin {equation}
\label {action}
 \lambda_r \ : Y \to  \left( \begin {array}{ccccc} r^2Y_{11} & Y_{12} & \dots &  &  
r^{4-2m} Y_{1m}\\
r^4 Y_{21} & r^2Y_{22} &  & \ & \\
\dots & \ & \dots & \ & \ \\
 & \ & \  & \  &  Y_{m-1,m} \\
r^{2m}Y_{m1} & \ & \dots & r^4Y_{m,m-1} & r^2Y_{mm}
\end {array} \right ).
\end {equation}

Consider the affine space
\begin {equation}
\label {smn}
 \ss_{m, n} = \left \{ \left( \begin {array}{ccccc} Y_{11} & I 
&  &  & \\
Y_{21} & 0 & I & \ & \\
\ldots & \ & \ & \ldots & \ \\
 & \ & \  & \  &  I \\
Y_{m1} & \ & \ & \ & 0
\end {array} \right ): Y_{m1} \in \gl(n) \text{ for } m>1\ , \ Y_{11} \in \sl(n)  \right 
\}.
\end {equation}

\begin {lemma}
$\ss_{m,n}$ is a $\lambda$-invariant slice to the adjoint orbit of $N_{m, n}$ in $\sl(mn).$
\end {lemma}

\noindent \textbf{Proof. } The adjoint orbit of $N_{m,n}$ has (complex) dimension $(mn)^2 - 
mn^2,$ and $\ss_{m,n}$ has complementary dimension $mn^2 -1.$  Also, $\lambda$-invariance is 
clear from (\ref{action}). Therefore the only thing to show is that $\ss_{m,n}$ intersects 
the tangent space to the adjoint orbit trivially. Let $N^+ + X$ be an element in their 
intersection. Then all but the first $n$ columns of $X$ are nonzero, and $X$ is of the form 
$[N^+, Z]$ for some $Z \in \ggg.$ A quick calculation shows that $Z$ must be upper 
triangular in our block form, which implies that $X = [N^+, Z] =0. \hfill \fin$
\medskip

We are interested the simultaneous resolution (\ref{slice}) with $\ss = 
\ss_{m,n}$ and $\pi = \pi(m,n) = (1^m(n-1)^m),$ and in particular in understanding the 
fibers of
$$ \chi^{\pi}|_{\ss}: \ss \cap \ggg^{\pi} \longrightarrow \hh^{\pi}/W^{\pi}.$$

We denote a typical element of $\hh^{\pi}/W^{\pi}$ by $\tau = (\lambdav, \muv).$ Here 
$\lambdav = ( \lambda_1, \dots, \lambda_m )$ and $\muv = (\mu_1, \dots, \mu_m )$ are 
unordered $m$-tuples of elements of $\cc,$ with repetitions allowed and such that $\sum 
\lambda_i + (n-1)\sum \mu_j =0.$ Strictly speaking, the element in $\hh^{\pi}/W^{\pi}$ 
associated to $\tau$ is the class of the diagonal matrix $D_{\tau} \in \sl(mn)$ with 
eigenvalues $\lambda_1, \dots, \lambda_m$ with multiplicity $1$ (called {\it thin} 
eigenvalues), and $\mu_1, \dots, \mu_m,$ each with multiplicity $n-1$ (called {\it thick} 
eigenvalues). The fiber $\ymnt=(\chi^{\pi}|_{\ss})^{-1}(\tau)$ is then the intersection of 
the adjoint orbit $\oo_{\tau}$ of $D_{\tau}$ with our chosen slice:
$$ \ymnt = (\chi^{\pi}|_{\ss_{m,n}})^{-1}(\tau) = \ss_{m,n} \cap \oo_{\tau}.$$

Inside of $\hh^{\pi}/W^{\pi}$ we have the bipartite configuration space $\conf^0_{\pi} = 
\confb^0_m$ mentioned in the introduction, which corresponds to all the $\lambda_i$ and 
$\mu_j$ being distinct. According to Proposition~\ref{luna}, the fibers $\ymnt$ are smooth 
for $\tau \in \confb^0_m.$  

The following two lemmas are generalizations of lemmas 18 and 19 in \cite{SS}:

\begin {lemma} 
\label {newl}
For any $Y \in \ss_{m,n}$ and $\mu \in \cc,$ projection to the first $n$ coordinates 
produces an injective map $\ker (\mu I - Y) \to \cc^n.$
\end {lemma}

\noindent \textbf{Proof. } Suppose the contrary is true, so that $\ker(\mu I - Y)$ has nonzero 
intersection with  $\{0\}^n \times \cc^{(m-1)n}$. Using
the $\cc^*$-action, one sees that the same holds for $\ker(r^2\mu I
- \lambda_r(Y)).$ In the limit $r \rightarrow 0$ we obtain a nonzero
element in $\ker(N^+) \cap (\{0\}^n \times \cc^{(m-1)n})$, which is a contradiction.
$\hfill \fin$
\medskip

\begin {lemma} 
\label {identh}
Let $\tau=(\lambdav, \muv) \in \hh^{\pi}/W^{\pi}$ be such that 
$\lambdav=(\lambda_1=0, \lambda_2, \dots, \lambda_m), \muv=(\mu_1=0, \mu_2, \dots, \mu_m).$ 
Denote by $\bar \tau=\bigl( (\lambda_2, \dots, \lambda_m), (\mu_2, \dots, \mu_m) \bigr).$ 
Then the subspace of $Y \in \ymnt$ such that $\ker(Y)$ is $n$-dimensional can be canonically 
identified with $\yy_{m-1, n, \bar \tau}.$ 
\end {lemma}

\noindent \textbf{Proof. } From the previous lemma we know that a vector in $\ ker(Y)$ 
is uniquely determined by its first $n$ entries. There are $n$ linearly independent such 
vectors if and only if $Y_{m1} = 0$. The subspace of $\ss_{m,n}$ with $Y_{m1}=0$ can be 
identified with $\ss_{m-1,n}$ in a straightforward way. If $\bar Y \in \yy_{m-1, n, \bar 
\tau} \subset \ss_{m-1,n},$ then $Y$ must fix a partial flag $\bar F \in 
\ff^{(1^{m-2}(n-1)^{m-2})}$ and act diagonally (by the components of $\bar \tau$) on the 
successive quotients. The corresponding matrix $Y \in \ss_{m,n}$ with $Y_{m1} =0$ fixes 
the partial flag $F$ obtained from $\bar F$ by taking direct sum with an arbitrary flag 
of type $(1, n-1)$ in $\cc^n \times \{0\}^{(m-1)n} \subset \ker(Y).$ Using 
Lemma~\ref{reorder}, we get that $Y$ is in $\ymnt.$ Conversely, starting with $F$ fixed 
by $Y \in \ymnt$ we can construct a $\bar F \in \ff^{(1^{m-2}(n-1)^{m-2})}$ fixed by 
$\bar Y$ by intersecting everything with the subspace $\{0\}^n \times \cc^{(m-1)n} 
\subset \cc^n. \hfill \fin$ \medskip

\subsection {The case $m=1$} \label {sec:mone} When $m=1,$ the slice $\ss_{1,n}$ is the whole 
$\sl(n,\cc).$ The subvariety $\slan$ consists of traceless matrices having 
an eigenspace of dimension at least $n-1,$ cf. Example~\ref{ex0}.  The spaces $\yy_{1,n, 
\tau}$ are the fibers of the map
\begin {equation}
\label {m1}
\chin : \slan \to \cc,
\end {equation}
which takes a matrix into its eigenvalue of multiplicity $\geq n-1.$ This map will play an 
important role in our paper, being part of the local model for constructing vanishing 
projective spaces in Section~\ref{sec:vvo}.

\subsection {A thick and a thin eigenvalue coincide} \label{sec:thickthin} As mentioned in 
Section~\ref{sec:props}, the fiber $\ymnt$ is smooth when $\lambda_i$ and $\mu_j$ are all 
distinct. In this subsection we study the structure and form of the singularities in a 
fiber $\ymnt,$ where the $\lambda_i$ and $\mu_j$ that appear in $\tau$ satisfy $\lambda_1 
= \mu_1 = \lambda,$ but are otherwise distinct.

Lemma~\ref{structure} tells us the structure of the fiber of $\chi^{\pi}: \gg^{\pi} \to 
\hh^{\pi}/W^{\pi}$ over $\tau.$ The partition $\sigma$ associated to a semisimple element 
in that fiber is $\sigma=(1^{m-1}(n-1)^{m-1}n).$ There is a unique breaking $b$ of $\pi = 
(1^m(n-1)^m)$ according to $\sigma,$ cf. Example~\ref{ex2}. It follows that there are two 
adjoint orbits in the fiber 
$(\chi^{\pi})^{-1}(\tau):$ a regular orbit $\oo^{\reg}$ consisting of matrices with a 
Jordan block of size two and $n-2$ blocks of size one for the eigenvalue $\lambda,$ and a 
subregular orbit $\oo^{\sub}$ of matrices having $n$ independent $\lambda$-eigenvectors. 
The orbit $\oo^{\reg}$ is open and dense in $(\chi^{\pi})^{-1}(\tau),$ while $\oo^{\sub}$ 
is closed. By Lemma~\ref{locmodel}, the points of $\oo^{\reg}$ are regular for 
$\chi^{\pi},$ and then from Proposition~\ref{big}(ii) the points of $\oo^{\reg} \cap 
\ss_{m,n}$ are regular for the restriction of $\chi^{\pi}$ to the slice $\ss_{m,n}.$

\begin {remark}
When $\tau$ is as above, we denote by $\ccc_{m,n, \tau} = \oo^{\sub} \cap \ss_{m,n}$ the 
singular set of $\ymnt.$ This 
is a regular fiber of the restriction of $\chi^{\sigma}$ to $\ss_{m,n},$ with $\sigma = 
(1^{m-1}(n-1)^{m-1}n).$ The union of all $\ccc_{m,n, \tau}$ is $\ccc_{m,n} = \ss_{m,n} \cap 
\gg^{\sigma, \reg}.$ By Proposition~\ref{luna}, we have a differentiable fiber bundle: 
\begin {equation}
\label {cmnt}
\ccc_{m,n} \to \hh^{\sigma, reg}/W^{\sigma} = \conf^0_{\sigma}.
\end {equation}
\end {remark} 

The structure of the map $\chi^{\pi}|_{\ss_{m,n}}$ near the singular stratum $\oo^{\sub} 
\cap \ss_{m,n}$ is described by the lemma below, which basically states that the local model 
is a trivially fibered version of the $m=1$ case, in other words of the map (\ref{m1}).

\begin {lemma}
\label {thickthin}
Let $D \subset \hh^{\pi}/W^{\pi}$ be a disk corresponding to $\lambdav = ( \lambda - 
(n-1)z, \lambda_2, \dots, \lambda_m)$ and $\muv =(\lambda + z, \mu_2, \dots, \mu_m)$ 
with $z \in \cc$ small. Then there is a neighborhood of $\oo^{\sub} \cap \ss_{m,n}$ inside 
$(\chi^{\pi})^{-1}(D) \cap \ss_{m,n},$ and an isomorphism of that with a neighborhood of 
$(\oo^{\sub} \cap \ss_{m,n}) \times \{0\}$ inside  $(\oo^{\sub} \cap \ss_{m,n}) \times 
\slan.$ (Here $0 \in \sl(n, \cc)$ is the zero matrix.) The isomorphism 
fits into a commutative diagram:
$$ \begin {CD}
(\chi^{\pi})^{-1}(D) \cap \ss_{m,n} @>{\text{local } \cong \ \text{near} \ \oo^{\sub} 
\cap \ss_{m,n}}>> 
(\oo^{\sub} \cap \ss_{m,n}) \times \slan \\
@V{\chi^{\pi}}VV   @V{\chin}VV \\
D @>{\quad\quad\quad z \quad\quad\quad}>> \cc
\end {CD} $$
\end {lemma}

\noindent \textbf{Proof. } Let us first look at a neighborhood 
$\oo^{\sub}$ inside $(\chi^{\pi})^{-1}(D),$ without restricting to $\ss_{m,n}.$ From 
Section~\ref{sec:semis} we know that $\chi^{\pi}$ looks locally like a fibered version 
of the product of $\chin,$ $n-1$ copies of $\chi^{(n-1)},$ and $n-1$ copies of 
$\chi^{(1)},$ see (\ref{cairo2}). Since $\chi^{(n-1)}$ and $\chi^{(1)}$ correspond to 
trivial partitions, they are all just the identity map for a point. Thus we know 
that $\chi^{\pi}$ looks like a fibered version of $\chin$ near $\oo^{\sub}.$ 

We need to show that when we restrict to $\ss_{m,n},$ this fibered structure is preserved 
and, moreover, that the normal data along $\oo^{\sub} \cap \ss_{m,n}$ is trivial. At every 
point $Y \in \oo^{\sub} \cap \ss_{m,n},$ we choose a subspace $R_Y \subset T_Y\ss_{m,n}$ 
which is complementary to $T_Y(\oo^{\sub} \cap \ss_{m,n})$ and depends holomorphically on 
$Y.$ This is possible because $\oo^{\sub} \cap \ss_{m,n}$ is affine, and therefore the 
relevant $Ext^1$ obstruction group is zero. The spaces $(Y + R_Y) \cap \gg^{\pi}$ form a 
local tubular neighborhood of $\oo^{\sub} \cap \ss_{m,n}$ inside $\gg^{\pi} \cap \ss_{m,n}.$ 
On the other hand, we also have the canonical slice $\ss^{ss} = \ss^{ss}_Y$ to $\oo^{\sub},$ 
and a canonical local isomorphism (\ref{caniso}) between $(Y+R_Y) \cap \gg^{\pi}$ and 
$\ss^{ss}_Y \cap \gg^{\pi},$ which only moves points inside their adjoint orbits. Hence the 
restriction of $\chi^{\pi}$ to $\ss^{ss}_Y \cap \gg^{\pi}$ serves as a local model for 
$\chi^{\pi}|_{\ss_{m,n}}$ near $\oo^{\sub} \cap \ss_{m,n}.$

It remains to check triviality of the normal data. By Lemma~\ref{hello}, the intersection 
$\ss^{ss}_Y \cap \gg^{\pi}$ can be identified via translation with 
$$Z(\gg_Y) \times \sl(E_Y(\lambda))^{(1(n-1))} \times \prod_{i=2}^m 
\sl(E_Y(\lambda_i))^{(1)} \times \prod_{i=2}^m \sl(E_Y(\mu_i))^{(n-1)}.$$

Here by $E_Y(\alpha)$ we denoted the eigenspace of $Y$ with eigenvalue $\alpha.$ Note that 
each of the $\sl(E_Y(\lambda_i))^{(1)}$ and $ \sl(E_Y(\mu_i))^{(n-1)}$ is just a point, 
while $Z(\gg_Y)$ can be identified with $\cc^{2m-2}$ in a canonical way. Finally, by 
Lemma~\ref{newl} we have a preferred isomorphism of $E_Y(\lambda)$ with $\cc^n,$ and thus 
of $\sl(E_Y(\lambda))^{(1(n-1))}$ with $\slan.$ This completes the proof of 
triviality for the normal data. It is easy to see that all the 
isomorphisms which we used depended holomorphically on $Y.$ $\hfill \fin$
\medskip

\begin {remark} \label{remcom} Consider the partition $\pi^+=(1^{m+n-1}(n-1)^{m-1}).$ Since 
$\pi^+ \preccurlyeq \pi,$ we have that $$(\chi^{\pi})^{-1}(D) \cap \ss_{m,n}\ \subset \ 
\gg^{\pi} \cap \ss_{m,n} \ \subset \ \gg^{\pi^+} \cap \ss_{m,n}.$$ We can proceed as in 
the proof of Lemma~\ref{thickthin} and 
study the structure of $\gg^{\pi^+} \cap \ss_{m,n}$ near $(\oo^{\sub} \cap \ss_{m,n}).$ There 
are now several breakings of $\pi^+$ according to $\sigma = (1^{m-1}(n-1)^{m-1}n),$ cf. 
Example~\ref{ex3}. One of them is $b^+$ given by $n=1+\cdots +1, \ (n-1), \cdots, (n-1), 1, 
\cdots, 1.$ 
Lemma~\ref{hello} says that a neighborhood of $(\oo^{\sub} \cap \ss_{m,n})$ in $\gg^{\pi^+} 
\cap \ss_{m,n}$ has several components. We denote the one corresponding to $b^+$ by 
$(\gg^{\pi^+} \cap \ss_{m,n})_{b^+}.$ It is isomorphic to $(\oo^{\sub} \cap \ss_{m,n}) \times 
\sl(n, \cc).$ (The normal data is trivial just as in Lemma~\ref{thickthin}.) Furthermore, 
this local isomorphism is compatible with the one in Lemma~\ref{thickthin}, in the sense that 
there is a commutative diagram 
$$ \begin {CD} 
\gg^{\pi} \cap \ss_{m,n} @>{\text{local } \cong }>> (\oo^{\sub} \cap \ss_{m,n}) 
\times \slan \\ @VVV @VVV \\ (\gg^{\pi^+} \cap \ss_{m,n})_{b^+} @>{\text{local 
} \cong }>> (\oo^{\sub} \cap \ss_{m,n}) \times \sl(n, \cc), 
\end {CD} $$ 
where the vertical maps are inclusions. \end {remark}

We can generalize Lemma~\ref{thickthin} to the case when several pairs of one thick and one 
thin eigenvalue come together. Let $\tau = (\lambdav, \muv)$ be a point in $\hh^{\pi}/W^{\pi}$ 
such that $\lambda_1 = \mu_1, \lambda_2 = \mu_2, \dots,  \lambda_k = \mu_k,$ and these 
are the only coincidences. There are $2^k$ orbits in the fiber of $\chi^{\pi}$ 
corresponding to which of the $k$ relevant Jordan blocks are semisimple. The smallest 
orbit $\oo^{\min}$ where all the $k$ blocks are semisimple is closed in 
$(\chi^{\pi})^{-1}(\tau).$ An adaptation of the arguments in the proof of 
Lemma~\ref{thickthin} gives:

\begin {lemma}
\label {multithickthin}
Let $P \subset \hh^{\pi}/W^{\pi}$ be a polydisk corresponding to $\lambdav = ( \lambda -
(n-1)z_1, \dots, \lambda_k - (n-1)z_k, \lambda_{k+1}, \dots, \lambda_m)$ and $\muv =(\lambda + 
z_1, \dots, \lambda_k+z_k, \mu_{k+1}, \dots, \mu_m)$
with the $z_k$'s small. Then there is a neighborhood of $\oo^{\min} \cap \ss_{m,n}$ 
inside $(\chi^{\pi})^{-1}(P) \cap \ss_{m,n},$ and an isomorphism of that with a 
neighborhood of
$(\oo^{\min} \cap \ss_{m,n}) \times \{0\}^k$ inside  $(\oo^{\min} \cap \ss_{m,n}) \times
\bigl( \slan \bigr)^k.$ The isomorphism fits into
a commutative diagram:
$$ \begin {CD}
(\chi^{\pi})^{-1}(P) \cap \ss_{m,n} @>{\text{local } \cong \ \text{near} \ \oo^{\min}
\cap \ss_{m,n}}>>
(\oo^{\min} \cap \ss_{m,n}) \times \bigl( \slan \bigr)^k \\
@V{\chi^{\pi}}VV   @V{(\chin, \dots, \chin)}VV \\
P @>{\quad\quad(z_1, \dots, z_k)\quad\quad}>> \cc^k
\end {CD} $$
\end {lemma}

\subsection {Two thin eigenvalues coincide} Although not necessary for the rest of the 
paper, it is instructive to consider the case when $\tau = (\lambdav, \muv)$ has $\lambda_1 
= \lambda_2 = \lambda,$ and all the other $\lambda_i$'s and $\mu_j$'s are distinct from 
each other and from $\lambda.$ Let us assume that $n > 3.$ The corresponding partition 
$\sigma$ is now $(1^{m-2}2(n-1)^m),$ and there is again a unique breaking of $\pi$ 
according to $\sigma,$ cf. Example~\ref{ex4}. We get two orbits $\oo^{\reg}$ and 
$\oo^{\sub}$ in the fiber, just as in Section~\ref{sec:thickthin}. The orbit $\oo^{\sub}$ 
is closed, and the local model around it is the map $\chi^{(1^2)} : \sl(2,\cc) \to \cc,$ 
which is basically the $A_1$-singularity $\cc^3 \to \cc, \ (a,b,c) \to a^2 + b^2 + c^2.$ 
The same arguments as in the proof of Lemma~\ref{thickthin} show that the local model for 
$\chi^{\pi}|_{\ss_{m,n}}$ near $\oo^{\sub} \cap \ss_{m,n}$ is a fibered version of the 
$A_1$-singularity. However, note that we do not have a canonical identification of the 
$\lambda$-eigenspace with $\cc^2,$ and so we cannot claim the triviality of the normal 
data.

\subsection {One thick and two thin eigenvalues coincide} \label {sec:ttt} A slightly more 
involved situation is when $\tau = (\lambdav, \muv)$ has $\lambda_1 = \lambda_2 = \mu_1 = 
\lambda,$ and these are the only coincidences. The partition $\sigma$ is 
$(1^{m-2}(n-1)^{m-1}(n+1)),$ and there is still a unique breaking of $\pi$ according to 
$\sigma,$ namely $(n+1) = (n-1) + 1 +1$ and the rest trivial, cf. Example~\ref{ex5}. By 
Lemma~\ref{structure} the 
fiber $(\chi^{\pi})^{-1}(\tau)$ consists of three orbits: $\oo^{\reg}$ made of matrices 
with a Jordan block of size $3$ and $n-2$ blocks of size $1$ for the eigenvalue $\lambda;$ 
a subregular orbit $\oo^{\sub}$ with a Jordan block of size $2$ and $n-1$ blocks of size 
$1$ for $\lambda$; and a minimal semisimple orbit with $n+1$ independent 
$\lambda$-eigenvectors. However, according to Lemma~\ref{newl}, the minimal orbit does not 
intersect $\ss_{m,n}.$ Thus $\ymnt$ is only made of two strata: an open and dense part 
$\oo^{\reg} \cap \ss_{m,n},$ and a closed part $\oo^{\sub} \cap \ss_{m,n}.$

The discussion from Section~\ref{sec:gen} can be applied here to find a local model for the 
map $\chi^{\pi}$ near $\oo^{\sub}.$ One can show that the same model is valid for the 
restriction of $\chi^{\pi}$ to $\ss_{m,n},$ along the lines of the proof of 
Lemma~\ref{thickthin}, and then study the normal data in the fibers (which turns out to be 
nontrivial). To be precise, for every $Y \in \oo^{\sub} \cap \ss_{m,n},$ we denote by $E_Y$ 
the kernel of $(\lambda I - Y)^2,$ by $F_{1,Y}$ the image of $E_Y$ under the map $\lambda I 
- Y,$ and by $F_{2,Y}$ the kernel of $\lambda I - Y.$ We let $\ee, \ff_1, \ff_2$ be the 
complex vector bundles over $\oo^{\sub} \cap \ss_{m,n}$ whose fibers over $Y$ are $E_Y, 
F_{1,Y}, F_{2,Y},$ respectively. Since $Y$ is affine, we can choose complements $T_{1,Y}$ 
of $F_{1,Y}$ in $F_{2, Y}$ and $T_{2,Y}$ of $F_{2,Y}$ in $E_Y$ that depend holomorphically 
on $Y.$ The map $\lambda I - Y$ induces a linear isomorphism from $T_{2,Y}$ to $F_{1, Y}.$ 
This gives a linear quadruple $T_Y$ in the sense of Remark~\ref{act}. Since the data $T_Y$ 
depends holomorphically on $Y,$ we call the whole thing a {\it holomorphic quadruple 
bundle} over $\oo^{\sub} \cap \ss_{m,n},$ and denote it by $\tt.$ Using Remark~\ref{act}, 
we have an associated space $\xx(T_Y)$ for each $Y;$ the corresponding fiber bundle over 
$\oo^{\sub} \cap \ss_{m,n}$ is denoted $\xx(\tt).$ The fibers of $\xx$ are isomorphic to 
the space $\xx_n$ which appears in (\ref{mapp}). There is also a map $q: \xx(\tt) \to 
\sym^2(\cc)$ which is given by (\ref{mapp}) in each fiber.

\begin {lemma}
\label {ttt}
Consider the bidisk $P \subset \hh^{\pi}/W^{\pi}$ corresponding to $\lambdav = ( 
\lambda + z_1, \lambda+z_2, \dots, \lambda_m)$ and $\muv =(\lambda + z_3, \mu_2, \dots, 
\mu_m)$ with $z_1, z_2, z_3 \in \cc$ small, $z_1 + z_2 + (n-1)z_3 = 0.$

Then there is a neighborhood of $\oo^{\sub} \cap \ss_{m,n}$ inside
$(\chi^{\pi})^{-1}(P) \cap \ss_{m,n},$ and an isomorphism of that with a neighborhood of
the zero-section inside  $\xx(\tt).$ The isomorphism fits 
into a commutative diagram:
$$ \begin {CD}
(\chi^{\pi})^{-1}(P) \cap \ss_{m,n} @>{\text{local } \cong \ \text{near} \ \oo^{\sub}
\cap \ss_{m,n}}>> \xx(\tt) \\
@V{\chi^{\pi}}VV   @VV{q=(q_1, q_2)}V \\
P @>{\quad\quad(z_1,z_2)\quad\quad}>> \sym^2(\cc).
\end {CD} $$
\end {lemma}

In short, this says that the map $\chi^{\pi}|_{\ss_{m,n}}$ looks locally near $\oo^{\sub} 
\cap \ss_{m,n}$ like a fibered version of the map (\ref{mapp}). The possible nontriviality
of the normal data is encoded in the structure of the fiber bundle $\xx(\tt).$ In 
particular, the bundles $\ee$ and $\ff_1$ can be nontrivial. On the other hand, note that 
by Lemma~\ref{newl}, the bundle $\ff_2$ is always trivial.

There is also an analogue of Remark~\ref{remcom}:
\begin {remark}
\label {remcom2}
One breaking of $\pi^+=(1^{m+n-1}(n-1)^{m-1})$ according to $\sigma=(1^{m-2}(n-1)^{m-1}(n+1))$
is $\beta^+$ given by $(n+1)=1+\cdots+1, \ (n-1), \cdots, (n-1), 1, \cdots, 1.$ Take a 
neighborhood of $\oo^{\sub} \cap \ss_{m,n}$ in $\gg^{\pi^+} \cap \ss_{m,n}$ as in 
Lemma~\ref{hello}. Denote by  $(\gg^{\pi^+} \cap \ss_{m,n})_{\beta^+}$ the part corresponding 
to $\beta^+.$ Then there is a commutative diagram made of inclusions and local isomorphisms:
$$ \begin {CD}
\gg^{\pi} \cap \ss_{m,n} @>{\text{local } \cong}>>  \xx(\tt)  \\ 
@VVV @VVV \\ (\gg^{\pi^+} \cap \ss_{m,n})_{\beta^+} 
@>{\text{local } \cong }>>  \sl(\ee).
\end {CD} $$
\end {remark}

\section {Parallel transport and various vanishing objects}
\label {sec:vvo}

\subsection {Parallel transport} \label{sec:parallel}

In this section we review the definition of rescaled parallel transport in a Stein fibration, 
following Seidel and Smith \cite[Section 4(A)]{SS}.

Let $p: Y \to T$ be a holomorphic map between complex manifolds, which is also a submersion 
with fibers $Y_t.$ Let $\gamma: [0,1] \to T$ be a path on the base. The parallel transport 
vector field $H_{\gamma}$ on the pullback $\gamma^*Y \to [0,1]$ consists of the sections of 
$TY|Y_{\gamma(s)}$ which project to $\gamma'(s)$ and are orthogonal to the vertical tangent 
bundle $T(Y_{\gamma(s)}).$ Thus,
\begin {equation}
\label {hgam}
H_{\gamma} = \frac{\nabla p}{\| \nabla p \|^2 } \gamma'(s).
\end {equation}

In some cases, for example if $p$ is proper or if we have 
suitable 
estimates on $H_{\gamma},$ integrating $H_{\gamma}$ yields a symplectic isomorphism $h_{\gamma} 
: Y_{\gamma(0)} \to Y_{\gamma(1)},$ called \emph{(naive) parallel transport.} In general, 
however, the lack of compactness means that integral lines may not exist everywhere, so 
that $h_{\gamma}$ is only defined on compact subsets $P \subset Y_{\gamma(0)}$ for a short time 
$\epsilon > 0$ which depends on $P.$

To fix this problem, we use \emph {rescaled parallel transport}. Assume that there exists a 
function $\psi: Y \to \rr$ with the following properties:
\begin {eqnarray}
\label {cond1}
&\bullet & \psi  \text{ is proper and bounded below. } \\
\label {cond2}
&\bullet & -dd^c \psi > 0, \text{ so that } \Omega=-dd^c\psi \text{ is a K\"ahler form on }Y. 
\\
\label{cond3}
& \bullet & \text{ Outside a compact set of }Y, \ \text{we have} \ \|\nabla \psi\| \leq
\rho \psi \text{ for some }\rho > 0; \\
\label {cond4}
& \bullet & \text{ The fiberwise critical set } \{y \in Y : d\psi|_{\ker(Dp)} = 0\}  \text{ 
maps properly to }T.
\end {eqnarray}

The Liouville vector field $Z_t = \nabla (\psi|_{Y_t})$ on $Y_t$ is well-defined for all times 
because of condition (\ref{cond3}). Given the path $\gamma,$ condition (\ref{cond4}) implies 
that there exists $c > 0$ such that the critical values of $\psi$ on the fibers $Y_{\gamma(s)}$ 
all lie in $[0,c).$ Using this and (\ref{cond1}) we can find $\sigma > 0$ such that the 
integral lines of $ \tilde H_{\gamma} = H_{\gamma} - \sigma Z_{\gamma(s)}$ stay inside 
$\psi^{-1}([0,c]).$ Integrating the flow $\tilde H_{\gamma}$ and then composing with the time 
$\sigma$ map of the Liouville flow on $Y_{\gamma(1)}$ yields a symplectic embedding
$$ h_{\gamma}^{\resc}: Y_{\gamma(0)} \cap \psi^{-1}([0,c]) \longrightarrow Y_{\gamma(1)}.$$

This is called rescaled parallel transport. It depends on $\sigma$ only up to isotopy in the 
class of symplectic embeddings. Replacing $c$ by a larger value yields a map defined on a 
bigger set, whose restriction to the smaller set is isotopic to the original one. As a 
consequence, the image $h_{\gamma}^{\resc}(L)$ is well-defined (up to Lagrangian isotopy) for 
any compact Lagrangian submanifold $L \subset Y_{\gamma(0)}.$

\subsection {Vanishing cycles for singular metrics.} \label{sec:vsm} Let $\Omega$ be an 
arbitrary 
K\"ahler form on $Y=\cc^n,$ and denote by $g$ the corresponding K\"ahler metric. Consider the 
projection
$$ p: \cc^n \to \cc, \ p(z_1, \dots, z_n) = z_1^2 + \dots + z_n^2,$$
and equip the fibers with the induced metrics and forms.

The real part $re(p)$ is a Morse function with a critical point at the origin. The 
\textit{stable manifold} $W$ is defined as the set of points $y \in \cc^n$ such that the flow 
line of $-\nabla re(p)$ starting at $y$ exists for all times $s \geq 0,$ and converges to 
zero as $s \to \infty.$ 

Since $\Omega$ is K\"ahler, the negative gradient flow of $re(p)$ is the same as the 
Hamiltonian vector field of the imaginary part $im(p).$ This preserves $im(p),$ hence $W$ 
lies in the preimage of $\rr_{\geq 0}.$ The intersection $C_t$ of $W$ with a fiber 
$Y_t=p^{-1}(t)$ for $t >0$ is called a \textit{vanishing cycle}. It is well-known that for 
$t$ small $C_t$ is a smooth manifold and, in fact, is diffeomorphic to $S^{n-1}.$ (See 
\cite{J}, for example.) An alternate definition is as follows. Consider the path $\gamma: 
[0,1] \to \cc, \ \gamma(s)=(1-s)t.$ Since $-\nabla re(p)$ is proportional to the naive 
parallel transport $H_{\gamma},$ we have \begin {equation} \label {altvc} C_t = \{ y \in Y_t 
\ | \ h_{\gamma|[0,s]} \text{ is defined near }y \text{ for } s <1, \text{ and } 
h_{\gamma|[0,s]}(y) \to 0 \text{ as } s \to 1 \}. \end {equation}

Because the gradient flow preserves $\Omega,$ by taking the limit $s \to 1$ we find that 
$\Omega$ vanishes on $TC_t.$ In other words, $C_t \subset Y_t$ is Lagrangian.

We need a generalization of this construction in the case when the K\"ahler form 
$\Omega$ is degenerate at the origin. In that situation the gradient flow $-\nabla re(p)$ can 
take a nonzero point into zero in finite time. However, we can still define a vanishing cycle 
$C_t$ using parallel transport maps. 

\begin {proposition}
\label {vansing}
Let $\Omega$ be a smooth two-form on $\cc^n$ that is real analytic in a neighborhood 
of zero and K\"ahler on $\cc^{n} - \{ 0 \}.$ For $t > 0,$ define $C_t$ by (\ref{altvc}).
Then, for $t > 0$ sufficiently small, $C_t$ is a Lagrangian $(n-1)$-sphere in $Y_t.$
\end {proposition}

(It is conceivable that the real analyticity consition is unnecessary. However, our proof 
is fundamentally different from the one in the nondegenerate case, and makes essential use of 
real analyticity.)

\medskip

\noindent \textbf {Proof of Proposition ~\ref{vansing}. } For $t \in \rr,$ we will denote by 
$Y_{t+i\rr}$ the preimage of $t$ under the map $re(p): Y \to \rr.$ This is the union of 
$Y_{t+i\tau}$ for all $\tau \in \rr.$ Note that $Y_{t+i\rr}$ is diffeomorphic to $S^{n-1} 
\times \rr^n$ for all $t \neq 0.$

Consider the projection from $\cc^n$ to its imaginary part 
$$\pim: \cc^n \to \rr^n, \ \pim(z_1, \dots, z_n) = (im(z_1), \dots, im(z_n)).$$ 

Denote by $\iota$ the inclusion of $Y_{i\rr}$ into $\cc^n.$ Take the preimage $V=(\pim 
\circ \iota)^{-1}(B) \subset Y_{i\rr}$ of a ball $B^n=B^n(\epsilon) \subset \rr^n$ centered 
at the 
origin. Given a point $y \in V -\{0\}$ with $p(y)=i\tau,$ we can consider the translate 
$\gamma_{\tau}$ of $\gamma$ by $i \tau,$ as a path in $\cc.$ Reverse parallel transport along 
the corresponding $\gamma_{\tau}$ can be applied to all points in $V-\{0\}$ for at least some 
small time $t > 0.$ (Up to a reparametrization of time, this is equivalent to going along the 
forward gradient flow of $re(p).$) The image of this reverse transport, together with the 
vanishing cycle $C=C_t$, forms an open set $U \subset Y_{t+i\rr}.$ We can define $h: U \to V$ 
by setting $h(y)=0$ for $y \in C,$ and $h(y)=h_{\gamma_{\tau}}(y)$ for $y\in V-C,$ with 
$\tau$ being the imaginary part of $p(y).$

If our parameters $\epsilon$ and $\tau$ were sufficiently small, the whole parallel transport 
of $U$ along the $\gamma_{\tau}$'s happens in the neighborhood where $\Omega$ is real 
analytic. The solution of an ODE with real analytic coefficients and initial data is a 
real analytic function by the Cauchy-Kovalevsky theorem. Therefore the map $\iota \circ h : 
U \to \cc^n$ is real analytic, and so must be 
$$f=\pim \circ h: U \to B^n.$$

The restriction of $h$ to $U - C$ is a diffeomorphism onto $V-\{ 0 \}.$ It follows 
that $f: U \to B^n$ is an $S^{n-1}$-fibration over $B^n - \{0\},$ outside of the vanishing 
cycle $C.$ We aim to show that $C=f^{-1}(0)$ is also an $(n-1)$-sphere.

Note that $C$ is an analytic subvariety of the real analytic manifold $Y_{t+i\rr} = S^{n-1} 
\times \rr^n.$ Let us recall a few facts about real analytic varieties, cf. \cite[Section 
6.3]{KP}, \cite{Lo}. Every analytic variety $A$ admits a Whitney stratification and, in 
particular, it has a top-dimensional stratum of dimension $d=\dim A.$ Every point $y \in A$ 
admits a Zariski tangent space $T_yA,$ and (assuming $A$ is connected) we have $\dim (T_yA) 
\geq \dim A$ for all $y,$ with equality if and only if $y$ is a regular point, i.e. part of 
the top-dimensional stratum.

In our case, since parallel transport preserves $\Omega,$ it follows that $\Omega$ must 
vanish on $T_yC$ for all $y \in C,$ cf. the remark following (\ref{altvc}). Since $C 
\subset Y_t$ and the 
restriction of $\Omega$ to the complex submanifold $Y_t \subset Y$ is nondegenerate, 
this implies that $\dim (T_yC) \leq n-1$ for all $y \in C.$ In particular, the dimension 
of $C$ is at most $n-1.$

On the other hand, each fiber $f^{-1}(x)$ for $x \neq 0$ represents a generator of 
$H_{n-1}(S^{n-1} \times \rr^n).$ Hence $f^{-1}(x)$ must intersect each $\{z\} \times \rr^n 
\subset S^{n-1} \times \rr^n$ nontrivially. Taking the limit $x \to 0,$ we find that the 
projection of $f^{-1}(0) = C \subset S^{n-1} \times \rr^n$ to the $S^{n-1}$ factor is 
surjective as well. As a consequence, $C$ must have dimension exactly $n-1.$

For $n >1,$ we also claim that $C$ is connected. If it were disconnected, then a small 
neighborhood $f^{-1}(B^n(\epsilon'))$ for $\epsilon' < \epsilon$ would also be 
disconnected, and the same would hold for its boundary $f^{-1}(S^{n-1}(\epsilon')).$ 
However, this boundary is diffeomorphic to $S^{n-1} \times S^{n-1}$, and we arrive at a 
contradiction.

Since $C$ is a connected analytic variety of dimension $n-1$ and $\dim (T_yC) \leq n-1$ at any 
$y \in C,$ we deduce that all of its points are regular. In other words, $C$ is an 
$(n-1)$-dimensional real analytic manifold. As we already noted, $\Omega$ vanishes on its 
tangent space.

To show that $C$ is a sphere, we interpolate between $\Omega^{(0)}=\Omega$ and a form 
$\Omega^{(1)}$ that is nondegenerate over all of $\cc^n.$ For example, if $N$ is the null 
space of $\Omega$ at the origin and $N^{\perp}$ is its orthogonal complement in the 
standard Hermitian metric, we can choose $\Omega^{(1)}$ to be constant over all of $\cc^n$ 
(with respect to the standard holomorphic coordinates), identical to $\Omega$ on 
$N^{\perp}$ (at the origin), and nondegenerate on $N.$ Then, all the forms in the family 
$\Omega^{(r)} = (1-r)\Omega^{(0)} + r\Omega^{(1)}, \ r\in (0,1],$ are real analytic and 
K\"ahler in a neighborhood of zero. For $t$ sufficiently small, we can thus relate the 
vanishing cycles $C=C_t^{(0)}$ to $C_t^{(1)}$ by a smooth family of manifolds. Note that 
$\Omega^{(1)}$ is nondegenerate, hence we already know that $C_t^{(1)}$ is a sphere. Since 
all the manifolds in a smooth family are diffeomorphic, we conclude that $C$ is also a 
sphere.

For $n=1,$ a similar argument shows that $C$ is a disjoint union of exactly two points. 
$\hfill \fin$

\begin {remark}
\label {multsc}
By multiplying $p$ with some scalar in $S^1,$ we can define vanishing cycles $C_t \subset 
Y_t$ for all sufficiently small $t \in \cc^*.$
\end {remark}

This discussion can be generalized to relative vanishing cycles. Take any complex manifold 
$X$ and consider the projection
$$ p: Y=X \times \cc^n \to \cc, \ p(x,z_1, \dots, z_n) = z_1^2 + \dots + z_n^2.$$

Equip $Y$ with any K\"ahler form $\Omega$, possibly degenerate on $X \times \{ 0 \}^n.$ The 
function $re(p)$ is now Morse-Bott, and its critical point set can be identified with $X.$ Let 
$K \subset X$ be a compact Lagrangian submanifold. For $t \in \cc^*,$ the \textit{relative 
vanishing cycle} $C_t$ associated to $K$ is defined as the set of points $y \in Y_t=p^{-1}(t)$ 
which are taken into $K$ by parallel transport along the path $\gamma: [0,1] \to \cc, \ 
\gamma(s)=(1-s)t,$ as $s \to 1.$ The same arguments as in the proof of 
Proposition~\ref{vansing} show:

\begin {proposition}
\label {relvs}
Assume that $\Omega$ is real analytic in a neighborhood of $X.$ Then for sufficiently small 
$t \in \cc^*, \ C_t$ is a Lagrangian submanifold of $Y_t$ diffeomorphic to $K \times 
S^{n-1}.$ 
\end {proposition}

The case when $\Omega$ is nondegenerate appears in \cite[Lemma 26]{SS}.

\subsection {Vanishing projective spaces} \label {sec:van}

We now describe a construction similar to that of vanishing cycles, but where instead of a 
sphere we obtain a complex projective space $\cpn$. Our main interest lies in a fibered 
version of this (treated in the next section). However, we decided to present the simpler 
situation first in order to make clear all the ideas involved.

The construction takes place in the space $Z=\slan$ from 
Section~\ref{sec:mone}, for $n > 2.$ Recall that
$$Z = \{ A \in \sl(n, \cc) \ | \text{ there exists } t \in \cc \text{ with } \dim \ker (A- tI) 
\geq n-1 \},$$
and the map $\chin$ takes a matrix $A$ to the corresponding $t \in \cc.$ 

The following alternate description is also helpful:

\begin {lemma}
\label {git}
Consider the linear action of $\cc^*$ on $\cc^{2n}$ with weights $1$ and $-1,$ each with 
multiplicty $n.$ Then the space $Z$ is the GIT quotient of $\cc^{2n}$ by this action.
\end {lemma}

\noindent \textbf {Proof. } Write $\cc^{2n}$ as the space of pairs $(v, w),$ with $v, w \in 
\cc^n,$ and the $\cc^*$ action given by 
\begin {equation} 
\label {cstar}
(v, w) \to (\zeta v, \zeta^{-1}w)
\end {equation}
 for $\zeta \in 
\cc^*.$ All orbits corresponding to $v, w \neq 0$ are stable. There is one semistable orbit 
given by $v=w=0,$ and two unstable orbits corresponding to only one of $v$ and $w$ being 
zero.

We define an algebraic map
$$ f: \cc^{2n} \to Z, \ f(v,w) = v^Tw - \frac{1}{n}(v\cdot w)I.$$

This is $\cc^*$-invariant and therefore factors through the GIT quotient $\cc^{2n}/\cc^*.$ It 
is easy to see that the induced map $\cc^{2n}/\cc^* \to Z$ is bijective. $\hfill \fin$

\medskip

Note that $Z$ has an isolated singularity at the origin, and the map $\chin$ is a 
submersion everywhere except over zero. Given a K\"ahler metric $\Omega$ on $Z - \{0\}$, we 
can study parallel transport with respect to $\chin.$ For $t \in \cc^*,$ we take 
the path $\gamma: [0,1] \to \cc, \ \gamma(s) = (1-s)t$ as in the previous section, and define 
\begin {equation}
\label {vcps}
 L_t = \{ y \in Z_t
\ | \ h_{\gamma|[0,s]} \text{ is defined near }y \text{ for all } s <1, \ 
h_{\gamma|[0,s]}(y) \to 0 \text{ as } s \to 1 \}. 
\end {equation}

Before discussing the general case, let us work with a specific example. Take $\Omega = 
\Omega_{st}$ to be the restriction of the standard form on $\cc^{n^2}$ under the inclusion $Z 
\hookrightarrow \gl(n, \cc) \cong \cc^{n^2}.$ For each $t \in \cc^*,$ consider the diagonal 
matrix $E_t=diag(t,t,\dots,t, (1-n)t) \in Z_t,$ and let
\begin {equation}
\label {ua}
\ua_t = \{ UE_tU^{-1} \ | \ U \in U(n) \}.
\end {equation}

Note that $\ua_t$ is diffeomorphic to $U(n)/(U(n-1) \times U(1)) \cong \cpn.$

\begin {lemma}
\label {ost}
The vanishing space $L_t \subset Z_t$ for $\Omega_{st}$ contains $\ua_t.$
\end {lemma}

\noindent \textbf {Proof. } We claim that parallel transport along $\gamma|[0,s]$ takes 
$UE_tU^{-1}$ into $UE_{t(1-s)}U^{-1}$ for $s<1$ and any $U \in U(n).$ From this it will follow 
that in the limit $s \to 1,$ the space $\ua_t$ is indeed sent to the origin.

It suffices to verify our claim infinitesimally, i.e. to show that the parallel transport
vector field at $UE_tU^{-1}$ is proportional to $UE_1U^{-1}.$ Since the K\"ahler metric is 
$U(n)$-invariant, it also suffices to check this at $U=I.$ There the statement is that $E_1$ 
is perpendicular to the tangent space to $Z_t$ at $E_t.$ The tangent space $T_{E_t}Z_t$ 
consists of the commutators $[E_t, B]$ for $B \in \gl(n,\cc).$ A simple computation shows that 
$[E_t, B]$ is always perpendicular to $E_1$ in $\cc^{n^2}. \hfill \fin$

\medskip

We now turn our attention to a more general class of K\"ahler metrics:

\begin {definition}
\label {eff}
Let $i$ be the standard inclusion of $Z=\slan$ in $\sl(n, \cc).$ A K\"ahler 
form 
$\Omega$ on $Z - \{0\}$ is called effective if it is of the form $i^*\omega,$ where $\omega$
is a K\"ahler form on $\sl(n, \cc)$ that is real analytic in a neighborhood of zero.  
\end {definition}

\begin {lemma}
\label {commens}
Let $\Omega$ be an effective K\"ahler form on $Z - \{0\},$ and $\Omega_{st}$ the standard 
form from Lemma~\ref{ost}. Then the corresponding K\"ahler metrics $g$ and $g_{st}$ are 
equivalent in 
the following sense: given a relatively compact open set $V \subset Z,$ there exists a 
constant $c > 0$ such that $c^{-1}\cdot g_{st} \leq g \leq c\cdot g_{st}$ on $V - \{0\}.$ 
\end {lemma}

\noindent \textbf {Proof. } This is a direct consequence of the fact that both metrics are 
obtained by pull-back from the nonsingular space $\sl(n, \cc). \hfill \fin$

\begin {lemma}
\label {exist}
Let $\Omega$ be an effective K\"ahler form on $Z - \{0\},$ and define $L_t$ as in 
(\ref{vcps}). Then for any sufficiently small $t \in \cc^*, \ L_t$ is a smooth manifold 
diffeomorphic to $\cpn.$
\end {lemma}

\noindent \textbf {Proof. } Let $f: Y=\cc^{2n} \to Z$ be the GIT quotient map from 
Lemma~\ref{git}. Observe that $f^*\Omega$ is nondegenerate on the union of all stable 
$\cc^*$-orbits, except in the directions of those orbits. Consider the Hamiltonian $H: 
\cc^{2n} \to \rr, \ H(v,w) = |v|^2-|w|^2$ associated to the action of $S^1 \subset \cc^*.$ The 
restriction of the map $|H|^2: \cc^{2n} \to \cc$ to any nonzero (stable or unstable) 
$\cc^*$-orbit is strictly plurisubharmonic. Let $g_H$ be the symmetric two-tensor 
associated to $-dd^c|H|^2$ by the complex structure. An explicit calculation shows that 
for any neighborhood $N$ of zero in $\cc^n,$ we can find some small $\epsilon > 0$ such that
$f^*g_{st} + \epsilon g_H > 0$ on $N-\{0\}.$ Lemma~\ref{commens} implies that the same must be 
true for $f^*g$ instead of $f^*g_{st}$ (with possibly a different $\epsilon.$) It follows that
$$ \tilde \Omega = f^*\Omega - \epsilon dd^c |H|^2$$ 
is K\"ahler on $N - \{0\}.$ 

We can now apply Proposition~\ref{vansing} to our situation; see also Remark~\ref{multsc}. 
(The fact that $\tilde \Omega$ is not K\"ahler everywhere on $\cc^{2n}$ does not make a 
difference.) In place of $p: \cc^{2n} \to \cc$ we take the map $(v,w) \to (v \cdot w)/n$ which 
is the composition of $f$ with $\chin.$ This is the same as the map $p$ from 
Section~\ref{sec:vsm}, up to a linear change of coordinates.

We obtain vanishing cycles $C_t \subset Y_t$ diffeomorphic to $S^{2n-1},$ for any $t \in 
\cc^*$ small. Note that $\tilde \Omega$ is $S^1$-equivariant, hence the same is true 
for parallel transport in $Y=\cc^{2n}.$ If we think of $Z$ as the quotient 
$H^{-1}(0)/S^1,$ then parallel transport along $\gamma$ in $(Y, \tilde \Omega)$ corresponds to 
parallel transport in $(Z, \Omega).$ The vanishing cycle $C_t$ is preserved by the $S^1$ 
action, and its quotient is $L_t.$ Since the action is free outside zero, we find that $L_t$ 
must be a smooth, connected $(2n-2)$-dimensional manifold. 

In particular, in the case when $\Omega$ is the standard $\Omega_{st}$ used in 
Lemma~\ref{ost}, we find that $\ua_t \cong \cpn$ is the whole of $L_t.$ To see that 
$L_t$ is diffeomorphic to $\cpn$ for a general $\Omega,$ interpolate between $\Omega$ and 
$\Omega_{st}$ by a family of effective forms, and use the fact that all manifolds in a 
smooth family are diffeomorphic.  $\hfill \fin$

\medskip

We call $L_t$ a {\it vanishing projective space} in $\slan.$ As seen in the 
proof of Lemma~\ref{exist}, $L_t$ is the quotient of an ordinary vanishing cycle 
$S^{2n-1}$ by an $S^1$ action. It is also a Lagrangian submanifold of $Z_t.$

\subsection {Relative vanishing projective spaces.} \label {sec:relvan}

Let $X$ be a smooth complex manifold. With an eye towards Lemma~\ref{thickthin}, we 
look at the product $Z= X \times \slan.$ A similar discussion to that in the 
previous section applies here, using Proposition~\ref{relvs}. Let us identify 
$X$ with $X \times \{0\}.$ Consider the map $$ 
\pi: Z \to \cc, \ (x, A) \to \chin(A),$$ and denote by $Z_t$ its fibers. Let $i$ be 
the standard inclusion of $Z$ into $X \times \sl(n,\cc).$ By analogy with Definition~\ref{eff}, we 
say that a K\"ahler form $\Omega$ on $Z - X = X \times (\slan - \{ 0\})$ is {\it 
effective} if it is obtained by pull-back from a form on $X \times \sl(n, \cc)$ that is real 
analytic around $X \times \{ 0 \}.$ 

For the remaining of this section we will assume that $Z-X$ is endowed with an effective 
K\"ahler form $\Omega.$ Given a compact Lagrangian $K \subset X,$ we obtain a {\it relative 
vanishing projective space} $L_t$ in $Z_t$ (for small $t$), diffeomorphic to $K \times \cpn.$ 

We denote by $g$ be the K\"ahler metric corresponding to $\Omega,$ and by $g_{st}$ the 
product of an arbitrary metric on $X$ with the pull-back of standard metric on $\sl(n,\cc).$ 
The analogue of Lemma~\ref{commens} holds true, so that $g$ and $g_{st}$ are equivalent on 
the complement of $X$ in any relatively compact subset of $Z.$

Fix a relatively compact open set $W \subset X$ and a ball $B$ in $\slan$ around 
the origin (in the standard pulled-back metric, say). We assume that the form $\Omega$ is real 
analytic on $V=W \times B \subset Z.$ We are interested in estimating the parallel transport 
maps in $V$ obtained from the map $\pi$ and the form $\Omega.$

\begin {lemma}
\label {piano}
The quantity $ \| \nabla \pi \| $ is bounded below on $V - X.$
\end {lemma}

\noindent \textbf {Proof. } Since $g$ and $g_{st}$ are equivalent, it suffices to establish 
the result when the gradient is taken with respect to $g_{st}.$ In that case it turns out 
that $ \| \nabla \pi \| $ is bounded below globally. Indeed, recall from (\ref{trivia}) 
that every matrix $A \in \slan$ is of the form $S + tUE_1U^{-1},$ where $S$ 
is a matrix of rank at most one, $U$ is unitary, $E_1$ is the diagonal matrix $diag(1, 
\dots, 1, 1-n)$, and $t \in \cc.$ At a point $(x,A) \in Z$ we have $\pi(x,A) = t.$ Fix a 
representation of $A$ as $S+tUE_1U^{-1},$ and consider the one-dimensional complex subspace 
$\{x\} \times \{ S+\tau UE_1U^{-1} \ | \ \tau \in \cc \} \subset Z.$ The gradient of the 
restriction of $\pi$ to this subspace is $UE_1U^{-1}.$ This always has norm 
$(n^2-n)^{1/2},$ because the metric is $U(n)$-invariant. It follows that the actual 
gradient of $\pi$ on $Z$ (with respect to $g_{st}$) has norm at least $(n^2-n)^{1/2}. 
\hfill \fin$

\medskip

Now let $K \subset U \subset X$ be a compact Lagrangian submanifold, and denote by $\delta > 
0$ its distance from $\partial V$ in the metric $g.$ Let also $\nu$ be the lower bound on $\| 
\nabla \pi \|$ in $V$ which is given by the previous lemma.

\begin {lemma}
\label {wdw}
For any $t$ with $0 < |t| < \nu \delta/2,$ the relative vanishing cycle $L_t$ is well-defined 
and lies in $V.$ 
\end {lemma}

\noindent \textbf {Proof. } Without loss of generality we can assume that $t$ is real and 
positive. Think in terms of parallel transport along $\gamma(s)=s.$ By the previous lemma, the 
the horizontal vector field $H_{\gamma} = \nabla \pi/\| \nabla \pi \|^2$ satisfies
\begin {equation}
\label {hga}
\| H_{\gamma} \| \leq \nu^{-1}.
\end {equation}
Say we are given a flow line of $H_{\gamma}$ defined for $s \in (0,t)$ and 
converging to a point in $K$ as $s \to 0.$ Then by integrating (\ref{hga}) we obtain that the 
whole flow line lies at distance at most $\nu^{-1}t < \delta/2$ from $K,$ hence it extends to 
$s=t.$ The desired claim follows. $\hfill \fin$

\medskip

Consider now the circle $\gamma: [0,2\pi] \to \cc^*, \ \gamma_t(s) = t \ exp(is).$ A 
similar argument to that in the proof of Lemma~\ref{wdw} shows that for $0 < |t| < \nu 
\delta/2,$ parallel transport $h_{\gamma_t}(y)$ is well-defined and lies in $V$ for any $y 
\in L_t.$ Furthermore, we have

\begin {lemma}
The relative vanishing cycle $L_t$ is Lagrangian isotopic to $h_{\gamma_t}(L_t)$ inside $V 
\cap Z_t.$
\end {lemma}

\noindent \textbf {Proof. } By our estimates, the vanishing cycle $L_{\tau}$ is well-defined 
for any $\tau = \gamma_t(s),$ and it is moved by parallel transport along $\gamma_t|_{[s, 
2\pi]}$ without going outside $V.$ The required isotopy is given by $h_{\gamma_t|{[s,
2\pi]}}(L_{\tau}). \hfill \fin$

\subsection {Iterating the relative $\cpn$ construction.}
Now we deal with the situation relevant for Lemma~\ref{multithickthin}. Take $Z = X \times 
\bigl ( \slan \bigr )^k,$ with the map
$$ \pi: Z \to \cc^k, \ \pi(x, A_1, \dots, A_k) = \bigl (\chin(A_1), \dots, 
\chin(A_k) \bigr ), $$
and an effective K\"ahler form (defined just like in the previous section).

Starting with a compact Lagrangian submanifold $K$ of $X= X \times \{0\}^k \subset Z,$ one can 
use the first component of $\pi$ to construct a relative vanishing projective space $L_{t_1}$ 
in $ X \times \slan \times \{0\}^{k-1},$ for $t_1 \in \cc^*$ small. Using the 
second component next, and then iterating the process produces a Lagrangian submanifold 
$L_{t_1, \dots, t_k} \subset \pi^{-1}(t_1, \cdot, t_k)$ diffeomorphic to $X \times (\cpn)^k.$ 
Parallel transport estimates show that this is well-defined whenever $0 \leq |t_j| \leq 
\sigma$ for all $1\leq j \leq k,$ where $\sigma$ is a bound independent of $j.$

\begin {lemma}
\label {changeorder}
Changing the order in which the components of $\pi$ are used to construct $L_{t_1, \dots, 
t_k}$ produces the same outcome up to Lagrangian isotopy, at least as long as all the $|t_j|$ 
are sufficiently small.
\end {lemma}

\noindent \textbf {Proof. } The statement is trivial when the K\"ahler form is the product of 
a form on $X$ with the standard forms on each $\slan.$ The case of a general 
effective K\"ahler form follows from this by a Moser Lemma argument similar to the one in 
Lemma 30 in \cite{SS}. $\hfill \fin$

\section {An analogue of the fibered $A_2$ singularity} \label {sec:fiba2}

Our goal for this section is to understand the behavior of vanishing projective spaces 
under parallel transport in the situation of Section~\ref{sec:ttt}. The case $n=2$ is a 
fibered $A_2$ singularity and was treated in \cite[Section 4(C)]{SS}. The general case 
discussed here runs parallel to that, but the geometry is somewhat more complicated.

\subsection{The non-fibered case} \label{sec:non} We start with analyzing the situation 
in Section~\ref{sec:exa}. There we had the space $\xx_n \subset \sl(n+1, \cc)$ 
consisting of the matrices of the form (\ref{sjmm}) which admit an $(n-1)$-dimensional 
eigenspace. A matrix $A \in \sl(n+1, \cc)$ has a ``thick'' eignevalue $d$ of multiplicity at 
least $n-1,$ and two other ``thin'' eigenvalues which we write as $d+z_1$ and $d+z_2.$ (Of 
course, there can be coincidences giving higher multiplicities.) The zero trace condition 
implies that $z_1 + z_2 = -(n+1)d.$ We also have the map (\ref{mapp}) given by
$$ \mapp: \xx_n \to \sym^2(\cc) \cong \cc^2.$$

As coordinates on $\cc^2$ we take $d$ and the product $z = z_1z_2,$ so that $\mapp(A) = 
(d,z).$ By Proposition~\ref{luna} the map $\mapp$ is a smooth fibration except over the 
critical value sets corresponding to $z=0$ (a thick-thin coincidence) or $z= ((n+1)d/2)^2$ (a 
thin-thin coincidence). 

We view $d$ as an auxiliary parameter. (Our choice of $d$ is somewhat different from the one 
in \cite{SS}.) We denote the restriction of $\mapp$ to some $\xx_{n,d} = \mapp^{-1}(\{d \} 
\times \cc)$ by $\mapp_d : \xx_{n,d} \to \{d \} \times \cc = \cc,$ and we write $\xx_{n,d,z}$ 
for a fiber of $\mapp_d.$ For $d \neq 0,$ the map $\mapp_d$ has two critical values $0$ and 
$\zeta_d = ((n+1)d/2)^2.$ These coalesce for $d =0.$

Let us make this picture more explicit. A matrix $A \in \xx_{n,d}$ is of the form
\begin {equation}
\label {sjmm2}
dI +   \left( \begin {array}{ccccc} \alpha & 1 & 0 & \cdots & 0\\
a_{11} + \alpha^2 & \alpha & a_{12} & \cdots & a_{1n}\\
a_{21} & 0 & a_{22} & \cdots & a_{2n} \\
 & \cdots & \  & \cdots  &  \ \\
a_{n1} & 0 & a_{n2} & \cdots & a_{nn}
\end {array} \right ).
\end {equation}

(As compared to (\ref{sjmm}), we changed the meaning of the notation $a_{ij}.$) The fact that 
$A - dI$ has rank at most two translates into the condition that the associated {\it reduced 
matrix}
\begin {equation}
\label {sjmm3}
A_{\red}=   \left( \begin {array}{cccc} 
a_{11} & a_{12} & \cdots & a_{1n}\\
a_{21} &  a_{22} & \cdots & a_{2n} \\
 \cdots & \  & \cdots  &  \ \\
a_{n1}  & a_{n2} & \cdots & a_{nn}
\end {array} \right )
\end {equation}
has rank at most one. Note that 
\begin {equation}
\label {alfa}
2\alpha + a_{22} + \cdots + a_{nn} = -(n+1)d. 
\end {equation}

If we denote $s=a_{22} + \cdots + a_{nn},$ an evaluation of the sum of $2$-by-$2$ minors of $A$ 
gives $z=-a_{11}+2\alpha s,$ or
\begin {equation}
\label {betha}
z = -a_{11} -(n+1)d\cdot s - s^2.
\end {equation}

Thus $\xx_{n,d}$ can be identified with the space of matrices $A_{\red}$ of rank at most one,
and the map $\mapp_d: \xx_{n,d} \to \cc$ is given by 
\begin {equation}
\label {rede}
\mapp_d(A_{\red}) = -a_{11} - (n+1)d\sum_{k \geq 2} a_{kk} - \bigl ( \sum_{k \geq 2}  
a_{kk} \bigr )^2.
\end {equation}

Over the critical value $z=0$ we have the element in $\xx_{n,d}$ corresponding to $A_{\red} 
=0.$ 

\begin {lemma}
\label {biholo}
There exists a biholomorphism $f$ between a neighborhood of $A_{\red}=0$ in $\xx_{n,d}$ and a 
neighborhood of zero in the space $\slan$ from Section~\ref{sec:mone}, having 
the property that 
\begin {equation}
\label {propr}
\chin \circ  f = \mapp_d.
\end {equation}

This biholomorphism extends to one 
between neighborhoods of zero in $\gl(n, \cc),$ in the sense that there is a commutative 
diagram
$$ \begin {CD}
\xx_{n,d} @>{\text{local } \cong \ \text{near} \ 0}>> \slan \\
@VVV   @VVV \\
\gl(n, \cc) @>{\text{local } \cong \ \text{near} \ 0}>> \gl(n, \cc),
\end {CD} $$
where the vertical maps are the natural inclusions.
\end {lemma}

\noindent \textbf {Proof. } The spaces $\xx_{n,d}$ and $\slan$ are clearly 
biholomorphic: one can subtract a suitable multiple of the identity from a reduced matrix 
in $\xx_{n,d}$ and make it traceless. However, this obvious biholomorphism does not satisfy 
the required condition (\ref{propr}). What we need is to find a biholomorphic change of 
coordinates near zero 
in $\xx_{n,d}$ that takes the expression 
\begin {equation}
\label {eke}
 -\frac{1}{n} (a_{11} + \cdots + a_{nn})
\end {equation}
into the right-hand side of (\ref{rede}). 

It is helpful to work with a different set of coordinates on $\xx_{n,d}.$ By Lemma~\ref{git}, 
the variety $\slan$ and hence $\xx_{n,d}$ can be identified with a GIT 
quotient $\cc^{2n}/\cc^*.$ We denote by $v_i, w_j \ (i, j = 1, \dots, n)$ the coordinates on 
$\cc^{2n},$ so that the entries $a_{ij}$ of a reduced matrix $A_{\red}$ (for an element of  
$\xx_{n,d}$) are $a_{ij} = v_i w_j.$ Consider the holomorphic map from $\cc^{2n}$ to itself 
which keeps all $v_i$ fixed, takes $w_1$ into $nw_1,$ and
$$ w_j \to n(n+1) d\cdot w_j + nw_j(v_2w_2 + \dots + v_nw_n) \ \text{for } j=2, \dots, n.$$

The differential of this map at zero is invertible, hence this is a biholomorphism between 
neighborhoods of zero in $\cc^{2n}.$ Since it is $\cc^*$-invariant with respect to 
(\ref{cstar}), it induces a local biholomorphism between the GIT quotients. Written in terms 
of $a_{ij},$ this takes (\ref{eke}) into (\ref{rede}), as desired. 

The same formulae in terms of $a_{ij}$ also give the required extension to a local 
biholomorphism in $\gl(n, \cc). \hfill \fin$ 

\begin {remark} \label {remc} More canonically, suppose that instead of $\xx_n$ we have a 
space $\xx(T) \subset \sl(E)$ constructed from a linear quadruple $T$ as in 
Remark~\ref{act}. A reduced matrix represents a linear map from $F_1 \oplus T_1$ to $T_2 
\oplus T_1.$ Using the isomorphism $T_2 \to F_1$ from the definition of a linear quadruple, 
we can think of reduced matrices as automorphisms of $F_2 = F_1 \oplus T_1.$ Hence,   
Lemma~\ref{biholo} provides a local biholomorphism between neighborhoods of zero in 
$\xx(T)_d$ and $\sl(F_2)^{(1(n-1))},$ respectively. This extends to a local 
self-biholomorphsim of $\gl(E).$ \end {remark}

\subsection {Effective K\"ahler forms}
\label {sec:ek}

By analogy with Definition~\ref{eff} we make the following

\begin {definition}
\label {efftwo}
Denote by $\Reg \xx_n$ the open dense set of smooth points of the variety $\xx_n,$ and by $i$  
the standard inclusion of $\xx_n$ in $\sl(n+1, \cc).$ A K\"ahler form
$\Omega$ on $\Reg \xx_n$ is called effective if it is of the form $i^*\omega,$ where 
$\omega$ is a K\"ahler form on $\sl(n+1, \cc)$ that is real analytic in a neighborhood of 
zero.
\end {definition}

If we denote by $(b_{ij}), \ i,j=1, \dots, n+1$ the entries of a matrix $B \in \sl(n+1, 
\cc),$ we call $$ -dd^c \Bigl ( |b_{11} - b_{22}|^2 + \sum_{(i,j) \neq (1,1), (2,2)} 
|b_{ij}|^2 \Bigr ) $$ the {\it standard K\"ahler form} on $\sl(n+1, \cc).$
Its pullback to $\xx_{n,d}$ is called the standard K\"ahler form on that space, and denoted 
$\Omega_{st}.$ (Strictly speaking, because of the zero trace condition, there is no 
``standard'' choice of basis for $\sl(n+1, \cc).$ Our choice is the most natural here 
because we deal with matrices of the form \ref{sjmm2}.)

Fix a small (relatively compact) neighborhood $V$ of zero in $\xx_n.$ The analogue of 
Lemma~\ref{commens} still holds. If we restrict $\Omega_{st}$ and any other effective 
K\"ahler form on $\Reg \xx_{n,d}$ to the regular set $\Reg V$, then the corresponding 
metrics $g$ and $g_{st}$ are equivalent.

For small $d \in \cc, $ we can define two natural K\"ahler forms on $\Reg \xx_{n,d}.$ One is 
the restriction of $\Omega_{st},$ which we call the {\it first standard K\"ahler form} and 
denote by $\Omega_{st(1)}.$ The other, denoted $\Omega_{st(2)}$ and called the {\it second 
standard K\"ahler form}, is the pullback of the usual K\"ahler form on $\cc^{n^2}$ under the 
inclusion of $\xx_{n,d}$ (viewed as the set of reduced matrices) into $\gl(n, \cc) = 
\cc^{n^2}.$

The only difference between $\Omega_{st(1)}$ and $\Omega_{st(2)}$ appears because of the 
entry $a_{11} + \alpha^2$ in (\ref{sjmm2}). Using (\ref{alfa}) we see that this quadratic 
discrepancy is negligible for $V$ and $d$ sufficiently small. More precisely, the first and 
second standard metrics $g_{st(1)}$ and $g_{st(2)}$ induced by $\Omega_{st(1)}$ and 
$\Omega_{st(2)}$ on $\xx_{n,d}$ are equivalent on $\Reg (V \cap \xx_{n,d}).$ Moreover, the 
constants appearing in the definition of equivalence for metrics can be taken to be 
independent of $d,$ for $d$ sufficinetly small. Note that from here it follows that the 
same statements are also true for $g_{st(2)}$ and any metric induced by the restriction of 
an effective K\"ahler form on $\Reg \xx_n.$

\subsection {Estimating parallel transport} \label {sec:seces}

We work on the same open set $V \subset \xx_n$ as in the previous subsection. Endow 
$\Reg \xx_n$ with an effective K\"ahler form $\Omega$ in the sense of 
Definition~\ref{efftwo}. It is easy to see that the restrictions of $\Omega$ to $\Reg 
\xx_{n,d} \iso \slan - \{ 0 \}$ are effective in the sense of Definition~\ref{eff}. 

Using Lemmas \ref{biholo} and \ref{exist} we can find a natural Lagrangian vanishing 
projective space (with respect to the map $\mapp_d$)
$$ L_{d,z} \cong \cpn \subset \xx_{n,d,z},$$ for any $z$ with $ |z| \ll |d|.$ If $z$ and $d$ 
are small enough this lies inside $V.$

\begin {figure}
\begin {center}
\input {path.pstex_t}
\end {center}
\caption {Swapping the thin eigenvalues.}
\label {fig:path}
\end {figure}
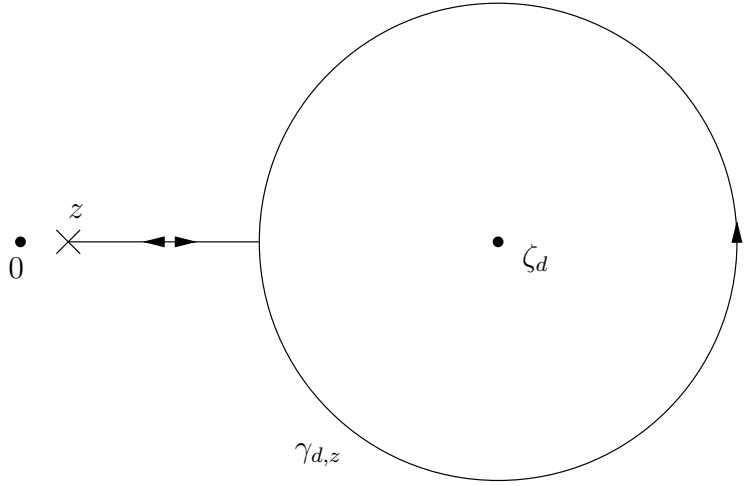

For $d, z > 0$ real ($z \ll d$), consider the path $\gamma_{d,z}$ in $\cc - \{0, \zeta_d 
\}$ going from $z$ to $\zeta_d/2$ on the real axis, then making a positive full circle 
around $\zeta_d,$ and going back to $z$ along the real axis. (See Figure~\ref{fig:path}.)
Parallel transport along this path basically corresponds to swapping the order of the two 
thin eigenvalues $d+z_1$ and $d+z_2.$

The image of $L_t$ under parallel transport along $\gamma_{d,z}$ is another Lagrangian 
projective space
\begin {equation}
\label {imago}
 h_{\gamma_{d,z}}(L_{d,z}) \subset \xx_{n,d,z}.
\end {equation}

In order to show that this is well-defined, we need to estimate the size 
of the parallel transport vector fields.

\begin {lemma}
\label {maines}
On each $\Reg (V \cap \xx_{n,d})$ we have
\begin {equation}
\label {deq}
 \| \nabla \mapp_d \|^2 \geq \nu \cdot |\mapp_d - \zeta_d|,
\end {equation}
where $\nu > 0$ is a constant that depends only on $n, V,$ and the K\"ahler form $\Omega.$
\end {lemma}

\noindent \textbf{Proof. } By the observation at the end of Section~\ref{sec:ek}, it suffices 
to prove the required bound using the second standard K\"ahler form on each $\Reg \xx_{n,d}.$ 

If we think of $\xx_{n,d}$ as the space of reduced matrices, we have an action of $U(n-1)$ on 
this space given by conjugation with the block matrix

$$ \left( \begin {tabular}{c|c}
$1$  &  $\begin {array}{ccc}
0 & \cdots & 0 \\
\end {array}$

 \\
\hline

$\begin {array}{c}
0 \\
\cdots \\
0 \\
\end {array}$ 
  &

$ \begin {array}{ccc} 
& & \\
& U & \\
& & \\
\end {array} $

\end {tabular} \right)
$$
for $U \in U(n-1).$ By using this action we can transform every matrix
$A_{\red}= (a_{ij})$ of rank at most one into one satisfying $a_{ij}=0$ for all $i \geq 3, 
j
\geq 2.$ Thus it suffices to prove the inequality (\ref{deq}) for matrices of this form, 
because both the K\"ahler metric $g_{st(2)}$ and the expression (\ref{rede}) are 
$U(n-1)$-invariant.

Since all the $2$-by-$2$ minors of $A_{\red}$ are zero, we can in fact assume that $a_{ij} 
= 0$ 
for $i \geq 3$ and any $j.$ (Indeed, if $a_{i1} =0$ for some $i \geq 3,$ then all columns 
except the first are zero, and by transposing we are back in the same case.)

Just as in the proof of Lemma~\ref{piano}, it is now helpful to restrict $\mapp_d$ to various 
subspaces where the size of its gradient is easier to estimate. The first restriction is to 
the space of reduced matrices with the property mentioned in the previous paragraph: 
\begin {equation}
\label {aij}
a_{ij} = 0 \ \text{ for } \ i \geq 3 \ \text{ and any  }  j. 
\end {equation}

Inside of this space we have the open dense subset $W$ consisting of matrices of the form
$$ \left( \begin {array}{cccc}
uw_1 & uw_2 & \cdots & uw_n\\
w_1 & w_2 & \cdots & w_n \\
0 & 0 & \cdots & 0 \\
& \cdots & \  & \cdots  \\
0 & 0 & \cdots & 0
\end {array} \right ), $$
with $u, w_j \in \cc.$ The only matrices that satisfy (\ref{aij}) and are not in $W$ have 
only zeros in the second row, but not so in the first. Note that if we establish 
(\ref{deq}) for 
matrices in the dense set $W,$ by taking limits the same inequality must be true for the 
remaining matrices too.

Therefore we just focus our attention on the set $W.$ In terms of the coordinates $u, w_j,$ 
the restriction of $q_d$ to $W \cong \cc^{n+1}$ is given by
$$ p_d = q_d|_W : W \to \cc, \ \ (u, w_j) \to  -uw_1 - (n+1)d\cdot w_2 - w_2^2.  $$

Since $\| \nabla p_d \| \leq  \| \nabla q_d \|,$ the inequality (\ref{deq}) would follow if we 
were able to prove
\begin {equation}
\label {what}
\| \nabla p_d \|^2 \geq \nu \cdot |p_d - \zeta_d|.
\end {equation}

The restriction of the second standard K\"ahler form to $W$ is
$$ -dd^c \Bigl ( (|u|^2 + 1) (|w_1|^2 + \dots + |w_n|^2) \Bigr ). $$

Let us consider the restriction of $p_d$ to the subspace $W_1 \cong \cc \subset W$ which 
corresponds to keeping all coordinates except $w_1$ fixed. We get
\begin {equation}
\label {in1}
\| \nabla p_d \|^2 \geq \| \nabla p_d|_{W_1} \|^2 = \frac{|u|^2}{|u|^2 + 1}.
\end {equation}

Doing the same thing with for the $w_2$-coordinate we obtain
\begin {equation}
\label {in2}
\| \nabla p_d \|^2 \geq  \frac{|(n+1)d + 2w_2|^2}{|u|^2 + 1}.
\end {equation}
 
Finally, using the $u$-coordinate:
\begin {equation}
\label {in3}
\| \nabla p_d \|^2 \geq \frac{|w_1|^2}{|w_1|^2 + \dots + |w_n|^2}.
\end {equation}

The matrices for which we seek to prove (\ref{what}) all lie inside the relatively compact 
set $V,$ which means that there is an a priori upper bound $R$ on the absolute values of all 
$w_j$'s. There is no such bound on $u.$ Nevertheless, if $|u| > 1,$ then (\ref{in1}) would 
automatically imply $\| \nabla p_d \| \geq 1/2,$ and then (\ref{what}) would follow from the a 
priori bound on its right-hand side. Thus we can assume $|u| \geq 1.$ Putting (\ref{in1}) 
and (\ref{in3}) together and using these bounds we obtain
$$ \| \nabla p_d \|^2 \geq \max \bigl( |u|^2/2, |w_1|^2/(nR^2) \bigr ) \geq (2n)^{-1/2}R^{-1} 
|uw_1|. $$

Combining this with (\ref{in2}) gives
\begin {eqnarray*}
 \| \nabla p_d \|^2 & \geq & \max \bigl((2n)^{-1/2}R^{-1} |uw_1|, 2|(n+1)d/2 + w_2|^2) \\
& \geq & \nu\cdot |uw_1 + ((n+1)d/2+w_2)^2| = \nu \cdot | p_d - \zeta_d |, 
\end {eqnarray*} 
as desired. $\hfill \fin$

\begin {lemma}
For $0 < z \ll d$ small, parallel transport along $\gamma_{d,z}$ is well-defined near the 
vanishing projective space $L_{d,z} \subset V,$ and the image (\ref{imago}) still lies in $V.$
\end {lemma}

\noindent \textbf{Proof. } Note that the length of $\gamma_{d,z}$ is $(1+\pi)\zeta_d,$ and all 
of its points are at distance at least $\zeta_d/2$ from $\zeta_d.$ Using Lemma~\ref{maines} 
and (\ref{hgam}), we deduce that any flow line of $H_{\gamma_{d,z}}$ contained in $W$ must 
satisfy $$ \int 
\|H_{\gamma_{d,z}} \| = \int \|\nabla q_d \|^{-1} \leq (1+\pi)\zeta_d \cdot \nu^{-1/2} 
(\zeta_d/2)^{-1/2} = C \cdot d, $$
where $C$ is a constant depending only on $n,$ $V,$ and $\Omega.$ The conclusion follows just 
as in the proof of Lemma~\ref{wdw}.
$\hfill \fin$

\subsection {Projection to the $\alpha$-coordinate} 
\label {sec:ac} Let us assume that $\Reg \xx_{n}$ is 
endowed with its standard K\"ahler form, so that on each $\xx_{n,d}$ (and its subspaces) we 
have the first standard K\"ahler form. In this situation, a more concrete picture of the 
vanishing cycles can be obtaines by projecting $\xx_{n,d,z}$ to the $\alpha$-coordinate, where 
$\alpha$ is given by the expression (\ref{alfa}) in terms of the $a_{ij}$'s.

Denote by $\xx_{n,d,z,\alpha}$ the subspace of $\xx_{n,d,z}$ corresponding to fixed $\alpha.$ 
Using (\ref{alfa}), (\ref{betha}), and the representation of $\xx_{n,d}$ in terms of the 
variables $v_i, w_j$ as in the proof of Lemma~\ref{biholo}, we get that $\xx_{n,d,z,\alpha}$ 
is isomorphic to the variety
$$ \mathscr{V}_{c_1, c_2} = \{ v_i, w_j \in \cc \ (i,j=1, \dots,n)\ | \ v_1w_1 = c_1, \ 
v_2w_2 + \dots + v_n w_n = c_2 \}/\cc^*.$$
Here the $\cc^*$-action is given by (\ref{cstar}), and
\begin {eqnarray*}
c_1 &=& -z-2(n+1)d \cdot \alpha - 4\alpha^2; \\
c_2 & = & -(n+1)d - 2\alpha.
\end {eqnarray*}

\begin {lemma}
The variety $ \mathscr{V}_{c_1, c_2} $ is nonsingular if and only if $c_1 \neq 0$ and $c_2 
\neq 0.$ 
\end {lemma}

\noindent \textbf{Proof. } If $c_1=0,$ the variety has two irreducible components 
corresponding to $v_1 = 0$ and $w_1=0,$ respectively, and the two components have a nontrivial 
intersection.

If $c_1 \neq 0,$ we can use the $\cc^*$-action to fix $v_1 =c_1, \ w_1 =1.$ We get that $ 
\mathscr{V}_{c_1, c_2}$ is the quadric given by the equation $v_2w_2 + \dots + v_n w_n = c_2.$
This is singular if and only if $c_2 =0. \hfill \fin$

\medskip

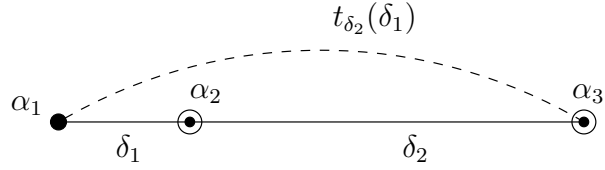
\begin {figure}
\begin {center}
\input {lag.pstex_t}
\end {center}
\caption {Projections of the Lagrangians.}
\label {fig:lag}
\end {figure}

Therefore, for fixed $n,d,z,$ there are typically three values of $\alpha$ for which 
$\xx_{n,d, z,\alpha}$ is singular: $\alpha_1 = -(n+1)d/2,$ and the two roots $\alpha_2$ and 
$\alpha_3$ of the equation $4\alpha^2 + 2(n+1)d\alpha + z=0.$ Note that for $0 < z \ll d,$ 
all these three values are real and negative. The leftmost one is $\alpha_1,$ then one of 
$\alpha_2$ and $\alpha_3$ (say $\alpha_2$) is slightly to the right of $\alpha_1,$ and 
the remaining value $\alpha_3$ is close to zero. (See Figure ~\ref{fig:lag}.)

Let us consider the $U(n-1)$-action on $\xx_{n,d}$ that appeared in the beginning of 
the proof of Lemma~\ref{maines}. The fibers $\xx_{n,d,z, \alpha}$ are preserved by this 
action, which in terms of the coordinates $\vec{v} =(v_2, \dots, v_n), \vec{w} =(w_2, 
\dots, w_n)$ is simply 
\begin {equation}
\label {un}
U \in U(n-1) \ : \ (\vec{v}, \vec{w}) \to (U\vec{v}, U^{-1}\vec{w}).
\end {equation}

This action is Hamiltonian with respect to $\Omega_{st(1)}.$ If $\mu: \xx_{n,d,z, \alpha} 
= \mathscr{V}_{c_1, c_2} \to \mathfrak{u}(n-1)$ denotes the moment map, the set 
$\mu^{-1}(0)$ is determined by the system of equations
$$ v_i \bar{v}_j = \bar{w}_i w_j \ \ (i,j =2, \dots, n). $$

If these equations are satisfied, then first of all $|v_i| = |w_i|$ for all $i,$ and then 
we know that $\bar{w}_i/v_i$ is independent of $i.$ This means that $w_i = \lambda \bar 
v_i$ for some $\lambda \in S^1.$ The equation $v_2w_2 + \dots + v_n w_n=c_2$ becomes 
$|v_2|^2 + \dots + |v_n|^2 = \lambda^{-1}c_2,$ and (unless $c_2=0$) the value of $\lambda$ 
is determined by $c_2$ in such a way as to make $\lambda^{-1}c_2$ a positive real number.

In light of this, the intersection
\begin {equation} \label {cndz}
\mathscr{C}_{n,d,z,\alpha} = \xx_{n,d,z,\alpha} \cap \mu^{-1}(0) 
\end {equation}
is easily seen to be a $(2n-3)$-dimensional sphere if $c_1, c_2 \neq 0,$ a point if $c_2 = 
0,$ and a copy of $\mathbb{CP}^{n-2}$ if $c_1=0, c_2 \neq 0.$ 

\begin {remark}
\label {klm}
If we are in the situation from Remark~\ref{remc}, with $\xx(T)$ instead of 
$\xx_n,$ then for $c_1, c_2 \neq 0,$ the space (\ref{cndz}) can be naturally identified 
with a zero-centered sphere in $\text{Hom}(F_1, T_1)=\text{Hom}(F_1, F_2/F_1).$ Indeed, 
if we use the 
$\cc^*$-action to set $w_1=1,$ then the $v_j$'s ($j \geq 2$) are elements on the 
first column of the reduced matrix $A_{\red}.$ This matrix represents a linear operator 
in $\gl(F_2),$ written in a basis in which the first component is in $F_1.$ 
\end {remark}

Given a path $\delta: [0,1] 
\to \cc$ such that $\delta(0)=\alpha_1, \ \delta(1) \in \{\alpha_2, \alpha_3\}$ and 
$\delta(t) \not \in \{\alpha_1, \alpha_2, \alpha_3\}$ for $t\in (0,1)$, we can construct a 
Lagrangian $\cpn \subset \xx_{n,d,z}:$ 
\begin {equation}
\label {bigla}
\Lambda_{\delta} = \bigcup_{t\in [0,1]} \mathscr{C}_{n,d,z, \delta(t)} 
\end {equation}
as the union of a point sitting over $\alpha_1,$ a family of $S^{2n-3}$'s sitting over 
$\delta(t), t \in (0,1),$ and a $\mathbb{CP}^{n-2}$ sitting over the other endpoint 
($\alpha_2$ or $\alpha_3$).

Suppose now that $0 < z \ll d.$

\begin {lemma}
\label {earth}
The vanishing projective space $L_{d,z} \subset \xx_{n,d,z}$ with respect to 
$\Omega_{st(1)}$ and its 
monodromy image (\ref{imago}) are Lagrangian isotopic to the projective spaces 
(\ref{bigla}) associated to the paths $\delta_1$ and $t_{\delta_2}(\delta_1)$ shown in 
Figure~\ref{fig:lag}. (Here $t_{\delta_2}$ stands for a half-twist around $\delta_2.$)
\end {lemma}

\noindent \textbf{Proof. } This is completely similar to the proof of Lemma 32 in 
\cite{SS}. A short computation shows that the parallel transport vector 
fields for $\mapp_d: \xx_{n,d} \to \cc$ are invariant with respect to the $U(n-1)$-action 
(\ref{un}), and $d\mu$ vanishes on them. The critical point corresponding to $A_{\red}=0$ 
is a fixed point of the $U(n-1)$-action, and lies in $\mu^{-1}(0).$ It follows that the 
vanishing cycles are $U(n-1)$-invariant and lie in $\mu^{-1}(0)$ as well. Any such 
Lagrangian is of the form $\Lambda_{\delta}$ for some path $\delta$ going from $\alpha_1$ 
to either $\alpha_2$ or $\alpha_3.$ This information, and the fact that $L_{d,z}$ must lie 
close to the critical point $A_{red}=0$ for $z$ small, determines the isotopy class of 
$L_{d,z}$ uniquely. A similar argument works for the mondromy image, but in that case the 
thin eigenvalues (and hence $\alpha_2$ and $\alpha_3$) are being swapped. $\hfill \fin$

\begin {remark}
\label {moon}
The Lagrangians $L_{d,z}$ and $h_{\gamma_{d,z}}(L_{d,z})$ (constructed using the first 
standard K\"ahler form) can be isotoped to intersect exactly in one point.
\end {remark}

\subsection {The fibered case.} \label {sec:fibcase}
We now consider the situation appearing in 
Lemma~\ref{ttt} and Remark~\ref{remcom2}. Our discussion parallels that at the end 
of Section 4(C) in \cite{SS}.

Let $X$ be a complex manifold and $\tt$ a holomorphic quadruple bundle over $X$ as in 
Section~\ref{sec:ttt}. Over $X$ there are associated holomorphic vector bundles $\ff_1 
\subset \ff_2 
\subset \ee$ of ranks $1, n, n+1,$ respectively. We also have a fiber bundle $\xx(\tt) 
\subset \sl(\ee)$ over $X$ with fibers isomorphic 
to $\xx_n,$ constructed using Remark~\ref{act}. We make the assumption that the 
bundle $\ff_2 \to X$ is trivial. (This is satisfied in Lemma~\ref{ttt}, cf. the 
observation preceding Remark~\ref{remcom2}.) Consider also the map $q: 
\xx(\tt) \to \sym^2(\cc) \cong \cc^2$ which is equal to (\ref{mapp}) in every $\xx_n$ 
fiber. 

We identify $X$ with the zero section in $\sl(\ee),$ and endow $\sl(\ee)$ with a K\"ahler 
metric $\omega.$ Let $K, K'$ be two closed Lagrangian submanifolds of $X,$ and assume that 
$\omega$ is real analytic in a neighborhood of $K, K' \subset \sl(\ee).$ Using parallel 
transport in the union of the critical value sets $Crit(q_d) \cong X \subset 
\xx(\tt)_{d,0}$ 
corresponding to $A_{\red} = 0$, we obtain Lagrangian submanifolds $K_d, K'_d$ in 
$Crit(q_d)$ for any small $d.$

Applying the results of Section~\ref{sec:relvan} and taking into account the triviality 
of $\ff_2$ and Remark~\ref{remc}, we get associated relative vanishing 
projective spaces for $\pi_d$ in the form of Lagrangian submanifolds
\begin {equation}
\label {eins}
 L_{d, z}, \ L'_{d, z} \subset \xx(\tt)_{d,z}
\end {equation}
for $0 < z \ll d.$ These are diffeomorphic to $K \times \cpn$ and $K' \times \cpn,$ 
respectively.

Fix a relatively compact open subset $W \subset X$ containing $K, K',$ and an open subset 
$V \subset \xx(\tt)$ which is the unit ball bundle over $W$ with respect to the given 
metric. 
For $0 < z \ll d,$ the Lagrangians $L_{d,z}, L'_{d,z}$ will lie in $V.$ The estimates for 
parallel transport from Section~\ref{sec:seces} apply in this situation as well, and show 
that there is another well-defined Lagrangian submanifold
\begin {equation}
\label {zwei}
L''_{d,z} = h_{\gamma_{d,z}}(L'_{d,z}) \subset \xx(\tt)_{d,z} \cap V.
\end {equation}

Let us now look at a particular construction for the K\"ahler form $\omega$. Given an 
arbitrary K\"ahler form on $X$ and a Hermitian metric on $\ee \to X$ (both with good local 
real analyticity properties), we can combine them to obtain an associated form $\omega$ on 
$\sl(\ee)$ which is K\"ahler at least in a neighborhood of the zero set of the moment map 
$\mu$ for (\ref{un}). (For details of this construction, we refer to Remark 33 in 
\cite{SS}.) Given $\alpha \in \cc$ and a path $\delta$ in $\cc$ as in 
Section~\ref{sec:ac}, we get a Lagrangian submanifolds $\Lambda_{d,z,K, \delta}$ by 
applying the construction (\ref{bigla}) fiberwise over a Lagrangian $K \subset X.$ The 
arguments in the proof of Lemma~\ref{earth} carry over to give the following 
generalization:

\begin {lemma}
\label {earth2}
Equip $\Reg \xx(\tt)$ and the smooth fibers $\xx(\tt)_{d,z}$ with the restriction of a 
K\"ahler form $\omega$ 
constructed as in the previous paragraph. Given closed Lagrangian submanifolds $K, K' 
\subset X,$ consider the relative vanishing projective spaces (\ref{eins}) and the 
mondromy image (\ref{zwei}). Then, up to Lagrangian isotopy,
$$ L_{d,z} = \Lambda_{d, z, K, \delta_1}, \ \ \ L''_{d,z} = \Lambda_{d,z, K', 
t_{\delta_2}\delta_1}.$$
Here the paths $\delta_1$ and $t_{\delta_2}\delta_1$ are as in Figure~\ref{fig:lag}.
\end {lemma}

\section {Floer cohomology}
\label{sec:floer}

Lagrangian Floer cohomology (\cite{F}, \cite{FOOO}) can be defined in various settings, and it 
comes with various amounts of structure depending on how restrictive our assumptions are. We 
work here in the setting of \cite[Section 4(D)]{SS}, except that we need to relax the spinness 
condition on Lagrangians. 

\subsection {Definition} \label {sec:flo} Let $M$ be a Stein manifold endowed with an 
exact K\"ahler form $\omega.$ We assume that $c_1(M) = 0$ and $H^1(M) =0.$ Let $L,L'$ be 
two closed, connected, oriented Lagrangian submanifolds of $M$ satisfying $H_1(L) = 
H_1(L') =0.$ We also assume that the pair of Lagrangians $(L, L')$ is relatively spin in 
the sense of \cite{FOOO}, meaning that there exists a class $st \in H^2(M; \zz/2\zz)$ that 
restricts to $w_2(L)$ on $L$ and to $w_2(L')$ on $L'.$ The choice of the class $st$ is 
called a {\it relative spin structure}.

By doing a small Lagrangian perturbation we can arrange so that the intersection $L \cap L'$ is 
transverse. The Floer cochain complex is then defined to be the abelian group
$$ CF(L, L') = \bigoplus_{x \in L \cap L'} O_x,$$
where $O_x$ is the orientation group of $x.$ This is the abelian group 
(noncanonically isomorphic to $\zz$) which is generated by the two possible orientations of 
$x,$ with relation that their sum is zero. Our assumptions allow one to define a relative 
Maslov index $\Delta gr(x,y) \in \zz$ for every $x,y \in L \cap L'.$ The relative index 
staisfies $\Delta gr(x,y) + \Delta gr(y,z) = \Delta gr(x,z)$ and induces a relative 
$\zz$-grading on the cochain complex.

The Lagrangian Floer cohomology groups $HF^*(L,L')$ are the cohomology groups of $CF^*(L, L')$ 
with respect to the differential $d_J$ defined on generators by
$$ d_J(x) = \sum_y n_{xy} y.$$
Here $n_{xy} \in \zz$ is the signed count of isolated solutions (modulo translation in $s$) to 
Floer's equation
\begin{equation} \label{floer}
\left\{
\begin{aligned}
 & u: \rr \times [0,1] \rightarrow M, \\
 & u(s,0) \in L, \quad u(s,1) \in L', \\
 & \partial_s u + J_t(u) \partial_t u = 0, \\
 & \lim_{s \rightarrow +\infty} u(s,\cdot) = x, \  \lim_{s \rightarrow -\infty} u(s,\cdot) = y,
\end{aligned}
\right.
\end{equation}
where $J = (J_t)_{0 \leq t \leq 1}$ is a generic smooth family of
$\omega$-compatible almost complex structures, which all agree
with the given complex structure outside a compact subset. The solutions to (\ref{floer}) are 
called pseudo-holomorphic disks.

The fact that $M$ is Stein implies that the solutions of (\ref{floer}) remain inside a 
fixed compact subset of $M.$ The exactness of $\omega$ and the fact that $H^1(L) = H^1(L') 
= 0$ ensure that no bubbling occurs and thus the moduli spaces of solutions to 
(\ref{floer}) have well-behaved compactifications. Finally, the condition that the 
Lagrangians are relatively spin is necessary in order to define orientations on the moduli 
spaces. The orientations depend on some additional data, cf. \cite[p.192]{FOOO}: the class 
$st \in H^2(M; \zz/2\zz)$ which we called a relative spin structure, as well as 
(in principle) choices of spin structures on certain bundles over the two-skeleta of $L$ 
and $L'.$ However, in our case the latter piece of information is vacuous: since our 
assumptions imply that $H^1(L; \zz/2\zz) = H^1(L';\zz/2\zz)=0,$ the spin structures are 
unique. (This follows from \cite[p. 81, Corollary 1.5]{LM}, for example.) 

An important property of the Floer cohomology groups $HF^*(L, L')$ in our setting is that 
they are invariant under Lagrangian isotopies of either $L$ and $L'.$ Furthermore, they 
are invariant under any smooth deformation of the objects involved in their definition 
(for example, the K\"ahler metric), as long as all the assumptions which we made are still 
satisfied.

\subsection {Absolute gradings.} \label{sec:abs} As defined in the previous subsection, 
the Floer cohomology $HF^*(L, L')$ groups are only relatively $\zz$-graded. However, in 
the presence of some additional data, one can improve this to an absolute grading. This 
improvement is due to Seidel \cite{Se}, who was inspired by the ideas of Kontsevich 
\cite{Ko}.

Since $c_1(Y) = 0,$ we can pick a complex volume form $\Theta$ on $Y,$ i.e. a nowhere vanishing 
section of the canonical bundle. This determines a square phase map
\begin {equation}
\label {sqp1}
 \theta : \lfrak \to
\cc^*/\rr_+, \ \theta(V) = \Theta(e_1 \wedge \dots \wedge e_n)^2
\end {equation}
for any orthonormal basis $e_1, \dots, e_n$ of $T_xL.$ 
We can identify $\cc^*/\rr_+$ with $S^1$ by the contraction $z \to z/|z|.$ The condition 
$H^1(\lll) = 0$ allows us to lift $\theta_{\lll}$ to a real-valued map. A {\it grading} on 
$\lll$ 
is a choice $\tilde \theta_{\lll} :\lll \to \rr$ of such a lift. If we choose a grading for 
$\lll'$ 
as well, then every
point $x$ in the intersection $\lll \cap \lll'$ (which was assumed to be
transverse) has a well-defined absolute Maslov index $gr(x) \in \zz$ \cite{Se}, such that 
$\Delta gr(x,y) = gr(x) - gr(y).$
In turn, this gives an absolute $\zz$ grading on the cochain 
complex and on cohomology. The result does not depend on the choice of $\theta,$
because the condition $H^1(Y) =0$ ensures that different choices are
homotopic in the class of smooth trivializations of the canonical bundle.

We should note that if we denote by $L \to L[1]$ the process which subtracts the 
constant $1$ from the grading, then we have
\begin {equation}
\label {grad}
 HF^*(L, L'[1]) = HF^*(L[-1], L') = HF^{*+1} (L, L').
\end {equation}

\subsection {A K\"unneth formula} Let us consider Floer cohomology in the geometric 
situation from Section~\ref{sec:relvan}. We work in the following context, which is 
slightly more general than the one in Section~\ref{sec:relvan}:
 
\begin{itemize} 
\item 
$Z$ is a complex affine variety equipped with a holomorphic function $\pi: Z \rightarrow 
\cc.$ There is a neighborhood $D$ of the origin in $\cc$ such that the open set 
$\pi^{-1}(D - \{ 0 \}) \subset Z$ is smooth, and the restriction of $\pi$ to this open set 
is a smooth submersion.

\item
We have a smooth complex subvariety $X \subset Z$ and an isomorphism between a 
neighbourhood of that subvariety and a neighbourhood of $X \times \{0\} \subset X \times 
\slan$, such that the following diagram commutes: 
\begin{equation} \label {eql}
\begin{CD}
 Z @>{\text{local $\iso$ defined near $X$}}>> X \times \slan \\
 @V{\pi}VV @V{\chin}VV \\
 \cc @>{\quad\quad\qquad\qquad\quad}>> \cc
\end{CD}
\end{equation}

\item The regular set $\Reg Z$ of $Z$ carries an exact K{\"a}hler form $\Omega$ that is 
effective near $X$ in the sense of Section~\ref{sec:relvan}. In particular, since $\Omega$ 
is locally obtained by pullback from $X \times \sl(n, \cc),$ this automatically equips $X$ 
with a K\"ahler form as well.

\item
We assume that the first Chern classes $c_1(Z_t)$ for 
$t\neq 0,$ as well as $c_1(X),$ are zero. Also, $H^1(Z_t) = 0$ for small $t \neq 0$, and
$H^1(X) = 0$.
\end{itemize}

We endow all the smooth fibres $Z_t$ with the restrictions of $\Omega$. Note that 
(\ref{eql}) implies that every $Z_t$ (for $t \neq 0$ small) has an open set $U_t \subset 
Z_t$ diffeomorphic to $X \times T^*\cpn.$ Let $K,K'$ be closed Lagrangian submanifolds of 
$X$ which have the properties necessary to define the Floer cohomology $HF(K,K').$ In 
particular, the pair $(K, K')$ comes with a relative spin structure $st \in H^2(X; 
\zz/2\zz).$ For sufficiently small $t \neq 0$ we have associated relative vanishing 
projective spaces $L_t,L_t' \subset U_t \subset Z_t,$ as in Section~\ref{sec:relvan}. 
These are products of $K,K'$ with $\cpn$, so their Floer cohomology $HF(L_t,L_t')$ is 
again well-defined, provided we make sure that the pair $(L_t, L'_t)$ is relatively spin. 
This last condition is not automatic, so we introduce an additional piece of data:

\begin {itemize}
\item 
For small $t \neq 0,$ we have a class $\widetilde{st} \in H^2(Z_t; \zz/2\zz)$ such that 
its restriction to each $U_t \iso X \times T^*\cpn \subset Z_t$ is isomorphic to the product 
of $st$ and $\xi_n.$ Here $\xi_n$ is zero for $n$ even, and equals the generator of 
$H^2(T^*\cpn; \zz/2\zz)$ for $n$ odd. In other words, $\xi_n$ is the pull-back of the 
second Stiefel-Whitney class of $\cpn$ under the projection $T^*\cpn \to \cpn.$  
\end {itemize}

We use $\widetilde{st}$ to define $HF(L_t, L_t').$ Note that under these hypotheses, the 
Floer cohomology $HF(L_t, L'_t)$ is independent of $t$ by the invariance principle 
mentioned at the very end of Section~\ref{sec:flo}.

\begin{lemma} \label{kunneth}
$HF(L_t,L_t') \iso HF(K,K') \otimes H^*(\cpn)$, where $H^*(\cpn)$ is given its standard 
grading.
\end{lemma}

\noindent {\bf Proof. } This is completely analogous to the proof of Lemma 36 in 
\cite{SS}. One can find a holomorphically weakly convex 
neighborhood $V$ of $X$ in $Z$ on which the local isomorphism (\ref{eql}) is well-defined. 
For $t$ small, the Lagrangians $L_t$ and $L'_t$ are in $V.$ Pseudoholomorphic disks in 
$Z_t$ with boundaries in $V \cap Z_t$ cannot go outside $V,$ so the Floer cohomology can 
be computed in $V \cap Z_t.$ On that set we can deform the K\"ahler form to be the product 
of the form on $X$ with the restriction of the standard form $\Omega_{st}$ from the 
example following (\ref{vcps}). Floer cohomology is not changed by this deformation, and 
in the new K\"ahler form $L_t$ and $L_t'$ become products $K \times \cpn$ and $K' \times 
\cpn,$ respectively, where the $\cpn$ factor is the same. The result now follows from a 
straightforward K\"unneth product formula in Floer cohomology, similar to the one in Morse 
theory. The Floer cohomology orientations were chosen to exactly match. $\hfill \fin$

\subsection {Floer cohomology in the setting of Section~\ref{sec:fibcase}} We now analyze
the following situation:

\begin{itemize}
\item
$Z$ is a complex affine variety equipped with a holomorphic map $\pi: Z \to \cc^2.$
The variety $Z$ is smooth over the set of pairs $(d,z) \in \cc^2$ with $0 < z \ll d$
sufficiently small, and the map $\pi$ is a submersion there.

\item
We have a smooth complex subvariety $X \subset Z$ and a holomorphic quadruple bundle 
$\tt \to \xx$ in the sense of Section~\ref{sec:ttt}. There are associated holomorphic 
bundles $\ff_1 \subset \ff_2 \subset \ee$ over $X$ of ranks $1, n, n+1,$ respectively, and 
we assume that $\ff_2$ is trivial. From here we obtain a fiber bundle $\xx(\tt) 
\to X$ as a 
subset of $\sl(\ee),$ and a map $q:\xx(\tt) \to \cc^2.$ We assume that there is an 
isomorphism 
between a neighborhood of $X$ in $Z$ and a neighborhood of the zero section in $\xx(\tt) 
\subset \sl(\ee),$ such that the following diagram commutes:
\begin{equation} \label{eql2}
\begin{CD}
 Z
 @>{\text{local $\iso$ defined near $X$}}>>
 {\xx(\tt)} \\
 @V{\pi}VV @V{q}VV \\
 {\cc^2} @>{\quad\qquad\qquad\qquad\quad}>> \cc^2
\end{CD}
\end{equation}

\item 
The regular set $\Reg Z$ carries an exact K{\"a}hler form $\Omega$ which in a neighborhood 
of $X$ is obtained by pull-back from a real analytic K\"ahler form on an open set in 
$\sl(\ee).$ In particular, $X$ also has an induced K\"ahler form. 

\item
We require that $c_1(Z_{d,z}) = 0$ and $H^1(Z_{d,z})=0$ for $0 < z \ll d.$ Also, 
$c_1(X) =0$ and $H^1(X) =0.$
\end{itemize}

Take closed Lagrangian submanifolds $K,K' \subset X$, satisfying the conditions in 
Section~\ref{sec:flo}, so that $HF(K,K')$ is well-defined. For sufficiently small $0 < z 
\ll d,$ the construction in Section~\ref{sec:fibcase} produces some new 
Lagrangians $L_{d,z},L_{d,z}', L''_{d,z} \subset Z_{d,z}.$ These can be assumed to 
lie inside an open subset $U_{d,z} \subset Z_{d,z}$ diffeomorphic to 
$X \times T^*\cpn,$ and where the isomorphism (\ref{eql2}) is well-defined. Given the 
relative spin structure $st$ for the pair $(K, K'),$ we make the following 
additional assumption.

\begin {itemize} \item For small $0 < z \ll d,$ we have a class $\widetilde{st} \in 
H^2(Z_{d,z}; \zz/2\zz)$ such that its restriction to each $U_{d,z} \iso X \times T^*\cpn 
\subset Z_{d,z}$ is isomorphic to the product of $st$ and $\xi_n.$ Here $\xi_n$ is zero 
for $n$ even, and equals the generator of $H^2(T^*\cpn; \zz/2\zz)$ for $n$ odd. \end 
{itemize}

Under these hypotheses, we can use $\widetilde{st}$ to define the Floer cohomology  
$HF(L_{d,z}, L''_{d,z}),$ for $0 < z \ll d$ small.

\begin {lemma}
\label {thom}
$HF(L_{d,z}, L''_{d,z}) \iso HF(K, K').$
\end {lemma}

\noindent {\bf Proof. } This goes just like the proof of Lemma 38 in \cite{SS}, so we only 
sketch the argument. We can assume that $K$ and $K'$ intersect transversely. Floer 
cohomology can be calculated inside a suitable holomorphically weakly convex subset of 
$Z_{d,z}$ where (\ref{eql2}) is defined. On that set we can deform the K\"ahler form into 
one that is induced by the form on $X$ and a Hermitian form on $\ee,$ as in the discussion 
at the end of Section~\ref{sec:fibcase}. Lemma~\ref{earth2} and Remark~\ref{moon} show 
that $L_{d,z}$ and $L''_{d,z}$ can be isotoped so that their intersection points exactly 
correspond to intersection points of $K$ and $K'.$ With a little more care we can also 
arrange that pseudoholomorphic disks in the two settings are in one-to-one correspondence 
as well. This implies that the Floer cohomology groups are isomorphic. $\hfill \fin$

\section {Definition of the invariants}
\label {sec:maind}

This section contains the construction of the Floer cohomology groups $\krss.$ We apply the 
symplectic geometric techniques which we developed so far to the map
\begin {equation}
\label {our}
\chi^{\pi}|_{\ss_{m,n}}: \ss_{m,n} \cap \gg^{\pi} \to \hh^{\pi}/W^{\pi} 
\end {equation}
from Section~\ref{sec:props}. Recall that $\gg = \sl(mn, \cc), \ \pi$ is the partition 
$(1^m(n-1)^m),$ and the fibers of $\chi^{\pi}|_{\ss_{m,n}}$ are denoted $\ymnt.$ The 
restriction of (\ref{our}) to the set $\ss_{m,n} \cap \gg^{\pi, \reg}$ lying over 
the bipartite configuration space $BConf^0_m$ is a symplectic fibration.

We extend the notation $\ymnt$ to $\tau = (\lambdav, \muv)$ such that the trace $T=\sum 
\lambda_j + (n-1) \sum \mu_j$ is not necessarily zero. (When the $\lambda$'s and $\mu$'s are 
distinct, the values of $\tau$ form a configuration space $BConf_m = Conf^{\pi}$ as defined in 
Section~\ref{sec:rpr}.) In this more general situation, by $\ymnt$ we mean the fiber of 
(\ref{our}) over the normalized configuration made of $\lambda_j - T/(nm)$ and $\mu_j - 
T/(nm).$

\subsection {Parallel transport in $\ss_{m,n} \cap \gg^{\pi, \reg}$} 
\label {sec:kh}
We seek a K\"ahler metric $\Omega = -dd^c \psi$ on the total space of our symplectic fibration
\begin {equation}
\label {ourf}
 \ss_{m,n} \cap \gg^{\pi, \reg} \to BConf^0_m
\end {equation}
satisfying the conditions (\ref{cond1})-(\ref{cond4}), so that rescaled parallel 
transport is well-defined.

The function $\psi$ will be the restriction of some $\hat \psi: \ss_{m,n} \to \rr$ whose 
construction is completely analogous to that in \cite[Section 5(A)]{SS}. Pick some $\alpha > 
m.$ For each $k=2,4, \dots, 2m,$ apply the function $\xi_k(z) = |z|^{2\alpha/k}$ to the 
coordinates of $\ss_{m,n}$ on which the action (\ref{action}) is by weight $k$ or, in other 
words, to the entries of $Y_{1k}$. Sum up these terms to get a function $\xi: \ss_{m,n} \to 
\cc.$ This is not smooth at the origin, but we can perturb it using compactly supported 
functions $\eta_k$ on $\cc$ such that $\psi_k = \eta_k + \xi_k$ is $C^{\infty}.$ We 
choose the $\psi_k$ to be real analytic in a neighborhood of the origin, and strictly 
plurisubharmonic everywhere. Adding up all 
the $\psi_k$'s we obtain the function $\hat \psi$ on $\ss_{m,n},$ whose restriction 
to $ \gg^{\pi, \reg}$ is $\psi.$ We set $\Omega = -dd^c \psi.$ 
\begin {lemma}
\label {rescpt}
The function $\psi$ satisfies (\ref{cond1})-(\ref{cond4}).
\end {lemma}

The proof of Lemma~\ref{rescpt} is identical to those of Lemmas 41 and 42 in \cite{SS}, so we 
omit the details. The main ingredient in the proof of (\ref{cond2}) is that by construction 
$\psi$ is asymptotically homogeneous for the action of $\rr_+ \subset \cc^*$ by (\ref{action}). 
We can then control its critical points using this $\rr_+$-action, after pulling back to the 
resolution $\tilde \ss_{m,n} \cap \tilde \gg^{\pi}$ from (\ref{big}).

According to the discussion in Section~\ref{sec:parallel}, the function $\psi$ defines rescaled 
parallel transport maps 
$$ h^{\resc}_{\beta}: \yy_{m,n,\beta(0)} \to \yy_{m,n, \beta(1)}$$
for every path $\beta:[0,1] \to BConf^0_m.$ (Strictly speaking, $h_{\beta}$ is only 
well-defined on compact subsets, but these can get arbitrarily large, cf. the note at the end 
of Section~\ref{sec:parallel}). We also extend the notation $h_{\beta}$ to paths in $Conf^0_m,$ 
with the understanding that each point of the path is translated by $T/(mn),$ where $T$ is 
its trace.

Rescaled parallel transport also exists for the symplectic fibration (\ref{cmnt}). The total 
space $\ccc_{m,n}$ lives inside $\ss_{m,n},$ so we can take the restriction of the function 
$\hat \psi$ and the corresponding K\"ahler form. Given a path $\bar \beta: [0,1] \to 
Conf^0_{\sigma}$ with $\sigma=(1^{m-1}(n-1)^{m-1}n)$ as in Section~\ref{sec:thickthin}, we get 
a rescaled parallel transport map $h^{\resc}_{\bar \beta}.$ Denote by $BConf^{0*}_{m-1} \subset 
Conf^0_{\sigma}$ the subset corresponding to the multiplicity $n$ eigenvalue being zero. We can 
identify $BConf^{0*}_{m-1}$ to an open subset of $BConf^0_{m-1}$ in the obvious way. Note that 
according to Lemma~\ref{identh}, a fiber $\ccc_{m,n,\tau}$ with $\tau=\{(0, \lambda_2, \dots, 
\lambda_m), (0, \mu_2, \dots, \mu_m)\}$ can be identified with $\yy_{m-1,n, \bar \tau},$ with 
$\bar \tau =\{(\lambda_2, \dots, \lambda_m), (\mu_2, \dots, \mu_m)\}.$ This identification is 
compatible with our choice of symplectic forms, provided we take the same $\alpha > 0$ and 
functions $\psi_k$ for both $m$ and $m-1.$ Thus parallel transport in $\ccc_{m,n}$ over paths 
in $BConf^{0*}_{m-1}$ is the same as parallel transport in $\ss_{m-1,n} \cap 
\sl((m-1)n,\cc)^{\pi'},$ where $\pi'=(1^{m-2}(n-1)^{m-2}).$

\subsection {The Lagrangians} \label{sec:tl} Let $\tau =\{(\lambda_1, \lambda_2, \dots, 
\lambda_m), (\mu_1, \mu_2, \dots, \mu_m)\}$ be a point in the bipartite configuration space 
$BConf_m.$

\begin {definition}
A crossingless matching $\mat$ with endpoints $\tau$ is a collection of $m$ disjoint embedded 
arcs $(\delta_1, \dots, \delta_m)$ in $\cc$ such that each $\delta_k: [0,1] \to \cc$ satisfies 
$\delta_k(0) = \lambda_k$ and $\delta_k(1) = \mu_{\nu(k)},$ where $\mu$ is a permutation of 
$\{1,2, \dots, m\}.$ \end {definition}

To each crossingless matching $\mat$ we associate a Lagrangian submanifold $L_{\mat} \subset 
\ymnt,$ diffeomorphic to the product of $m$ copies of $\cpn,$ and unique up to Lagrangian 
isotopy. Since we had not fixed an ordering of $(\mu_1, \mu_2, \dots, \mu_m),$ we can assume 
that $\nu$ is the identity permutation, so that $\delta_k$ joins $\lambda_k$ and $\mu_k.$ We 
choose a path $[0,1) \to BConf_m$ that starts at $\tau$ and moves the endpoints of $\delta_1$ 
towards each other along that arc, while keeping all the other components of $\tau$ constant. 
The endpoints $\lambda_1$ and $\mu_{\nu(1)}$ should collide as $s \to 1.$ We assume that 
$\delta_1$ is a straight line near its midpoint $\lambda$, and that the colliding points move 
towards the midpoint with the same speed as $s \to 1.$ We can translate this whole picture into 
$BConf^0_m \subset \hh^{\pi}/W^{\pi}$ as before. We obtain 
a path $\gamma: [0,1] \to 
\hh^{\pi}/W^{\pi}$ with $\gamma(1)= \{(\lambda', \lambda'_2, \dots, \lambda'_m), (\lambda', 
\mu'_2, \dots, \mu'_m)\},$ such that $\lambda_k' = \lambda_k + \lambda' - \lambda$ and $\mu'_k 
= \mu_k - \lambda' - \lambda,$ for all $k=2, \dots, m.$

The construction of $L_{\mat}$ proceeds by induction on $m,$ with $n$ being kept fixed. In 
the case $m=1$ we have the vanishing $\cpn$ from Section~\ref{sec:van} in the fibers over 
$\gamma(1-s)$ for small $s.$ We use reverse (rescaled) parallel transport along $\gamma$ 
to move this back to the fiber $\yy_{1,n, \tau}$ over $\gamma(0).$ For the inductive step, 
denote by $\bar \mat$ the crossingless matching of $(2m-2)$ points which is obtained from 
$\mat$ after removing $\delta_1.$ By assumption, we have a Lagrangian $L_{\bar \mat} 
\subset \yy_{m-1, n, \bar \tau},$ where $\bar \tau$ consists of the endpoints of $\bar 
\mat.$ According to Lemma~\ref{identh}, the space $\yy_{m-1, n,\bar \tau}$ can be 
identified with a fiber of the singular set fibration $\ccc_{m,n} \to Conf^0_{\sigma}.$ 
Using parallel transport in $\ccc_{m,n},$ we can move the Lagrangian into the singular 
locus of $\yy_{m,n, \gamma(1)}.$ The local models from Lemma~\ref{thickthin} and 
Remark~\ref{remcom} tell us that 
we can apply the relative $\cpn$ construction from Section~\ref{sec:relvan} and get a 
Lagrangian submanifold in $\yy_{m,n, \gamma(1-s)}$ for small $s.$ Using reverse parallel 
transport along $\gamma$ we move this to the fiber over $\gamma(0).$ The result is the 
desired Lagrangian $L_{\mat} \subset \ymnt.$ Note that the necessary real analyticity 
condition on the K\"ahler metrics (cf. Definition~\ref{eff} and Section~\ref{sec:relvan}) 
is satisfied because of the way we chose the $\psi_k$'s in Section~\ref{sec:kh}.

Let us note a few properties of the Lagrangians that arise from this construction. Observe that 
the construction can be done in families as well. If $\mat(s)\  (s \in [0,1])$ is a 
smooth family of crossingless matchings with endpoints $\beta(s),$ one can construct 
Lagrangians $L_{\mat(s)} \subset \yy_{m,n, \beta(s)}$ depending smoothly on $s.$ Parallel 
transport along $\beta|_{[s,1]}$ can be used to carry them into a common fiber. We obtain a 
Lagrangian isotopy
\begin {equation}
\label {lagi}
L_{\mat(1)}\ \iso \ h_{\beta}^{\resc}(L_{\mat(0)}).
\end {equation}

One corollary is that a smooth family of crossingless matchings with fixed endpoints always 
yields a family of isotopic Lagrangians.

\begin {figure}
\begin {center}
\input {match.pstex_t}
\end {center}
\caption {Matchings $\mat$ and $\mat'$.}
\label {fig:match}
\end {figure}
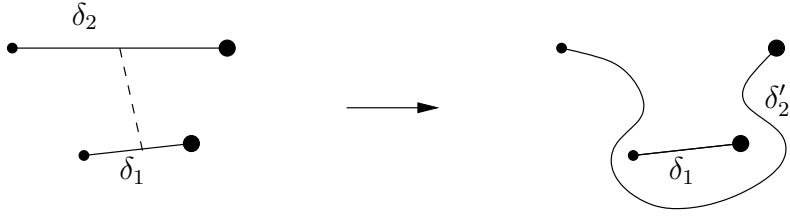

Note that in the first step of the construction of $L_{\mat}$ (corresponding to $m=1$), 
the matching is simply a path from $\lambda_1$ to $\delta_1,$ and therefore its isotopy 
class is unique. This implies that the procedure does not depend essentially on 
$\delta_1.$ A quick consequence is the following:

\begin {lemma}
\label {mm}
Let $\mat = (\delta_1, \delta_2, \delta_3, \dots, \delta_m)$ and $\mat'=(\delta_1, \delta'_2, 
\delta_3, \dots, \delta_m)$ be two crossingless matchings related to each other 
as shown in Figure~\ref{fig:match}: $\delta_2'$ is obtained from $\delta_2$ by sliding it 
over $\delta_1,$ using a dashed path (shown dashed in the figure) that does not intersect any 
of the other arcs. Then $L_{\mat}$ and $L_{\mat'}$ are Lagrangian isotopic.
\end {lemma}

Recall that in the local model given by Lemma~\ref{multithickthin}, we can change the order in 
an iteration of the vanishing $\cpn$ procedure, cf. Lemma~\ref{changeorder}. Applying this to 
our construction, we obtain:

\begin {lemma}
\label {indy}
Up to Lagrangian isotopy, $L_{\mat}$ is independent of the ordering of the components of the 
matching $\mat.$
\end {lemma}

\subsection {Floer cohomology} \label{sec:fc}

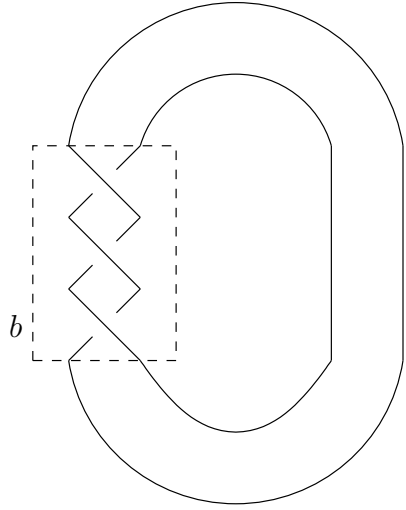
\begin {figure}
\begin {center}
\input {braid.pstex_t}
\end {center}
\caption {The trefoil as a braid closure.}
\label {fig:braid}
\end {figure}

Given an oriented link $\lk \subset S^3$ we can represent it as the closure of an $m$-stranded 
braid $b \in Br_m.$ (See Figure~\ref{fig:braid} for a presentation of the left-handed trefoil.)
Note that the braid group $Br_m$ is the fundamental group of the configuration space 
$Conf_{(1^m)}$ of $m$ unordered points in the plane.

\begin {figure}
\begin {center}
\input {std.pstex_t}
\end {center}
\caption {The standard crossingless matching $\mat_0$.}
\label {fig:std}
\end {figure}
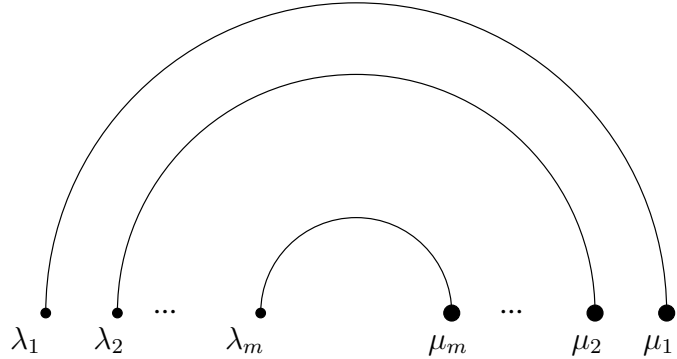

Consider the standard crossingless matching $\mat_0$ shown in Figure~\ref{fig:std}, where all 
the $\lambda$'s and the $\mu$'s are on the real line, with the ``thin points'' $\lambda$'s to 
the left of the ``thick points'' $\mu$'s, and all the arcs lying in the upper half-plane. Let 
$D \subset \cc$ be a disk that contains all the thin points and none of the thick ones. We 
denote by $U \subset BConf_m$ the open subset corresponding to bipartite configurations $\tau = 
\{(\lambda_1, \lambda_2, \dots, \lambda_m), (\mu_1,\mu_2, \dots, \mu_m)\},$ where the 
$\mu_k$'s are fixed to be the ones chosen for $\mat_0,$ while the $\lambda_k$'s are free to 
move inside the disk $D.$ Using a deformation retraction $\cc \to D,$ we can identify $Br_m$ 
with the fundamental group of $U.$

In this fashion, the braid $b$ representing $\lk$ induces a path $\beta: [0,1] \to U$ with 
$\beta(0) = \beta(1)$ being the configuration $\tau_0$ formed from the endpoints of $\mat_0.$ 
We consider the pair of Lagrangians $(L, L')$ in $M= \yy_{m,n, \tau_0},$ where $L$ is 
$L_{\mat_0}$ and $L'= h_{\beta}^{resc}(L)$ is obtained from $L$ by parallel transport. We 
denote by $w$ writhe of the braid $b,$ and set
\begin {equation}
\label {defff}
\krss(\lk) = HF^{*+(n-1)(m+w)}(L, L').
\end {equation}

The fact that these groups are link invariants will be shown in the next section. For now, 
let us check that we are in the setting of Section~\ref{sec:flo}, so that the Floer cohomology
groups on the right-hand side of (\ref{defff}) are well-defined. 

The complex manifold $M=\yy_{m,n, \tau_0}$ is an affine variety, and therefore Stein. The 
K\"ahler form $\omega = -dd^c \psi$ is exact by construction. The simultaneous resolution 
(\ref{slice}) shows that $M$ is deformation equivalent to the nilpotent fiber $\nnn_{\rho 
\pi}$ from (\ref{nnn}). It follows from the work of Maffei \cite{Maf} that $\nnn_{\rho \pi}$ 
is a quiver variety in the sense of Nakajima \cite{N1}, and hence hyper-K\"ahler. This 
implies that $c_1 = 0$ for $\nnn_{\rho \pi}$ and then the same must be true for $M.$ 

As noted in Section~\ref{sec:top}, the fiber $\nnn_{\rho \pi}$ is homotopy equivalent to the 
Spaltenstein variety $\Pi^{-1}(N_{m,n})$, whose rational cohomology was studied in 
\cite[Section 3.4]{BM}. Borho and MacPherson proved that, up to an even dimension shift, 
$H^*(\Pi^{-1}(N_{m,n}); \qq)$ is equal to a space of anti-invariants in the rational 
cohomology of the Springer variety of complete flags fixed by $N_{m,n}.$ Since the 
odd-dimensional Betti numbers of Springer varieties are zero (cf. \cite{CP}), the same must 
be true for those of Spaltenstein varieties. In particular, it follows that $H^1(M) =0.$

The Lagrangians $L$ and $L'$ are diffeomorphic to $(\cpn)^m$, so they satisfy $H_1 = 0$ as 
required. The Stiefel-Whitney class $w_2((\cpn)^m)$ is zero for $n$ even, while for $n$ 
odd it is the class $$(w,w,\dots, w) \in H^2((\cpn)^m; \zz/2\zz)\ \iso \ \bigoplus_{i=1}^m 
H^2(\cpn; \zz/2\zz),$$ where $w$ is the generator of $H^2(\cpn; \zz/2\zz) \ \iso \ 
\zz/2\zz.$

When $n$ is even, the Lagrangians $L$ and $L'$ are spin. The following lemma shows that the 
pair $(L, L')$ is relatively spin when $n$ is odd:

\begin {lemma}
\label {class}
Every element $Y \in \ss_{m,n} \cap \gg^{\pi, \reg}$ has exactly $m$ one-dimensional 
eigenspaces $E(\lambda_1), E(\lambda_2), \dots, E(\lambda_m),$ coresponding to eigenvalues 
$\lambda_j$ for
$j=1,\dots, m.$ Let $V$ be the complex line bundle over $\ss_{m,n} \cap \gg^{\pi, \reg}$ 
whose fibers are $E(\lambda_1) \otimes E(\lambda_2) \otimes \dots \otimes E(\lambda_m).$
Then the restriction of $w_2(V)$ to any Lagrangian $L_{\mat}$ coming from a matching $m$ is the 
class $(w,w, \dots, w) \in H^2((\cpn)^m; \zz/2\zz).$ 
\end {lemma}

\noindent \textbf{Proof.} Lemma~\ref{indy} says that the ordering of the factors is not 
essential in the construction of $L_{\mat}.$ Therefore, it suffices to show that the first 
component of $w_2(V)$ is $w.$ Since continuous deformations do not change the Stiefel-Whitney 
class, this is equivalent to showing that $w_2(V|_C)=w,$ where $C$ is a vanishing $\cpn$ 
appearing in the last step of the construction of $\cpn.$ More precisely, in the notation of 
Section~\ref{sec:tl}, we pick a point $Y_0$ in the singular locus of on $\yy_{m, n, 
\gamma(1)},$ and define $C$ to be the associated vanishing $\cpn$ in a nearby fiber $\yy_{m,n, 
\gamma(1-s)},$ arising from the local model in Lemma~\ref{thickthin}.

Recall that $\gamma(1)= \{(\lambda', \lambda'_2, \dots, \lambda'_m), (\lambda', \mu'_2, \dots, 
\mu'_m)\}.$ Let $N$ be a small neighborhood of $\gamma(1)$ in $\hh^{\pi}/W^{\pi}$ and $\tilde 
N$ its preimage in $\ss_{m,n} \cap \gg^{\pi}$ under (\ref{our}). Note that for $k=2, \dots, m,$ 
there are well-defined line bundles $E_k$ over $\tilde N$ whose fibers are the eigenspaces 
corresponding to the eigenvalue close to $\lambda_k'.$ There is also a line bundle $E_1$ over 
$\tilde N \cap \gg^{\pi, \reg}$ corresponding to the eigenspace with eigenvalue close to 
$\lambda',$ but $E_1$ cannot be extended to all of $\tilde N.$

Inside $\tilde N,$ the vanishing space $C$ is the boundary of a cone with vertex $x.$ Since the 
cone is contractible, it follows that all $E_k$ can be trivialized over $C,$ for $k=2, \dots, 
m.$ Consequently, $V|_C$ is isomorphic to $E_1|_C.$ Since we only care about the topology, we 
can move $C$ continuously into some $C' \ \iso \ \cpn$ that lies in the canonical semisimple 
slice $\ss^{ss}$ at $Y_0.$ It follows from the discussion in the proof of Lemma~\ref{thickthin} 
that $\ss^{ss} \cap \gg^{\pi}$ (after translation by $Y_0$) can be identified with $\cc^{2m-2} 
\times \slan.$ We can transport $C'$ further into some $C'' \ \iso \ \cpn \subset 
\{0\} \times \slan.$ The pullback of $E_1|_{C''}$ under this identification 
is the tautological bundle over some $\cpn \subset \slan,$ which has $w_2 = w. 
\hfill \fin$ 

\medskip

In order to define orientations on the moduli spaces of pseudoholomorphic curves, we need 
to choose the relative spin structure $st \in H^2(M; \zz/2\zz).$ We let $st$ be zero when 
$n$ is even, and $w_2(V)$ for $n$ odd, where $V$ is the bundle from Lemma~\ref{class} (or, 
more precisely, its restriction to $\yy_{m,n, \tau_0} \subset \ss_{m,n} \cap \gg^{\pi, 
\reg}).$

One other thing needed in the definition of Floer cohomology is the choice of orientations 
on the two Lagrangians. We make this inductively: when $m=1$ we give $\cpn \iso 
U(n)/(U(1) \times U(n-1))$ its natural complex structure and hence a complex orientation; 
then, at each step in the recursive definition of $L_{\mat}$ the new $\cpn$ factor can be 
identified with the one appearing in the $n=1$ case, and therefore we can also give it its 
complex structure and orientation. This defines an orientation on $L_{\mat}$ for any 
matching.

Finally, in order to have a well-defined absolute grading for the Floer groups, we need to 
endow the Lagrangians $L, L'$ with gradings as in Section~\ref{sec:abs}. Looking at the 
fibration (\ref{ourf}), we start by choosing arbitrary trivializations for the canonical 
bundles of the total space and the base. This produces a family of trivializations of the 
canonical bundles on the fibers $\ymnt.$ If we choose a grading for the Lagrangian 
$L_{\mat_0} \subset \yy_{m,n, \tau_0},$ we can then transport it in a continuous way to 
gradings on $h^{\resc}_{\beta|[0,s]}(L_{\mat_0}).$ In particular when $s=1$ we get a grading 
on the Lagrangian $L'.$ Relation (\ref{grad}) implies that the resulting $HF(L, L')$ does not 
depend on our choice of grading for $L.$

\section {Invariance under Markov moves} 
\label {sec:markov}

The Floer cohomology groups $\krss$ defined in
(\ref{defff}) are invariant under continuous deformations of the objects involved. We make 
the following observation:

\begin {lemma} All the noncanonical choices made in the definition of $\krss,$ with the 
exception of the braid element $b \in Br_m,$ are parametrized by weakly contractible 
spaces. \end {lemma} 

\noindent \textbf{Proof.} There were some choices made in the 
construction of the K\"ahler form $\Omega = -dd^c \psi:$ the constant $\alpha \in (m, 
\infty)$ and the functions $\psi_k: \cc \to \rr.$ The functions $\psi_k$ are required to 
satisfy certain properties, but these are preserved under linear interpolation. Therefore, 
the corresponding parameter space is a convex subset of $C^{\infty}(\cc, \rr),$ and hence 
contractible.

Another choice was made in Section~\ref{sec:fc}, where a particular path $\beta: [0,1] \to 
U$ was selected as a representative for a given homotopy class in $\pi_1(U)=Br_m.$ The 
space of these representatives is weakly contractible because $U$ is homotopy equivalent to 
the configuration space $\conf_m,$ and the latter was shown to be aspherical in \cite{FN}.

There are other choices involved in the definition of Floer cohomology, such as almost 
complex structures and perturbations. It is well-known that the respective parameter spaces 
are weakly contractible. $\hfill \fin$
\medskip

Thus, in order to prove that $\krss$ are link invariants, it suffices to show that they are 
independent of the presentation of the link $\lk$ as a braid closure. Two braids have the 
same closure if they are related by a sequence of Markov moves. The proof of 
Theorem~\ref{the} will be completed once we show invariance under the Markov moves.

\subsection {Markov $I$}
The type $I$ Markov move consists in replacing a braid $b$ by $s_k^{-1}bs_k,$ where $s_k$ 
is one of the standard generators $s_1, \dots, s_{m-1}$ of $Br_m$, the positive 
half-twist between the $k$th and $(k+1)$th strand. Choose a representative $\sigma_k$ for $s_k$ 
in the form of a loop in the subset $U \subset BConf_m$ from Section~\ref{sec:fc}. We also 
choose a representative $\sigma_{2m-k}$ for the half-twist between the thick points
$\mu_k$ and $\mu_{k+1},$ in the form of a loop in $BConf_m$ that fixes the thin points; in 
fact, we can assume that $\sigma_{2m-k}$ is the identity on the disk $D$ that contains the thin 
points. This implies that the parallel transport maps $h^{\resc}_{\sigma_{2m-k}}$ and 
$h^{\resc}_{\beta}$ do not interfere with each other and therefore commute.

The proof of Markov $I$ invariance is exactly as in \cite{SS}. Using (\ref{lagi}) and 
Lemma~\ref{mm}, we get that the the image of $L_{\mat_0}$ under parallel transport along 
$\sigma_{2m-k}^{-1} \circ \sigma_k$ is isotopic to $L_{\mat_0}.$ This fact, together with the
symplectomorphism invariance of Floer cohomology, leads to the string of identities
\begin{align*}
 & HF(L_{\mat_0},h^{\resc}_{\sigma_k^{-1}} h^{\resc}_\beta
 h^{\resc}_{\sigma_k}(L_{\mat_0})) 
\  \iso \ HF(h^{\resc}_{\sigma_k}(L_{\mat_0}),
 h^{\resc}_\beta h^{\resc}_{\sigma_k}(L_{\mat_0})) \\
 & \iso \ HF(h^{\resc}_{\sigma_{2m-k}}(L_{\mat_0}),
 h^{\resc}_\beta h^{\resc}_{\sigma_k}(L_{\mat_0})) 
\  \iso \ HF(L_{\mat_0},
 h^{\resc}_{\sigma_{2m-k}^{-1}} h^{\resc}_\beta
 h^{\resc}_{\sigma_k}(L_{\mat_0})) \\
 & \iso \ HF(L_{\mat_0},h^{\resc}_\beta
 h^{\resc}_{\sigma_{2m-k}^{-1} \circ \,\sigma_k}(L_{\mat_0})) 
 \ \iso \ HF(L_{\mat_0},h^{\resc}_\beta(L_{\mat_0})).
\end{align*} 

\subsection {Markov $II$} Given a braid $\bar b \in Br_{m-1},$ the type $II^+$ Markov 
move consists in adding a strand plus a half-twist of that strand with its neighbor,
so that the result is $b = s_{m-1}(\bar b \times 1) \in Br_m.$ There is also a type $II^-$ 
Markov move where instead of $s_{m-1}$ we have the negative half-twist $s_{m-1}^{-1},$ but 
since its treatment is completely similar to that for Markov $II^+$, we will just focus on 
$II^+.$ 

Consider the standard configuration $\tau_0 = (\lambda_1, \dots, \lambda_m, \mu_m, \dots, 
\mu_1)$ as in Figure~\ref{fig:std}. The eigenvalues $\lambda_{m-1}$ and $\lambda_m$ are 
those exchanged by $s_{m-1}.$ If $m \geq 3,$ we can assume that $\lambda_{m-1}, 
\lambda_m,$ and $\mu_m$ are small and 
\begin {equation}
\label {ser}
\lambda_{m-1} + \lambda_m + (n-1) \mu_m =0.
\end {equation}

In other words, $(\lambda_{m-1}, \lambda_m)$ is a point in a small bidisk $P$ around a 
configuration where two thin and one thick eigenvalues coincide, as in Lemma~\ref{ttt}. 
Let $\bar \tau_0$ be the configuration $(\lambda_1, \dots, \lambda_{m-2}, 0, \mu_{m-1}, 
\mu_{m-2}, \dots, 0)$ and $\bar \beta$ a loop in $BConf^0_{m-1}$ based at $\bar \tau_0$ 
and representing the braid $\bar b.$ Consider the Lagrangians
$$ K = L_{\bar \mat}, \ \ K' = h_{\bar \beta}^{\resc}(L_{\bar \mat}) $$
in $\yy_{m-1,n, \bar \tau_0}.$ The local model described in Lemma~\ref{ttt} and 
Remark~\ref{remcom2} allows us to apply the discussion in Section~\ref{sec:fibcase} and 
construct relative vanishing projective spaces $L_{d,z}, L'_{d,z}$ in the fiber $\yy_{m,n, 
\tau_0}.$ Here $d = \mu_m$ and $z = (\mu_m - \lambda_{m-1})(\mu_m - \lambda_m)$ as in 
Section~\ref{sec:non}.  

Note that $L=L_{\mat_0}$ is isotopic to $L_{d,z}$ by construction. The second Lagrangian 
$L'= h_{\beta}^{\resc}(L) $ is obtained from this by parallel transport. A deformation 
argument as in \cite[Section 5(D)]{SS} shows that $L'$ is isotopic to the Lagrangian 
obtained by taking $K',$ doing the relative vanishing $\cpn$ construction to get 
$L'_{d,z},$ and then applying the twist $s_{n-1}.$ This last step corresponds to going 
around the loop $\gamma_{d,z}$ from Figure~\ref{fig:path}. This new Lagrangian is exactly 
$L''_{d,z}$ from (\ref{zwei}).

At this point we can apply Lemma~\ref{thom} and find that $HF(L, L') \iso HF(K, K').$ In 
other words, the Floer groups associated to the braid $\bar b$ and its Markov $II^+$ 
transform $b$ are isomorphic. The shift by $(n-1)(m+w)$ in (\ref{defff}) was 
chosen so that the absolute gradings match as well. (See \cite[Section 6(A)]{SS} for more 
details in the $n=2$ case.) 

We should note that the discussion above was only applicable for $m \geq 3.$ When $m=2,$ 
the relation (\ref{ser}) cannot be arranged. However, we can bring $\lambda_{m-1}, 
\lambda_m$ and $\mu_m$ close to a nonzero value instead, and then use suitable parallel 
transport in $\ccc_{m,n}$ to bring that value back to zero.

\begin {remark}
Instead of doing the twist $s_{n-1}$ in Markov $II^+$ we could just compare the situations 
for the braids $\bar b$ and $\bar b \times 1.$  The Floer cohomology 
of $L= L_{d,z}$ and $L'_{d,z}$ can be computed using Lemma~\ref{kunneth}. Let $O$ 
denote the unknot. Taking into 
account the absolute gradings, we get that our link invariants satisfy: 
$$ \krss ( \lk  \amalg O)= \krss (\lk) \otimes H^{*+n-1}(\cpn).$$

In particular, the invariant of the unlink of $p$ components is the tensor product of $p$ 
copies of $H^{*+n-1}(\cpn).$
\end {remark}

\section {Calculation for the trefoil}
\label {sec:trefoil}

In this section we prove Proposition~\ref{3f}. We work with left-handed trefoil knot 
$\lk,$ which is the closure of the braid $b=s_1^3 \in Br_2.$ We need to consider Floer 
cohomology inside the space $Y=\yy_{2,n,\tau_0},$ where $\tau_0$ is the standard 
configuration as in Figure~\ref{fig:std}.

Just as in the proof of Markov $II$ invariance, we bring the first three eigenvalues 
$\lambda_1, \lambda_2, \mu_2$ close together. This allows us to work inside the setting of 
Section~\ref{sec:fibcase} once again. For the base space $X$ we have a fiber of the map 
$\chin: \slan \to \cc,$ which contains a vanishing projective space 
$K \iso \cpn.$ This corresponds to the first step in the construction of the Lagrangian 
$L_{\mat},$ i.e. to parallel transport by bringing $\lambda_1$ and $\mu_1$ together. Over 
$X$ we have three bundles $\ff_1 \subset \ff_2 \subset \ee.$ Here the fiber of $\ff_1$ 
over $A \in X$ is the eigenspace for the thin eigenvalue of $A,$ while $\ff_2$ is the whole 
space on which $A$ acts. (Hence, in fact $\ff_2$ is a trivial bundle.) We can choose a 
compatible holomorphic quadruple bundle $\tt \to X$ as in Section~\ref{sec:ttt}, and then 
get a fiber bundle $\xx(\tt) \subset \sl(\ee) \to X.$ This describes the local structure of 
$Y$ near $X.$ Locally near $X$ the K\"ahler form can be deformed to be one obtained by 
combining the standard form on $X$ (from Section~\ref{sec:van}) with a Hermitian form on 
$\ee.$ Recall that if $X$ is equipped with the standard form then $K$ can be described 
explicitly by (\ref{ua}).

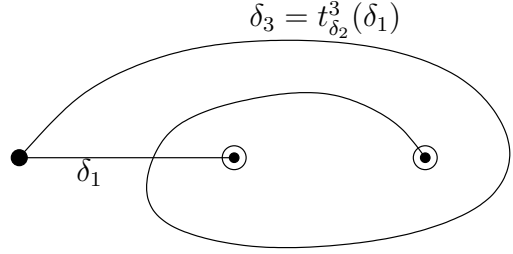
\begin {figure}
\begin {center}
\input {t3.pstex_t}
\end {center}
\caption {Projections of the two Lagrangians.}
\label {fig:t3}
\end {figure}

Applying the same reasoning as in the proof of Markov $II,$ we find that, up to 
isotopy, the Lagrangian $L = L_{\mat_0} \subset Y$ is obtained from $K$ by taking 
fiberwise over it the space $\Lambda_{\delta_1}$ from (\ref{bigla}). Also, the second 
Lagrangian $L' = h_{\beta}^{\resc}(L)$ is obtained from $K$ in the same way, but using 
fibers $\Lambda_{\delta_3}$ instead. Here $\delta_3 = t^3_{\delta_2}(\delta_1)$ as in 
Figure~\ref{fig:t3}.

The paths $\delta_1$ and $\delta_3$ intersect in one endpoint and one interior point. By 
the discussion in Section~\ref{sec:ac}, the corresponding intersections of 
$\Lambda_{\delta_1}$ and $\Lambda_{\delta_2}$ are a point and a $(2n-3)$-dimensional 
sphere, respectively. This implies that $L$ and $L'$ intersect in the disjoint union of a 
$\cpn$ and a $S^{2n-3}$-bundle over $\cpn.$ Using Remark~\ref{klm}, the latter turns out 
to be exactly the unit tangent bundle of $\cpn.$ Indeed, the restriction of $\ff_1$ to $K 
\iso \cpn$ is the tautological complex line bundle $\xi,$ and we have $T\cpn = 
\text{Hom}(\xi, \cc^n/\xi).$

The intersection $L \cap L'$ is clean in the sense of Pozniak \cite{P}. This gives a 
Morse-Bott long exact sequence
$$ \cdots \to H^{*-2}(UT\cpn) \to HF^*(L, L') \to H^*(\cpn) \to H^{*-1}(UT\cpn) \to 
\cdots $$

The absolute gradings can be understood from the geometry of the fibers. Observe that the 
differential $H^*(\cpn) \to H^{*-1}(UT\cpn) $ is zero for simple algebraic reasons. After 
taking into account the shift factor $(n-1)(m+w) = -n+1$ in (\ref{defff}), the proof of 
Proposition~\ref{3f} is completed. Explicitly, we have 
$$ \krsk(\lk)=
\begin{cases}
\zz & \text{ for $k=n-1$;}\\
\zz^2 & \text{ for $k=n+1+2j, \ 0 \leq j \leq n-2$;}\\
\zz/n\zz & \text{ for $k=3n-1$;}\\
\zz & \text{ for $k=3n+2j, \ 0 \leq j \leq n-2$;}\\
0 & \text{ otherwise.}
\end{cases}
$$

\section {Other link invariants}
\label {sec:other}

This section is more speculative in nature. We suggest some further avenues of research, 
by explaining how variants of our construction could lead to other link invariants.

\subsection {Other Lie algebras and representations} In Section~\ref{sec:motiv} we noted 
that the link polynomial $P_{(n)}$ from (\ref{skey}) is associated to the standard 
$n$-dimensional representation $V$ of $\sl(n, \cc).$ The procedure described there can in 
fact be applied to any complex simple Lie algebra $\ggg,$ the result being a polynomial 
invariant of decorated links, where the link components are decorated by 
finite-dimensional irreducible representations of $\ggg$ (\cite{RT}). In the case of 
$P_{(n)},$ all components are decorated by $V.$

Let us consider the situation in which the Lie algebra is still $\sl(n, \cc),$ but the 
components of the link are decorated by various exterior products of the standard 
representation $V.$ Represent the link as the closure of a braid $b$, and suppose that 
$b$ has $m_i$ strands colored in the representation $\Lambda^kV,$ for $0 < k < n.$ Let 
the total number of strands be $m = m_1 + \dots + m_{n-1}.$ Form the representation $$ W = 
\bigotimes_{k=1}^{n-1} \ (\Lambda^kV)^{\otimes m_k}.$$

To the braid $b$ one associates a map $F_b(q):W \to W,$ and then the  
polynomial invariant is obtained as the image of $1$ under a map
$$\begin {CD}
\cc \to  W \otimes W^* @>{F_b(q) \times id}>> W \otimes W^* \to \cc. \end {CD} $$
similar to (\ref{poly}). Note that the dual of $\Lambda^k V$ is $\Lambda^{n-k}V,$ so that
$$ W \otimes W^* = \bigotimes_{k=1}^{n-1} \ (\Lambda^kV)^{\otimes (m_k + m_{n-k})}.$$

In \cite[Section 11]{KR1}, Khovanov and Rozansky sketched the construction of a 
homology theory that categorifies this particular polynomial invariant. We conjecture that 
a similar theory can be constructed from Floer homology, and that it is equivalent to
that of Khovanov-Rozansky. One should still use intersections of transverse slices and 
adjoint orbits in $\sl(nm, \cc),$ but the partition $\pi$ from Section~\ref{sec:av} should 
be chosen as
$$ \pi = \bigl ( 1^{m_1 + m_{n-1}} 2^{m_2 + m_{n-2}} \dots (n-1)^{m_{n-1} + m_1}  \bigr ). 
$$ 

The second partition $\rho$ should still be $(n^m).$ To define the two Lagrangians one 
would need a generalization of the vanishing $\cpn$ construction, in the form of vanishing 
Grassmannians $Gr(n,k),$ for $0 < k < n.$ The resulting Lagrangians should be 
diffeomorphic to products of $m_k$ copies of the $Gr(n, k),$ over all $k.$ We expect that 
many of the arguments in this paper would go through in this more general setting as well. 
However, since the local models near degenerations are more complicated, it is possible 
that additional technical difficulties could appear, especially in the proof of Markov $II$ 
invariance.

More categorifications of quantum polynomial invariants were constructed by Gukov and 
Walcher in \cite{GW}. They mainly work with representations of the Lie algebras $\so(n, 
\cc),$ but also discuss representations of some exceptional Lie algebras. We expect that 
one can define similar link homology invariants using Floer theory. One should use 
transverse slices and adjoint orbits in a suitable Lie algebra, instead of those in 
$\sl(nm, \cc).$

\subsection {An involution} In some unpublished work (\cite{SS2}), Seidel and Smith define 
an involution on their manifold $\yy_{m,2, \tau},$ and do Floer theory in the 
fixed point set of that involution. (See also \cite[Section 7]{M}.) They show that the 
resulting Floer homology groups are isomorphic to $\widehat{HF} (\Sigma(\lk) \# S^1 
\times S^2),$ where $\widehat{HF}$ stands for a variant of Heegaard Floer homology 
(\cite{OS0}), and $\Sigma(\lk)$ for the double cover of $S^3$ branched over the link 
$\lk.$

There is a similar involution $\iota$ in our theory, for any $n \geq 2.$ On the slice 
(\ref{smn}), it takes each block $Y_{k1}$ into its transpose $Y_{k1}^T.$ A linear algebra 
exercise shows that $\iota$ does not change the conjugation class of a matrix, and 
therefore descends to an involution on each $\ymnt.$ By doing all constructions in this 
paper $\zz/2\zz$-equivariantly, we can arrange so that the involution also acts on the 
Lagrangians $L$ and $L'.$ We can then try to apply Floer's theory to the fixed point sets 
of $\iota|_L, \iota|_{L'}$, viewed inside the fixed point set of $\iota|_{\ymnt}.$ We 
expect the resulting Floer cohomology groups to form a series of link invariants $\kr_{\nn 
\iota}^* (\lk),$ for $n \geq 2,$ with the property that $\kr_{\nn \iota}^* (\text{unknot}; 
\zz/2\zz) = H^{*+(n-1)/2}(\mathbb{RP}^{n-1}; \zz/2\zz).$

\end{document}

%% file: path.pstex_t
\begin{picture}(0,0)%
\includegraphics{path.pstex}%
\end{picture}%
\setlength{\unitlength}{3947sp}%
\begingroup\makeatletter\ifx\SetFigFont\undefined%
\gdef\SetFigFont#1#2#3#4#5{%
  \reset@font\fontsize{#1}{#2pt}%
  \fontfamily{#3}\fontseries{#4}\fontshape{#5}%
  \selectfont}%
\fi\endgroup%
\begin{picture}(4632,3014)(2011,-3668)
\put(5251,-2311){\makebox(0,0)[lb]{\smash{{\SetFigFont{12}{14.4}{\rmdefault}{\mddefault}{\updefault}{\color[rgb]{0,0,0}$\zeta_d$}%
}}}}
\put(2401,-2011){\makebox(0,0)[lb]{\smash{{\SetFigFont{12}{14.4}{\rmdefault}{\mddefault}{\updefault}{\color[rgb]{0,0,0}$z$}%
}}}}
\put(2026,-2386){\makebox(0,0)[lb]{\smash{{\SetFigFont{12}{14.4}{\rmdefault}{\mddefault}{\updefault}{\color[rgb]{0,0,0}$0$}%
}}}}
\put(3826,-3511){\makebox(0,0)[lb]{\smash{{\SetFigFont{12}{14.4}{\rmdefault}{\mddefault}{\updefault}{\color[rgb]{0,0,0}$\gamma_{d,z}$}%
}}}}
\end{picture}%

%% file: lag.pstex_t
\begin{picture}(0,0)%
\includegraphics{lag.pstex}%
\end{picture}%
\setlength{\unitlength}{3947sp}%
\begingroup\makeatletter\ifx\SetFigFont\undefined%
\gdef\SetFigFont#1#2#3#4#5{%
  \reset@font\fontsize{#1}{#2pt}%
  \fontfamily{#3}\fontseries{#4}\fontshape{#5}%
  \selectfont}%
\fi\endgroup%
\begin{picture}(3698,1065)(961,-2155)
\put(3601,-2086){\makebox(0,0)[rb]{\smash{{\SetFigFont{12}{14.4}{\rmdefault}{\mddefault}{\updefault}{\color[rgb]{0,0,0}$\delta_2$}%
}}}}
\put(3526,-1261){\makebox(0,0)[rb]{\smash{{\SetFigFont{12}{14.4}{\rmdefault}{\mddefault}{\updefault}{\color[rgb]{0,0,0}$t_{\delta_2}(\delta_1)$}%
}}}}
\put(976,-1786){\makebox(0,0)[lb]{\smash{{\SetFigFont{12}{14.4}{\rmdefault}{\mddefault}{\updefault}{\color[rgb]{0,0,0}$\alpha_1$}%
}}}}
\put(1801,-2086){\makebox(0,0)[rb]{\smash{{\SetFigFont{12}{14.4}{\rmdefault}{\mddefault}{\updefault}{\color[rgb]{0,0,0}$\delta_1$}%
}}}}
\put(4501,-1711){\makebox(0,0)[lb]{\smash{{\SetFigFont{12}{14.4}{\rmdefault}{\mddefault}{\updefault}{\color[rgb]{0,0,0}$\alpha_3$}%
}}}}
\put(2101,-1711){\makebox(0,0)[lb]{\smash{{\SetFigFont{12}{14.4}{\rmdefault}{\mddefault}{\updefault}{\color[rgb]{0,0,0}$\alpha_2$}%
}}}}
\end{picture}%

%% file: match.pstex_t
\begin{picture}(0,0)%
\includegraphics{match.pstex}%
\end{picture}%
\setlength{\unitlength}{3947sp}%
\begingroup\makeatletter\ifx\SetFigFont\undefined%
\gdef\SetFigFont#1#2#3#4#5{%
  \reset@font\fontsize{#1}{#2pt}%
  \fontfamily{#3}\fontseries{#4}\fontshape{#5}%
  \selectfont}%
\fi\endgroup%
\begin{picture}(4909,1330)(2063,-3632)
\put(2476,-2461){\makebox(0,0)[lb]{\smash{{\SetFigFont{12}{14.4}{\rmdefault}{\mddefault}{\updefault}{\color[rgb]{0,0,0}$\delta_2$}%
}}}}
\put(2776,-3436){\makebox(0,0)[lb]{\smash{{\SetFigFont{12}{14.4}{\rmdefault}{\mddefault}{\updefault}{\color[rgb]{0,0,0}$\delta_1$}%
}}}}
\put(6226,-3436){\makebox(0,0)[lb]{\smash{{\SetFigFont{12}{14.4}{\rmdefault}{\mddefault}{\updefault}{\color[rgb]{0,0,0}$\delta_1$}%
}}}}
\put(6826,-2986){\makebox(0,0)[lb]{\smash{{\SetFigFont{12}{14.4}{\rmdefault}{\mddefault}{\updefault}{\color[rgb]{0,0,0}$\delta_2'$}%
}}}}
\end{picture}%

%% file: braid.pstex_t
\begin{picture}(0,0)%
\includegraphics{braid.pstex}%
\end{picture}%
\setlength{\unitlength}{3947sp}%
\begingroup\makeatletter\ifx\SetFigFont\undefined%
\gdef\SetFigFont#1#2#3#4#5{%
  \reset@font\fontsize{#1}{#2pt}%
  \fontfamily{#3}\fontseries{#4}\fontshape{#5}%
  \selectfont}%
\fi\endgroup%
\begin{picture}(2502,3168)(811,-3820)
\put(826,-2761){\makebox(0,0)[lb]{\smash{{\SetFigFont{12}{14.4}{\rmdefault}{\mddefault}{\updefault}{\color[rgb]{0,0,0}$b$}%
}}}}
\end{picture}%

%% file: std.pstex_t
\begin{picture}(0,0)%
\includegraphics{std.pstex}%
\end{picture}%
\setlength{\unitlength}{3947sp}%
\begingroup\makeatletter\ifx\SetFigFont\undefined%
\gdef\SetFigFont#1#2#3#4#5{%
  \reset@font\fontsize{#1}{#2pt}%
  \fontfamily{#3}\fontseries{#4}\fontshape{#5}%
  \selectfont}%
\fi\endgroup%
\begin{picture}(4198,2253)(2911,-3055)
\put(4276,-2986){\makebox(0,0)[lb]{\smash{{\SetFigFont{12}{14.4}{\rmdefault}{\mddefault}{\updefault}{\color[rgb]{0,0,0}$\lambda_m$}%
}}}}
\put(3451,-2986){\makebox(0,0)[lb]{\smash{{\SetFigFont{12}{14.4}{\rmdefault}{\mddefault}{\updefault}{\color[rgb]{0,0,0}$\lambda_2$}%
}}}}
\put(2926,-2986){\makebox(0,0)[lb]{\smash{{\SetFigFont{12}{14.4}{\rmdefault}{\mddefault}{\updefault}{\color[rgb]{0,0,0}$\lambda_1$}%
}}}}
\put(6901,-2986){\makebox(0,0)[lb]{\smash{{\SetFigFont{12}{14.4}{\rmdefault}{\mddefault}{\updefault}{\color[rgb]{0,0,0}$\mu_1$}%
}}}}
\put(6451,-2986){\makebox(0,0)[lb]{\smash{{\SetFigFont{12}{14.4}{\rmdefault}{\mddefault}{\updefault}{\color[rgb]{0,0,0}$\mu_2$}%
}}}}
\put(5551,-2986){\makebox(0,0)[lb]{\smash{{\SetFigFont{12}{14.4}{\rmdefault}{\mddefault}{\updefault}{\color[rgb]{0,0,0}$\mu_m$}%
}}}}
\end{picture}%

%% file: t3.pstex_t
\begin{picture}(0,0)%
\includegraphics{t3.pstex}%
\end{picture}%
\setlength{\unitlength}{3947sp}%
\begingroup\makeatletter\ifx\SetFigFont\undefined%
\gdef\SetFigFont#1#2#3#4#5{%
  \reset@font\fontsize{#1}{#2pt}%
  \fontfamily{#3}\fontseries{#4}\fontshape{#5}%
  \selectfont}%
\fi\endgroup%
\begin{picture}(3151,1573)(1293,-2738)
\put(1876,-2311){\makebox(0,0)[rb]{\smash{{\SetFigFont{12}{14.4}{\rmdefault}{\mddefault}{\updefault}{\color[rgb]{0,0,0}$\delta_1$}%
}}}}
\put(3751,-1336){\makebox(0,0)[rb]{\smash{{\SetFigFont{12}{14.4}{\rmdefault}{\mddefault}{\updefault}{\color[rgb]{0,0,0}$\delta_3=t_{\delta_2}^3(\delta_1)$}%
}}}}
\end{picture}%